\documentclass[10pt]{article}

\usepackage[margin=1.08in]{geometry}
\usepackage{amsmath,amsthm,amssymb,mathtools}
\usepackage{bm}
\usepackage{booktabs}
\usepackage{array}
\usepackage{multirow}
\usepackage{graphicx}
\usepackage{float}
\usepackage{caption}
\usepackage{subcaption}
\usepackage{algorithm}
\usepackage{algpseudocode}
\usepackage{enumitem}
\usepackage{microtype}
\usepackage[dvipsnames]{xcolor}
\usepackage[hypertexnames=false]{hyperref}
\usepackage[nameinlink,capitalise,noabbrev]{cleveref}
\usepackage{stmaryrd}
\usepackage{placeins}
\graphicspath{{Figure/}}

\hypersetup{
  colorlinks=true,
  linkcolor=MidnightBlue,
  citecolor=ForestGreen,
  urlcolor=BrickRed
}
\allowdisplaybreaks
\numberwithin{equation}{section}

\newtheorem{theorem}{Theorem}[section]
\newtheorem{lemma}[theorem]{Lemma}
\newtheorem{proposition}[theorem]{Proposition}

\newtheorem{definition}[theorem]{Definition}
\newtheorem{remark}[theorem]{Remark}

\newcommand{\R}{\mathbb R}
\newcommand{\E}{\mathbb E}
\newcommand{\dd}{\,\mathrm d}
\newcommand{\Th}{\mathcal T_h}
\newcommand{\Eh}{\mathcal E_h}
\newcommand{\Ehi}{\mathcal E_h^{\mathrm i}}
\newcommand{\Ehb}{\mathcal E_h^{\partial}}
\newcommand{\Pmat}{P}
\newcommand{\avg}[1]{\left\{#1\right\}}
\newcommand{\jump}[1]{\left[\!\left[#1\right]\!\right]}
\newcommand{\norml}[1]{\left\|#1\right\|_{\ell^2}}
\newcommand{\transpose}{^{\mathsf T}}
\newcommand{\hhat}{\widehat h}
\newcommand{\Bhat}{\widehat B}
\newcommand{\qxhat}{\widehat q^{x}}
\newcommand{\qyhat}{\widehat q^{y}}
\newcommand{\uhat}{\widehat u}
\newcommand{\vhat}{\widehat v}
\newcommand{\hD}{\widehat h_{\Delta}}
\newcommand{\BD}{\widehat B_{\Delta}}

\newcommand{\qxD}{\widehat q^{x}_{\Delta}}
\newcommand{\qyD}{\widehat q^{y}_{\Delta}}
\newcommand{\uD}{\widehat u_{\Delta}}
\newcommand{\vD}{\widehat v_{\Delta}}
\newcommand{\fq}{\widehat f_q}
\newcommand{\fqu}{\widehat f_{qu}}
\newcommand{\fqv}{\widehat f_{qv}}
\newcommand{\kpce}{K_{\mathrm{PCE}}}
\newcommand{\Std}{\operatorname{Std}}
\newcommand{\Var}{\operatorname{Var}}
\newcommand{\Pcor}{\Pmat^{\mathrm{cor}}}
\newcommand{\hK}{\Delta_K}
\newcommand{\err}{\mathrm{err}}
\newcommand{\Ccfl}{C_{\mathrm{CFL}}}
\newcommand{\Mtvb}{M_{\mathrm{TVB}}}
\newcommand{\Vk}{\mathbb V_\Delta^k}
\newcommand{\Nxi}{N_\xi}

\begin{document}

\title{A Structure-Preserving Discontinuous Galerkin--Stochastic Galerkin
Method for the Two-Dimensional Shallow Water Equations with Random Inputs}

\author{
Yekaterina Epshteyn\thanks{Department of Mathematics, University of Utah,
155 S 1400 E, Salt Lake City, UT 84112, USA
(\href{mailto:epshteyn@math.utah.edu}{epshteyn@math.utah.edu}).}
\and
Akil Narayan\thanks{Department of Mathematics and Scientific Computing and
Imaging Institute, University of Utah, 72 S Central Campus Dr, Salt Lake
City, UT 84112, USA (\href{mailto:akil@sci.utah.edu}{akil@sci.utah.edu}).}
\and
Yinqian Yu\thanks{Department of Mathematics, University of Utah,
155 S 1400 E, Salt Lake City, UT 84112, USA
(\href{mailto:yinqian.yu@utah.edu}{yinqian.yu@utah.edu}).}
}
\date{}
\maketitle


\begin{abstract}
The two-dimensional shallow water equations are a fundamental model for
geophysical free-surface flows, and in many applications their inputs, in
particular the bottom topography and the initial state, are only partially
known. To address this uncertainty we develop a high-order discontinuous Galerkin--stochastic Galerkin
(DG--SG) method for the two-dimensional shallow water equations with random
inputs. Building on the hyperbolicity-preserving stochastic Galerkin (SG)
formulation of [Dai, Epshteyn, Narayan, J.~Comput.~Phys.~2022] and on the
skew-symmetric discontinuous Galerkin framework of [Fu, J.~Sci.~Comput.~2022],
we derive a two-dimensional skew-symmetric SG system that is hyperbolic under
the same positive-definiteness condition as the conservative formulation and
that admits a companion entropy balance law. The resulting semi-discrete
DG--SG scheme on triangular meshes is locally conservative for the water
height and the discharge, well-balanced for stochastic lake-at-rest states, and entropy conservative or entropy stable
depending on the choice of interface dissipation. The proof of these properties rests on the
symmetry of the SG product operator and on an elementwise projection of both
velocity components. A fully discrete algorithm is obtained by combining
strong-stability-preserving Runge--Kutta time stepping with a
hyperbolicity-preserving time step restriction, a desingularized velocity
projection, and a characteristic total-variation-bounded limiter on
equilibrium variables. To monitor the entropy stability on domains with
outflow boundaries, we formulate the augmented energy of [Epshteyn, Narayan,
Yu, Comput.~Methods Appl.~Mech.~Engrg.~2025] for the DG weak form; it
accounts for the energy flux through the boundary and thereby separates the
dissipation of the scheme from the physical exchange of energy with the
exterior. Numerical experiments, including two convergence studies,
stochastic lake-at-rest perturbations, a nearly dry submerged plateau,
circular and curved dam breaks with uncertain water height, and problems
with two and with ten independent random variables, confirm the designed order of
accuracy, the well-balanced property, and the monotone decay of the
augmented energy.
\end{abstract}

\noindent\textbf{Keywords:}
stochastic shallow water equations; stochastic Galerkin method;
discontinuous Galerkin method; well-balanced method; entropy stability;
hyperbolicity preservation.

\medskip
\noindent\textbf{2020 Mathematics Subject Classification:}
35L65; 35Q35; 35R60; 65C20; 65M60; 65M70.

\section{Introduction}\label{sec:introduction}

The Saint--Venant system of shallow water equations (SWE) is a classical and
widely used model for geophysical free-surface flows in regimes where the
horizontal length scale is much larger than the vertical one, such as rivers,
lakes, estuaries, and coastal regions \cite{de1871theorie,kurganov2018finite,XING2017361}.
Despite its broad applicability, the deterministic model is often insufficient
in realistic settings: the bottom topography is measured with limited
resolution, and the initial water surface and discharge are frequently
inferred from sparse data. Incorporating such uncertainty is essential for
improving predictive capability, and it motivates the study of the SWE with
random (parametric) inputs, in which the uncertain inputs are described
through random variables; we refer to this model as the stochastic SWE.

We employ polynomial chaos expansions (PCE) to represent uncertainty.
Introduced by Wiener \cite{wiener1938homogeneous} and subsequently
generalized to a broad family of probability measures
\cite{xiu2002wiener,ghanem1991stochastic,le2010spectral,xiu2010numerical},
PCE encodes a random field through a truncated spectral expansion in the
random variables. Non-intrusive methods rely on an ensemble of realizations
and therefore reuse existing deterministic solvers: Monte Carlo and
multilevel Monte Carlo methods \cite{mishra2012multi,mishra2012multilevel}
estimate statistics of the solution directly from random samples, whereas
stochastic collocation methods \cite{xiu2005high,nobile2008sparse} construct
an interpolant, or the expansion coefficients, from solutions at prescribed
nodes in the random space. In this paper we instead
consider the intrusive stochastic Galerkin (SG) approach
\cite{ghanem1991stochastic,xiu2010numerical,sullivan2015introduction}, which
projects the residual of the governing equations onto the stochastic
polynomial space and yields a single coupled deterministic system for the
expansion coefficients.
The intrusive approach is chosen for two reasons.
First, a surrogate assembled from collocation solutions at different nodes
does not, in general, inherit the structural properties of the individual
deterministic solutions, and there is no mathematical guarantee that it
does, whereas the SG coefficients solve one coupled system for which such
properties can be proved. Second, for solutions that depend smoothly on the
random inputs the truncated PCE converges spectrally, so that a few modes
suffice, whereas sampling-based methods require a deterministic solve for
every sample or collocation node, and their number grows rapidly with the
accuracy required and with the number of random variables.
For problems whose solutions depend smoothly on the
random inputs, such as elliptic equations with random coefficients, SG
methods have a well-developed approximation theory
\cite{babuska2004galerkin,ernst2012convergence}, and they have been applied
to a variety of problems \cite{xiu2009efficient,eigel2014adaptive,hu2016stochastic,shu2017stochastic,tryoen2010intrusive}.
For nonlinear hyperbolic systems, however, the SG projection of the flux
does not in general preserve hyperbolicity, and the SG system may lose
well-posedness \cite{despres2013robust,poette2009uncertainty,gerster2019hyperbolic,jin2019study}.
A variety of strategies have been proposed to address this issue, including
entropy-variable and Roe-variable formulations
\cite{poette2009uncertainty,despres2013robust,gerster2020entropies,pettersson2014stochastic},
formulations based on a symmetrized form of the system
\cite{wu2017stochastic},
operator splitting \cite{chertock2015well,chertock2015operator},
filtering and hyperbolicity-preserving limiters
\cite{kusch2020filtered,schlachter2018hyperbolicity,durrwachter2020hyperbolicity},
kinetic and Haar-type formulations \cite{jin2019study,gerster2025haar},
cross-mode stabilized stochastic finite element formulations
\cite{chen2023cross}, and well-balanced and entropy stable SG
discretizations
\cite{jin2016well,zhong2022entropy,bender2024entropy,offner2026high}.
In the linear setting, SG approximations of a class of linear conservation
laws were studied in \cite{pulch2012generalised}, and bi-orthogonal
polynomial bases that decouple the SG equations for transport problems with
random velocity were proposed in \cite{zhou2012galerkin}.

This manuscript builds on prior work specifically for the SWE system. In particular we leverage a hyperbolicity-preserving SG formulation that was introduced in
\cite{dai2021hyperbolicity} in one spatial dimension and extended to two
dimensions in \cite{dai2022hyperbolicity}. Its key ingredient is a specific
Galerkin treatment of the non-polynomial flux components $(q^x)^2/h$,
$(q^y)^2/h$ and $q^xq^y/h$
through an SG product matrix $\Pmat(\hhat)$ built from the PCE coefficients
of the water height; the SG system is hyperbolic whenever $\Pmat(\hhat)$ is
positive definite. Energy conservative and energy stable finite-volume
schemes for this formulation were developed in \cite{dai2024energy} in one
dimension and in \cite{epshteyn2025energy} in two dimensions, together with
well-balanced and hyperbolicity-preserving reconstructions. These schemes
are second-order accurate. A detailed review and comparison of these
schemes, as well as new developments in one-dimensional settings, can be found in
\cite{epshteyn2026}.

In the deterministic setting, high-order
discontinuous Galerkin (DG) methods
\cite{cockburn1998runge,cockburn2001runge,hesthaven2008nodal} combine the
local conservation and geometric flexibility of finite-volume methods with
elementwise high-order polynomial approximation, and a large literature is
devoted to well-balanced
\cite{bermudez1994upwind,LEVEQUE1998346,audusse2004fast,xing2005high,xing2006high,xing2010positivity}
and entropy stable
\cite{tadmor1987numerical,tadmor2003entropy,fjordholm2009energy,fjordholm2011well,fjordholm2012arbitrarily}
discretizations of the SWE. Entropy stable high-order DG-type methods have
been developed for general systems of conservation laws and compressible
flows \cite{carpenter2014entropy,chen2017entropy,chan2018discretely} and for
the SWE in particular
\cite{gassner2016well,wintermeyer2017entropy,wintermeyer2018entropy,wu2021high}.
Of particular relevance to the present work is
the velocity-based skew-symmetric DG formulation of Fu \cite{fu2022high},
which is locally conservative, entropy stable, well-balanced and positivity
preserving, and which introduces the discrete velocity through an elementwise
$L^2$ projection rather than through a pointwise division.

High-order DG discretizations in physical space have been combined with SG
approximations before. For general uncertain hyperbolic systems, a DG scheme
with a multi-element SG ansatz and a hyperbolicity-preserving limiter was
developed in \cite{durrwachter2020hyperbolicity}. For the SWE, entropy
conservative DG methods for an SG system in Roe-type variables were proposed
in one dimension in \cite{bender2024entropy}. The recent work
\cite{offner2026high} compares that formulation with the
hyperbolicity-preserving formulation of
\cite{dai2021hyperbolicity,dai2022hyperbolicity}, shows that the two coincide
for stochastic bases whose product matrices share state-independent
eigenvectors, such as Haar-type wavelet bases, and extends the entropy
conservative DG method to two dimensions on tensor-product meshes. The
finite-volume schemes of \cite{dai2022hyperbolicity,epshteyn2025energy} use
the same SG formulation as the present paper but are limited to second order
on Cartesian grids.

In this paper we combine the hyperbolicity-preserving SG structure of
\cite{dai2022hyperbolicity} with the skew-symmetric DG structure of
\cite{fu2022high} into a high-order discontinuous Galerkin--stochastic
Galerkin (DG--SG) method for the two-dimensional stochastic SWE on
triangular meshes.
From \cite{fu2022high} we take only the idea of the skew-symmetric
formulation. The SG system is a different hyperbolic system, in which the
PCE coefficients of all unknowns are coupled, and it is not obvious a priori
that this idea applies to it; the skew-symmetric SG formulation, the DG
scheme and its entropy analysis are derived here for this system.
The one-dimensional case is treated in
our companion paper \cite{epshteyn2026structure}; the extension to two
dimensions on unstructured triangular meshes requires additional care in
several places, both in the analysis and in the implementation. First, the SG projection of the mixed flux $huv$ admits two
different representations, $\Pmat(\qxhat)\vhat$ and $\Pmat(\qyhat)\uhat$,
and the hyperbolicity-preserving construction selects one of them for each
momentum equation; the entropy analysis must account for this asymmetry.
Second, two velocity components are projected simultaneously on every
element, and the resulting discrete velocities enter the scheme both as
solution variables and as test functions. Third, on unstructured meshes the
interface dissipation, the hyperbolicity-preserving time step restriction,
the limiter and the boundary treatment must be formulated with respect to
the normal SG Jacobian on each face. Finally, in two dimensions the energy
exchanged through outflow boundaries is a line integral of the entropy flux,
so monitoring entropy stability on non-periodic domains requires an
augmented energy that accounts for this flux.

Compared with the closely related works above, our method differs in both
the SG formulation and the spatial discretization. For the polynomial chaos
bases considered here, such as Legendre and Jacobi polynomials, the two SG
formulations compared in \cite{offner2026high} are in general different; we
work directly with the conservative-variable SG formulation and the velocity
projection of \cite{fu2022high}, and the formulation and its analysis hold for
general distributions of the random inputs. The resulting scheme requires neither summation-by-parts
operators nor a flux-differencing form, and it is formulated on unstructured
triangular meshes. Moreover, the fully discrete algorithm includes slope
limiting, a computable hyperbolicity-preserving criterion and a treatment of
nearly dry regions; limiting and positivity preservation near dry states are
identified as open issues in \cite{offner2026high}.

\subsection{Contributions of this paper}

The main contributions of this paper are as follows.
\begin{itemize}[leftmargin=2em]
\item \emph{Continuous formulation.} We derive a two-dimensional
  \emph{skew-symmetric} SG formulation of the stochastic SWE. It is related
  to the hyperbolicity-preserving SG balance law of
  \cite{dai2022hyperbolicity} by an invertible state-dependent
  transformation, so that it is hyperbolic under the same
  positive-definiteness condition $\Pmat(\hhat)\succ0$, and it possesses the
  same companion entropy balance law with respect to a modified entropy
  variable.
\item \emph{Semi-discrete scheme.} We construct a high-order DG--SG scheme
  for this formulation on triangulations, with elementwise projected
  velocity components. We prove that it is locally conservative for the
  water height and for the discharge, well-balanced for stochastic
  lake-at-rest states, and entropy conservative with central fluxes or
  entropy stable with Lax--Friedrichs-type fluxes. The proofs use only the
  symmetry of the triple products that define the SG product operator and
  the velocity projection.
\item \emph{Fully discrete algorithm.} For completeness, we describe a fully
  discrete algorithm built on the semi-discrete scheme by adapting
  established deterministic DG techniques to the SG system: a desingularized
  velocity projection, a hyperbolicity-preserving time step restriction with
  linear scaling, a characteristic total-variation-bounded (TVB) limiter on
  equilibrium variables, and reflecting, periodic and outflow boundary
  conditions. To monitor the entropy stability on domains with outflow
  boundaries, we formulate the augmented energy of \cite{epshteyn2025energy}
  for the DG weak form; it serves as a diagnostic that separates the
  interior dissipation of the scheme from the physical energy exchange at
  the domain boundaries.
\item \emph{Numerical experiments.} We present two convergence studies, the
  stochastic two-dimensional test problems used for the finite-volume SG
  schemes of \cite{epshteyn2025energy}, two dam-break problems with uncertain
  water height, and problems with two and with ten independent,
  non-identically distributed random variables.
\end{itemize}

Throughout the paper we use the terms
entropy and energy interchangeably, since for the SWE the mathematical
entropy is the total (kinetic plus potential) energy.

The paper is organized as follows. In \cref{sec:sg-system} we review PCE
and the hyperbolicity-preserving SG formulation of the two-dimensional SWE,
derive the skew-symmetric SG system, and identify its companion entropy law.
In \cref{sec:dg-scheme} we construct the semi-discrete DG--SG scheme and
prove its local conservation, well-balanced, and entropy conservation and
stability properties. \Cref{sec:implementation} summarizes the fully discrete
algorithm; the details of the desingularization, the hyperbolicity-preserving
time step restriction, the limiter, the boundary conditions, and the
augmented energy are given in \cref{sec:appendix-algorithm}. Numerical experiments are presented in \cref{sec:numerics},
and \cref{sec:conclusion} contains concluding remarks.
\section{Skew-symmetric stochastic Galerkin formulation}
\label{sec:sg-system}

In this section we introduce a skew-symmetric stochastic Galerkin formulation
of the two-dimensional shallow water equations with random inputs. Our main results here are the skew-symmetric SG formulation \eqref{eq:sg-skew}, the establishment of its hyperbolicity (Proposition~\ref{prop:hyperbolicity}), and the companion balance law \eqref{eq:sg-entropy-density}--\eqref{eq:continuum-entropy-law} for the entropy.

\subsection{Parameterized two-dimensional shallow water equations}

Let $\Omega\subset\R^2$ be a bounded physical domain and let
$\bm\xi:\Theta\to\Gamma\subset\R^d$ be a random vector, defined on a
complete probability space $(\Theta,\mathcal F,\mathsf P)$, with probability density
$\rho$ on $\Gamma$. The random vector collects all uncertain inputs of the
model, for instance parameters of the bottom topography or of the initial
data, so that the solution depends on the spatial variables $(x,y)$, on time
$t$, and on $\bm\xi$. The parameterized state is
\[
  U=(h,q^x,q^y)\transpose,
  \qquad q^x=hu,\qquad q^y=hv,
\]
where $h=h(x,y,t,\bm\xi)$ is the water height, $(u,v)$ is the depth-averaged
velocity, $(q^x,q^y)$ is the discharge, and $B=B(x,y,\bm\xi)$ is a
time-independent bottom topography. For every fixed value of $\bm\xi$ the
state satisfies the deterministic shallow water equations
\begin{subequations}\label{eq:parameterized-swe}
\begin{align}
 h_t+(hu)_x+(hv)_y&=0,\\
 q^x_t+(hu^2)_x+(huv)_y+g h(h+B)_x&=0,\\
 q^y_t+(huv)_x+(hv^2)_y+g h(h+B)_y&=0,
\end{align}
\end{subequations}
where $g$ is the gravitational constant. Here the pressure flux and the
bottom source have been combined into the gradient of the water surface
$w:=h+B$; this form of the gravitational term is the natural one for the
preservation of the lake-at-rest equilibrium $w\equiv\mathrm{const}$,
$u=v=0$, which will be our notion of the well-balanced property.

Following \cite{fu2022high}, we also consider the skew-symmetric form of
\eqref{eq:parameterized-swe},
\begin{subequations}\label{eq:skew-parameterized}
\begin{align}
 h_t+(hu)_x+(hv)_y&=0,\label{eq:skew-h}\\
q^x_t+(hu^2)_x+(huv)_y+g h(h+B)_x
 -\frac12h_tu-\frac12(hu)_xu-\frac12(hv)_yu&=0,\label{eq:skew-qx}\\
q^y_t+(huv)_x+(hv^2)_y+g h(h+B)_y
 -\frac12h_tv-\frac12(hu)_xv-\frac12(hv)_yv&=0,\label{eq:skew-qy}
\end{align}
\end{subequations}
which is obtained by subtracting $\tfrac u2$ and $\tfrac v2$ times the
continuity equation from the two momentum equations. For smooth solutions
with positive water height, \eqref{eq:skew-parameterized} and
\eqref{eq:parameterized-swe} are equivalent. The additional terms vanish
identically at the continuous level, but they are the reason why, after
discretization, the kinetic energy contributions of the volume terms cancel
exactly \cite{fu2022high}. Skew-symmetric forms of systems of conservation
laws go back to \cite{tadmor1984skew}; for the SWE they were used in
entropy conservative DG methods in \cite{gassner2016well}. Our goal is to
extend this mechanism to the stochastic Galerkin setting in two spatial
dimensions.

\subsection{Polynomial chaos and the product matrix}

We briefly recall PCE; see
\cite{xiu2010numerical,le2010spectral,sullivan2015introduction} for further
details. Let $L^2_\rho(\Gamma)$ be the space of square-integrable functions
with respect to the density $\rho$, with inner product
\[
 \langle\phi,\psi\rangle_\rho
 :=\int_\Gamma\phi(\bm\xi)\psi(\bm\xi)\rho(\bm\xi)\dd\bm\xi .
\]
Assuming that $\rho$ has finite moments of all orders, there is a system of
polynomials that is orthonormal in $L^2_\rho(\Gamma)$; under mild conditions
on $\rho$ this system is complete, i.e., it is a basis of $L^2_\rho(\Gamma)$
\cite{ernst2012convergence}. We fix a finite multi-index set $\Lambda$ with
$|\Lambda|=\kpce$ elements, an enumeration $j=1,\dots,\kpce$ of its elements,
and the corresponding $\kpce$-dimensional polynomial subspace
$\mathbb P_\Lambda:=\operatorname{span}\{\phi_j\}_{j=1}^{\kpce}$, where the
basis polynomials are orthonormal and the first one is constant,
\[
 \langle\phi_j,\phi_\ell\rangle_\rho
 =\int_\Gamma\phi_j(\bm\xi)\phi_\ell(\bm\xi)
 \rho(\bm\xi)\dd\bm\xi=\delta_{j\ell},
 \qquad \phi_1\equiv1 .
\]
For a single random variable, $\Lambda=\{0,\dots,\kpce-1\}$ and $\phi_j$ has
degree $j-1$; for several random variables, $\Lambda$ is a set of
multi-indices, for instance a total-degree or a tensor-product set, and if
$\bm\xi$ has independent components then the basis is built from products of
univariate orthonormal polynomials.
For a scalar random field $z$ we write its $\kpce$-term PCE approximation as
\[
  z(x,y,t,\bm\xi)\simeq
  z_\Lambda(x,y,t,\bm\xi)
  =\sum_{j=1}^{\kpce}\widehat z_j(x,y,t)\phi_j(\bm\xi),
  \qquad
  \widehat z_j=\langle z,\phi_j\rangle_\rho,
\]
and we collect the coefficients in the vector
$\widehat z=(\widehat z_1,\dots,\widehat z_{\kpce})\transpose\in\R^{\kpce}$. The
statistics of $z_\Lambda$ follow directly from the coefficients: for example, since
$\phi_1\equiv1$, the mean and the variance are
\[
 \E[z_\Lambda]=\widehat z_1,
 \qquad
 \Var[z_\Lambda]=\sum_{j=2}^{\kpce}\widehat z_j^{\,2}.
\]

The SG formulation requires expansions of products, which can be computed
from the triple-product matrices
\[
  (M_j)_{\ell m}
  =\langle\phi_j,\phi_\ell\phi_m\rangle_\rho .
\]
These matrices define the linear SG product operator
$\Pmat:\R^{\kpce}\to\R^{\kpce\times\kpce}$,
\begin{equation}\label{eq:P-definition}
  \Pmat(\widehat z)=\sum_{j=1}^{\kpce}\widehat z_jM_j,
  \qquad
  [\Pmat(\widehat a)\widehat b]_\ell
  =\sum_{j,m=1}^{\kpce}\widehat a_j\widehat b_m
  \langle\phi_j,\phi_m\phi_\ell\rangle_\rho ,
\end{equation}
so that $\Pmat(\widehat a)\widehat b$ is the vector of PCE coefficients of
the Galerkin projection of the product $a_\Lambda b_\Lambda$ onto $\mathbb P_\Lambda$.
Each $M_j$ is symmetric, hence so is $\Pmat(\widehat z)$ for every
$\widehat z$, and the symmetry of the triple product in its three arguments
gives the commutativity identities
\begin{align}
 \Pmat(\widehat a)\widehat b
 &=\Pmat(\widehat b)\widehat a,\label{eq:P-comm}\\
 [\Pmat(\widehat a)\widehat b]\cdot\widehat c
 &=[\Pmat(\widehat a)\widehat c]\cdot\widehat b
 =[\Pmat(\widehat b)\widehat c]\cdot\widehat a.
\label{eq:P-triple}
\end{align}
These identities play the role of the commutativity of scalar
multiplication in the entropy calculations below; see also
\cite[Lemma~2.1]{dai2024energy}.

In the deterministic equations the velocity is $u=q/h$. In the SG setting
the role of the division by $h$ is played by the inverse of the 
matrix $\Pmat(\hhat)$: the SG velocity coefficients are defined by
\begin{equation}\label{eq:sg-velocity}
 \uhat=\Pmat(\hhat)^{-1}\qxhat,\qquad
 \vhat=\Pmat(\hhat)^{-1}\qyhat,
\end{equation}
which are the coefficients $\uhat,\vhat$ such that the Galerkin projections
of $h_\Lambda u_\Lambda$ and $h_\Lambda v_\Lambda$ coincide with $q^x_\Lambda$ and $q^y_\Lambda$. These relations
are well defined under the positive-definiteness condition
\begin{equation}\label{eq:hyperbolicity-condition}
 \Pmat(\hhat(x,y,t))\succ0 ,
\end{equation}
which is the stochastic counterpart of the positivity $h>0$ of the water
height and which, as we recall next, is also the condition under which the SG
system is hyperbolic.

\subsection{The two-dimensional skew-symmetric SG system}
\label{sec:skew-sg-system}

We recall how the hyperbolicity-preserving SG formulation of the
two-dimensional SWE is obtained \cite{dai2022hyperbolicity}, following the
standard Galerkin procedure in the stochastic space. Starting from
\eqref{eq:parameterized-swe}, the solution variables $(h,q^x, q^y)$ are replaced by their
truncated polynomial chaos approximations $(h_\Lambda, q^x_\Lambda, q^y_\Lambda)$:
\begin{equation}
\begin{aligned}
 z(x,y,t,\bm\xi)&\simeq z_\Lambda(x,y,t,\bm\xi)
  =\sum_{j=1}^{\kpce}\widehat z_j(x,y,t)\phi_j(\bm\xi), \hskip 20pt z \in \{h, q^x, q^y\},
\end{aligned}
\end{equation}
and the bottom topography is replaced by its stochastic projection
\begin{equation}
 B(x,y,\bm\xi)\simeq B_\Lambda(x,y,\bm\xi)
 =\sum_{j=1}^{\kpce}\widehat B_j(x,y)\phi_j(\bm\xi).
\end{equation}
The residuals are then projected onto the same stochastic space. The
polynomial terms of \eqref{eq:parameterized-swe} are projected exactly by the
product operator \eqref{eq:P-definition}. For the nonlinear, non-polynomial
terms $(q^x)^2/h$, $(q^y)^2/h$, and $q^xq^y/h$, the hyperbolicity-preserving
construction of \cite{dai2022hyperbolicity} uses the projected velocity
coefficients \eqref{eq:sg-velocity} and a directional choice for the two
representations of the mixed flux, which we make explicit below. The
resulting SG formulation is the system of balance laws for the unknown $\widehat U$:
\begin{align}\label{eq:sg-balance-law}
  \widehat U&=(\hhat\transpose,(\qxhat)\transpose,(\qyhat)\transpose)\transpose
  \in\R^{3\kpce}, & 
 \widehat U_t+\widehat F(\widehat U)_x
  +\widehat G(\widehat U)_y&=\widehat S(\widehat U,\Bhat), 
\end{align}
where
\begin{subequations}\label{eq:sg-fluxes-source}
\begin{align}
 \widehat F(\widehat U)
 &=\begin{pmatrix}
 \qxhat\\
 \Pmat(\qxhat)\uhat+\tfrac12g\Pmat(\hhat)\hhat\\
 \Pmat(\qxhat)\vhat
 \end{pmatrix},
 &
 \widehat G(\widehat U)
 &=\begin{pmatrix}
 \qyhat\\
 \Pmat(\qyhat)\uhat\\
 \Pmat(\qyhat)\vhat+\tfrac12g\Pmat(\hhat)\hhat
 \end{pmatrix},\label{eq:sg-fluxes}\\
 \widehat S(\widehat U,\Bhat)
 &=\begin{pmatrix}
 \bm 0\\
 -g\Pmat(\hhat)(\Bhat)_x\\
 -g\Pmat(\hhat)(\Bhat)_y
 \end{pmatrix}.
 \label{eq:sg-source}
\end{align}
\end{subequations}
In the deterministic case $\kpce=1$ the fluxes and the source reduce to their
classical counterparts. A feature specific to two space dimensions is that
the deterministic mixed flux $huv$ has two stochastic representations, and
they are in general different:
\begin{equation}\label{eq:mixed-flux-representations}
 \Pmat(\qxhat)\vhat\ne\Pmat(\qyhat)\uhat.
\end{equation}
In \eqref{eq:sg-fluxes} the first expression is used in the $x$-flux of the
$y$-momentum equation, whereas the second is used in the $y$-flux of the
$x$-momentum equation. This directional choice is part of the
hyperbolicity-preserving construction in \cite{dai2022hyperbolicity}. With
this choice, each momentum equation transports a velocity component with the
discharge: up to the pressure terms, the fluxes of $\qxhat$ and $\qyhat$ in
a direction $\bm n$ are $\Pmat(\qxhat n_x+\qyhat n_y)\uhat$ and
$\Pmat(\qxhat n_x+\qyhat n_y)\vhat$, respectively. This is the structure
that the entropy analysis below requires. For instance, by the
triple-product symmetry \eqref{eq:P-triple},
\[
 \bigl(\Pmat(\qxhat)\vhat\bigr)_x\cdot\vhat
 -\tfrac12\Pmat\bigl((\qxhat)_x\bigr)\vhat\cdot\vhat
 =\tfrac12\bigl(\Pmat(\qxhat)\vhat\cdot\vhat\bigr)_x
\]
is an exact derivative, whereas the same expression with
$\Pmat(\qyhat)\uhat$ in place of $\Pmat(\qxhat)\vhat$ is in general not.

When \eqref{eq:hyperbolicity-condition} holds, the flux Jacobians of
\eqref{eq:sg-balance-law} are well defined. In terms of $\kpce\times\kpce$ blocks,
they are
\begin{equation}\label{eq:sg-flux-jacobians}
\begin{aligned}
 \frac{\partial\widehat F}{\partial\widehat U}
 &=\begin{pmatrix}
 0&I&0\\
 g\Pmat(\hhat)-\Pmat(\qxhat)\Pmat(\hhat)^{-1}\Pmat(\uhat)
 &\Pmat(\qxhat)\Pmat(\hhat)^{-1}+\Pmat(\uhat)&0\\
 -\Pmat(\qxhat)\Pmat(\hhat)^{-1}\Pmat(\vhat)
 &\Pmat(\vhat)&\Pmat(\qxhat)\Pmat(\hhat)^{-1}
 \end{pmatrix},\\[0.4em]
 \frac{\partial\widehat G}{\partial\widehat U}
 &=\begin{pmatrix}
 0&0&I\\
 -\Pmat(\qyhat)\Pmat(\hhat)^{-1}\Pmat(\uhat)
 &\Pmat(\qyhat)\Pmat(\hhat)^{-1}&\Pmat(\uhat)\\
 g\Pmat(\hhat)-\Pmat(\qyhat)\Pmat(\hhat)^{-1}\Pmat(\vhat)
 &0&\Pmat(\qyhat)\Pmat(\hhat)^{-1}+\Pmat(\vhat)
 \end{pmatrix}.
\end{aligned}
\end{equation}
For a unit vector $\bm n=(n_x,n_y)\transpose$, the corresponding normal
Jacobian is
\begin{equation}\label{eq:sg-normal-jacobian}
 A_{\bm n}(\widehat U)
 =n_x\frac{\partial\widehat F}{\partial\widehat U}
  +n_y\frac{\partial\widehat G}{\partial\widehat U}.
\end{equation}
The main result of \cite{dai2022hyperbolicity} is that, under
\eqref{eq:hyperbolicity-condition}, $A_{\bm n}(\widehat U)$ is diagonalizable
with real eigenvalues for every unit vector $\bm n$, i.e., the SG system
\eqref{eq:sg-balance-law} is hyperbolic. The spectral radius of
$A_{\bm n}$ will be used in \cref{sec:dg-scheme,sec:time-stepping} as the
wave speed in the numerical dissipation and in the time step restriction.

Our approach to developing spatially high-order DG schemes for the SG shallow water system employs not the standard form of the balance law, but its skew-symmetric form. We therefore state the standard notions of hyperbolicity and of a companion balance law \cite{lax1973hyperbolic,leveque2002finite,dafermos2016hyperbolic} in a form that applies to such systems.

\begin{definition}[Hyperbolicity for a first-order system]
\label{def:first-order-hyperbolicity}
A first-order system of PDEs is called hyperbolic if, after left
multiplication by an invertible matrix $L(U)$ that depends only on the state
and contains no temporal or spatial derivatives, it can be written in the
quasilinear form
\begin{equation}\label{eq:general-quasilinear-system}
 U_t+A^x(U)U_x+A^y(U)U_y=S(U),
\end{equation}
and, for every unit vector $\bm n=(n_x,n_y)\transpose$, the normal matrix
\begin{equation}\label{eq:general-normal-matrix}
 A_{\bm n}(U)=n_xA^x(U)+n_yA^y(U)
\end{equation}
is diagonalizable with real eigenvalues
\cite{lax1973hyperbolic,leveque2002finite}.
\end{definition}

\begin{definition}[Companion balance law for a first-order system]
\label{def:first-order-companion-law}
Consider an arbitrary first-order system that, through left multiplication
by an invertible matrix $L(U)$ depending only on the state, can be transformed
to the balance law
\begin{equation}\label{eq:general-balance-law}
 U_t+F(U)_x+G(U)_y=S(U).
\end{equation}
For smooth solutions $U$, the system has a companion balance law if there are a convex
scalar function $E(U)$ and fluxes $H^x(U)$ and $H^y(U)$ such that
\begin{equation}\label{eq:general-companion-law}
 E(U)_t+H^x(U)_x+H^y(U)_y=0.
\end{equation}
Equivalently, the tuple $(E,H^x,H^y)$ satisfies the compatibility relation
\begin{equation}\label{eq:general-companion-compatibility}
 \frac{\partial E}{\partial U}
 \bigl(F(U)_x+G(U)_y-S(U)\bigr)
 =H^x(U)_x+H^y(U)_y.
\end{equation}
\end{definition}

If $V=(\partial E/\partial U)\transpose$ is the entropy variable for
\eqref{eq:general-balance-law}, then the multiplier that recovers the same
companion law from the left-transformed first-order system is
\begin{equation}\label{eq:modified-entropy-variable-general}
 \widetilde V=(L\transpose)^{-1}V.
\end{equation}
We emphasize that, in general, $\widetilde V$ is not the derivative of $E$ with
respect to $LU$, because $L$ may depend on $U$. What holds instead is the
identity
\begin{equation}\label{eq:modified-entropy-time-identity}
 \widetilde V\transpose L U_t
 =V\transpose U_t=E(U)_t,
\end{equation}
which is all that is needed to obtain \eqref{eq:general-companion-law} from
the transformed system.

\begin{remark}[Deterministic case]
For the deterministic SWE, the skew-symmetric system
\eqref{eq:skew-parameterized} is obtained from \eqref{eq:parameterized-swe}
by left multiplication with the invertible matrix
$L(U)=\bigl(\begin{smallmatrix}1&0&0\\-\frac u2&1&0\\-\frac v2&0&1\end{smallmatrix}\bigr)$,
and it is therefore hyperbolic whenever $h>0$. With the energy
$E=\frac12h(u^2+v^2)+\frac12gh^2+ghB$ and its entropy variable
$V=(g(h+B)-\frac12(u^2+v^2),u,v)\transpose$, the modified entropy variable
\eqref{eq:modified-entropy-variable-general} is
$\widetilde V=(g(h+B),u,v)\transpose$: multiplying \eqref{eq:skew-h} by
$g(h+B)$, \eqref{eq:skew-qx} by $u$, and \eqref{eq:skew-qy} by $v$ and adding
the results recovers the energy balance law of the SWE \cite{fu2022high}.
The stochastic construction below follows the same pattern.
\end{remark}

We now apply the Galerkin procedure of \eqref{eq:sg-balance-law} to the
skew-symmetric system \eqref{eq:skew-parameterized}: the solution variables
and the bottom are replaced by their truncated PCE approximations, the
products are projected through the operator $\Pmat$, and the velocities are
the projected velocities \eqref{eq:sg-velocity}. This yields the
skew-symmetric SG formulation of the two-dimensional SWE,
\begin{subequations}\label{eq:sg-skew}
\begin{align}
 \hhat_t+(\qxhat)_x+(\qyhat)_y&=0,\label{eq:sg-h}\\
 \qxhat_t+\bigl(\Pmat(\qxhat)\uhat\bigr)_x
 +\bigl(\Pmat(\qyhat)\uhat\bigr)_y
 +g\Pmat(\hhat)(\hhat+\Bhat)_x
 &-\frac12\Pmat(\hhat_t)\uhat
 -\frac12\Pmat((\qxhat)_x)\uhat
 -\frac12\Pmat((\qyhat)_y)\uhat=0,\label{eq:sg-qx}\\
 \qyhat_t+\bigl(\Pmat(\qxhat)\vhat\bigr)_x
 +\bigl(\Pmat(\qyhat)\vhat\bigr)_y
 +g\Pmat(\hhat)(\hhat+\Bhat)_y
 &-\frac12\Pmat(\hhat_t)\vhat
 -\frac12\Pmat((\qxhat)_x)\vhat
 -\frac12\Pmat((\qyhat)_y)\vhat=0.\label{eq:sg-qy}
\end{align}
\end{subequations}
In the deterministic case $\kpce=1$, \eqref{eq:sg-skew} reduces to
\eqref{eq:skew-parameterized}. The following proposition shows that
\eqref{eq:sg-skew} inherits the hyperbolicity of \eqref{eq:sg-balance-law}.

\begin{proposition}[Hyperbolicity of \eqref{eq:sg-skew}]
\label{prop:hyperbolicity}
Assume \eqref{eq:hyperbolicity-condition}. The skew-symmetric SG system
\eqref{eq:sg-skew} is obtained from the hyperbolicity-preserving SG balance
law \eqref{eq:sg-balance-law} by left multiplication with the block matrix
\begin{equation}\label{eq:skew-transformation}
 L_m(\widehat U)=
 \begin{pmatrix}
 I&0&0\\
 -\tfrac12\Pmat(\uhat)&I&0\\
 -\tfrac12\Pmat(\vhat)&0&I
 \end{pmatrix}.
\end{equation}
The matrix $L_m$ is invertible, and the skew-symmetric SG system is
hyperbolic in the sense of Definition~\ref{def:first-order-hyperbolicity}. In particular, the normal matrix
\eqref{eq:general-normal-matrix} of the skew-symmetric SG system is the
normal Jacobian \eqref{eq:sg-normal-jacobian}, which is diagonalizable with real eigenvalues.
\end{proposition}

\begin{proof}
The block lower-triangular matrix \eqref{eq:skew-transformation} has identity
diagonal blocks and inverse
\[
 L_m(\widehat U)^{-1}=
 \begin{pmatrix}
 I&0&0\\
 \tfrac12\Pmat(\uhat)&I&0\\
 \tfrac12\Pmat(\vhat)&0&I
 \end{pmatrix}.
\]
Multiplying the height residual in \eqref{eq:sg-balance-law} by
$-\tfrac12\Pmat(\uhat)$ and adding it to the $x$-momentum residual gives
\eqref{eq:sg-qx}, because $\Pmat(\uhat)\qxhat_{\,x}=\Pmat((\qxhat)_x)\uhat$ by
\eqref{eq:P-comm} and $\Pmat(\uhat)\hhat_t=\Pmat(\hhat_t)\uhat$;
multiplication by $-\tfrac12\Pmat(\vhat)$ gives \eqref{eq:sg-qy} in the same
way. Thus the two systems are related by an invertible state-dependent left
transformation containing no derivatives. Under
\eqref{eq:hyperbolicity-condition}, the balance law
\eqref{eq:sg-balance-law} is hyperbolic \cite{dai2022hyperbolicity}; that is,
the normal Jacobian \eqref{eq:sg-normal-jacobian} is diagonalizable with
real eigenvalues for every unit $\bm n$. The conclusion follows directly
from Definition~\ref{def:first-order-hyperbolicity}.
\end{proof}

\begin{remark}
The condition \eqref{eq:hyperbolicity-condition} is a condition on the whole
coefficient vector $\hhat$ rather than on the mean coefficient alone. A
sufficient and computable criterion is the positivity of the truncated
water-height expansion $h_\Lambda$ at the nodes of a positive stochastic quadrature
rule of sufficiently high degree \cite[Theorem~3.2]{dai2021hyperbolicity};
this criterion is the basis of the hyperbolicity-preserving time step
restriction in \cref{sec_hpc}.
\end{remark}

We now identify the companion balance law of \eqref{eq:sg-skew}. Define the
SG entropy (energy) density
\begin{equation}\label{eq:sg-entropy-density}
 E(\hhat,\qxhat,\qyhat)
 =\frac12g\norml{\hhat}^2
  +g\hhat\transpose\Bhat
  +\frac12\qxhat\cdot\uhat
  +\frac12\qyhat\cdot\vhat ,
\end{equation}
which is the expected total energy
$\E[\frac12h_\Lambda(u_\Lambda^2+v_\Lambda^2)+\frac12gh_\Lambda^2+gh_\Lambda B_\Lambda]$
of the truncated fields. For a time-independent bottom, the
corresponding physical-space entropy fluxes are
\begin{subequations}\label{eq:sg-entropy-flux}
\begin{align}
 H^x
 &=\frac12\Pmat(\qxhat)\uhat\cdot\uhat
   +\frac12\Pmat(\qxhat)\vhat\cdot\vhat
   +g\qxhat\cdot(\hhat+\Bhat),\\
 H^y
 &=\frac12\Pmat(\qyhat)\uhat\cdot\uhat
   +\frac12\Pmat(\qyhat)\vhat\cdot\vhat
   +g\qyhat\cdot(\hhat+\Bhat).
\end{align}
\end{subequations}
Under \eqref{eq:hyperbolicity-condition}, $E$ is convex in
$\widehat U$ \cite{epshteyn2025energy}. Its entropy variable for the balance
law \eqref{eq:sg-balance-law} is
\begin{equation}\label{eq:sg-entropy-variable}
 \widehat V
 =\frac{\partial E}{\partial\widehat U}\transpose
 =\begin{pmatrix}
 g(\hhat+\Bhat)-\tfrac12\Pmat(\uhat)\uhat
                  -\tfrac12\Pmat(\vhat)\vhat\\
 \uhat\\
 \vhat
 \end{pmatrix},
\end{equation}
and the modified entropy variable \eqref{eq:modified-entropy-variable-general}
associated with the skew transformation \eqref{eq:skew-transformation} is
\begin{equation}\label{eq:sg-modified-entropy-variable}
 \widetilde{\widehat V}
 =(L_m\transpose)^{-1}\widehat V
 =\begin{pmatrix}
 g(\hhat+\Bhat)\\
 \uhat\\
 \vhat
 \end{pmatrix}.
\end{equation}
Consequently, multiplying \eqref{eq:sg-h} by $g(\hhat+\Bhat)$,
\eqref{eq:sg-qx} by $\uhat$, and \eqref{eq:sg-qy} by $\vhat$, adding the
results and using \eqref{eq:P-triple}, gives for smooth solutions
\begin{equation}\label{eq:continuum-entropy-law}
 E_t+\nabla\cdot\bm H=0,
 \qquad
 \bm H=(H^x,H^y)\transpose.
\end{equation}
Thus $(E,H^x,H^y)$ is a companion entropy-flux
tuple for both the balance law \eqref{eq:sg-balance-law} and its
skew-symmetric form \eqref{eq:sg-skew}. Observe that the modified entropy
variable \eqref{eq:sg-modified-entropy-variable} is linear in the water
surface and in the projected velocities. This is what makes the
skew-symmetric form attractive for a Galerkin discretization: the entropy
variables are themselves admissible test functions, and the identity
\eqref{eq:continuum-entropy-law}, which holds separately in each coordinate
direction, is the template for the discrete entropy balance derived in the
next section.
\section{Structure-preserving DG--SG discretization}
\label{sec:dg-scheme}

In this section we derive a DG--SG discretization of the skew-symmetric SG
system \eqref{eq:sg-skew}
on triangular meshes, and we establish its local conservation,
well-balanced, and entropy conservation and stability properties through
appropriate choices of the numerical fluxes. Throughout this section and the
rest of the paper, the subscript $\Delta$ denotes quantities after spatial
discretization.

\subsection{Triangulation, traces, and approximation spaces}

Let $\Th$ be a shape-regular triangulation of $\Omega$ into triangles,
let $\Eh$ denote the set of all faces (edges) of $\Th$, and let
$\Ehi\subset\Eh$ and $\Ehb\subset\Eh$ denote the sets of interior and
boundary faces, respectively.
The subscript $h$ of the mesh $\Th$ and of the face sets is the customary
mesh label and does not refer to the water height.
For a polynomial degree $k\ge0$, let $\mathbb P^k(K)$ denote the subspace of at-most degree-$k$ polynomials on region $K$, and
\[
 \Vk:=\bigl\{\vartheta\in L^2(\Omega):\ \vartheta|_K\in\mathbb P^k(K)\ \ \forall K\in\Th\bigr\}
\]
be the discontinuous finite-element space of piecewise polynomials of total
degree at most $k$ \cite{cockburn2001runge,hesthaven2008nodal}. On an interior face $e=K^-\cap K^+$ we fix a unit normal
$\bm n=(n_x,n_y)\transpose$ pointing from $K^+$ into $K^-$ (so that $\bm n$
is the outward normal of $K^+$), denote by $a^\mp$ the traces of a function
$a$ from $K^\mp$, and define the average and jump
\[
 \avg{a}=\frac12(a^-+a^+),\qquad
 \jump{a}=a^+-a^-.
\]
On a boundary face the interior trace $a^-$ is taken from the element
containing the face, and the exterior trace $a^+$ is a ghost state
determined by the boundary condition; see \cref{sec:boundary}.
The DG--SG unknowns are vector-valued functions whose components are the
PCE coefficients,
\[
 \hD,\qxD,\qyD,\uD,\vD\in [\Vk]^{\kpce},
\]
where $\kpce$ is the number of PCE terms; we write $\widehat h_{\Delta,j}$
for the $j$th component of $\hD$, and similarly for the other unknowns. The
products $\Pmat(\cdot)$ and the dot products below act in the SG coefficient
space $\R^{\kpce}$ at every point of $\Omega$, and the bottom $\BD\in
[\Vk]^{\kpce}$ is the DG projection of the PCE coefficients of $B$.

\subsection{Weak formulation}

The semi-discrete DG--SG formulation of \eqref{eq:sg-skew} reads as follows:
find $(\hD,\qxD,\qyD,\uD,\vD)$, functions of time with values in
$[\Vk]^{\kpce}$, such that the equations below hold. For compactness
the equations are written after summing over all elements and faces, and
they are imposed componentwise in SG space. The height equation is: for every
$\bm\psi\in[\Vk]^{\kpce}$,
\begin{align}
 \sum_{K\in\Th}\int_K(\hD)_t\cdot\bm\psi\dd\bm x
 &-\sum_{K\in\Th}\int_K
  \left(\qxD\cdot\bm\psi_x+\qyD\cdot\bm\psi_y\right)\dd\bm x
 +\sum_{e\in\Ehi}\int_e\fq\cdot\jump{\bm\psi}\dd s=0,
\label{eq:dg-h}
\end{align}
where we choose the normal height flux as
\begin{equation}\label{eq:flux-h}
 \fq
 =\avg{\Pmat(\hD)\uD}n_x
  +\avg{\Pmat(\hD)\vD}n_y
  +\frac{\alpha_1}{2}\jump{\hD+\BD}.
\end{equation}
Here and below $\alpha_1\ge0$ is a dissipation coefficient, a single-valued
function on the faces: $\alpha_1=0$ gives the entropy conservative scheme,
and $\alpha_1>0$ the entropy stable scheme (\cref{thm:ec-es}). In the
numerical experiments $\alpha_1$ is the spectral radius of the normal SG
Jacobian \eqref{eq:sg-normal-jacobian}, evaluated from the two traces on the
face. The $x$-momentum equation is written componentwise:
for every $j=1,\dots,\kpce$ and every $\widehat\zeta^{x}_{\Delta,j}\in \Vk$,
\begin{align}
  &\sum_{K\in\Th}\int_K(\widehat q^{x}_{\Delta,j})_t\widehat\zeta^{x}_{\Delta,j}\dd\bm x
  -\sum_{K\in\Th}\int_K\bigl(\Pmat(\qxD)\uD\bigr)_j(\widehat\zeta^{x}_{\Delta,j})_x\dd\bm x
  -\sum_{K\in\Th}\int_K\bigl(\Pmat(\qyD)\uD\bigr)_j(\widehat\zeta^{x}_{\Delta,j})_y\dd\bm x\notag\\
  &+\sum_{e\in\Ehi}\int_e(\fqu)_j\jump{\widehat\zeta^{x}_{\Delta,j}}\dd s\notag + \sum_{K\in\Th}\int_Kg\bigl(\Pmat(\hD)(\hD+\BD)_x\bigr)_j\widehat\zeta^{x}_{\Delta,j}\dd\bm x\notag\\
  &-\sum_{e\in\Ehi}\int_eg\bigl(\Pmat(\jump{\hD+\BD})\avg{\widehat\zeta^{x}_{\Delta,j}\hD}\bigr)_j n_x\dd s\notag -\frac12\sum_{K\in\Th}\int_K\bigl(\Pmat((\hD)_t)\uD\bigr)_j\widehat\zeta^{x}_{\Delta,j}\dd\bm x\notag\\
  &+\frac12\sum_{K\in\Th}\int_K\bigl(\Pmat(\qxD)(\widehat\zeta^{x}_{\Delta,j}\uD)_x\bigr)_j\dd\bm x
  +\frac12\sum_{K\in\Th}\int_K\bigl(\Pmat(\qyD)(\widehat\zeta^{x}_{\Delta,j}\uD)_y\bigr)_j\dd\bm x\notag\\
  &-\frac12\sum_{e\in\Ehi}\int_e\bigl(\Pmat(\fq)\jump{\widehat\zeta^{x}_{\Delta,j}\uD}\bigr)_j\dd s=0,
  \label{eq:dg-qx}
\end{align}
where the subscript $j$ denotes the $j$th SG component, $\jump{\cdot}$ and
$\avg{\cdot}$ of a product act on the product, and the numerical momentum
flux is
\begin{equation}\label{eq:flux-qu}
  \fqu
  :=\Pmat\left(\avg{\Pmat(\hD)\uD}n_x+\avg{\Pmat(\hD)\vD}n_y\right)\avg{\uD}
  +\frac12\alpha_1\jump{\Pmat(\hD+\BD)\uD}.
\end{equation}
Analogously, the $y$-momentum equation is: for every $j=1,\dots,\kpce$ and
every $\widehat\zeta^{y}_{\Delta,j}\in \Vk$,
\begin{align}
 & \sum_{K\in\Th}\int_K(\widehat q^{y}_{\Delta,j})_t\widehat\zeta^{y}_{\Delta,j}\dd\bm x
  -\sum_{K\in\Th}\int_K\bigl(\Pmat(\qxD)\vD\bigr)_j(\widehat\zeta^{y}_{\Delta,j})_x\dd\bm x
  -\sum_{K\in\Th}\int_K\bigl(\Pmat(\qyD)\vD\bigr)_j(\widehat\zeta^{y}_{\Delta,j})_y\dd\bm x\notag\\
  &+\sum_{e\in\Ehi}\int_e(\fqv)_j\jump{\widehat\zeta^{y}_{\Delta,j}}\dd s\notag + \sum_{K\in\Th}\int_Kg\bigl(\Pmat(\hD)(\hD+\BD)_y\bigr)_j\widehat\zeta^{y}_{\Delta,j}\dd\bm x\notag\\
  &-\sum_{e\in\Ehi}\int_eg\bigl(\Pmat(\jump{\hD+\BD})\avg{\widehat\zeta^{y}_{\Delta,j}\hD}\bigr)_j n_y\dd s\notag -\frac12\sum_{K\in\Th}\int_K\bigl(\Pmat((\hD)_t)\vD\bigr)_j\widehat\zeta^{y}_{\Delta,j}\dd\bm x\notag\\
  &+\frac12\sum_{K\in\Th}\int_K\bigl(\Pmat(\qxD)(\widehat\zeta^{y}_{\Delta,j}\vD)_x\bigr)_j\dd\bm x
  +\frac12\sum_{K\in\Th}\int_K\bigl(\Pmat(\qyD)(\widehat\zeta^{y}_{\Delta,j}\vD)_y\bigr)_j\dd\bm x\notag\\
  &-\frac12\sum_{e\in\Ehi}\int_e\bigl(\Pmat(\fq)\jump{\widehat\zeta^{y}_{\Delta,j}\vD}\bigr)_j\dd s=0,
  \label{eq:dg-qy}
\end{align}
with the numerical flux
\begin{equation}\label{eq:flux-qv}
  \fqv
  :=\Pmat\left(\avg{\Pmat(\hD)\uD}n_x+\avg{\Pmat(\hD)\vD}n_y\right)\avg{\vD}
  +\frac12\alpha_1\jump{\Pmat(\hD+\BD)\vD}.
\end{equation}

The face sums in \eqref{eq:dg-h}, \eqref{eq:dg-qx} and \eqref{eq:dg-qy} run
over the interior faces only: the formulation is written for a mesh without
boundary elements, i.e., every face is an interior face, or equivalently,
periodic boundary conditions are imposed. On physical boundaries the same
face terms are added on the boundary faces $e\in\Ehb$, with the exterior
traces of the solution replaced by ghost states; see \cref{sec:boundary}.

The scheme is closed by the definition of the discrete velocities. They are
not obtained by inverting $\Pmat(\hD)$ pointwise at the quadrature nodes,
which would in general produce a velocity outside the DG space; instead, on
every element $K\in\Th$, the PCE coefficients of both velocity components
are defined by the elementwise projection
\begin{subequations}\label{eq:velocity-projection}
\begin{align}
 \int_K[\Pmat(\hD)\uD-\qxD]\cdot\bm\zeta^x\dd\bm x&=0,
 &&\forall\bm\zeta^x\in[\Vk]^{\kpce},\\
 \int_K[\Pmat(\hD)\vD-\qyD]\cdot\bm\zeta^y\dd\bm x&=0,
 &&\forall\bm\zeta^y\in[\Vk]^{\kpce}.
\end{align}
\end{subequations}
Equation \eqref{eq:velocity-projection} states that $\qxD$ and $\qyD$ are
the $L^2$ projections of $\Pmat(\hD)\uD$ and $\Pmat(\hD)\vD$ onto
$[\Vk]^{\kpce}$ on every element, which is the discrete counterpart
of \eqref{eq:sg-velocity} and the two-dimensional, stochastic analog of
the velocity definition in \cite{fu2022high}. It is well defined whenever
the elementwise mass matrix of $\Pmat(\hD)$ is positive definite, which is
the case under the discrete hyperbolicity condition $\Pmat(\hD)\succ0$ at
the quadrature nodes (see \cref{sec_hpc}). Because $\uD$ and $\vD$ belong to
the DG space, they are admissible test functions in \eqref{eq:dg-qx} and
\eqref{eq:dg-qy}; this fact, together with the symmetry of $\Pmat$, is what
drives the discrete entropy calculation below.

\subsection{Local conservation and well-balanced property}

We first record the local conservation properties of the DG--SG scheme,
which are the stochastic, two-dimensional analogs of the local
conservation statement for the deterministic velocity-based DG method in
\cite{fu2022high} and of the one-dimensional DG--SG result
\cite{epshteyn2026structure}. For an element
$K\in\Th$ we denote by $\bm n^K=(n^K_x,n^K_y)\transpose$ its outward unit
normal, by $a^{\mathrm{in}}$ and $a^{\mathrm{ext}}$ the traces on
$\partial K$ of a function $a$ from inside $K$ and from the neighboring
element (or the ghost state on a boundary face), respectively, and by
\[
 \jump{a}_K:=a^{\mathrm{ext}}-a^{\mathrm{in}}
\]
the jump across $\partial K$ in the outward direction. A numerical flux
with the superscript $K$, such as $\fq^{\,K}$, is the flux
\eqref{eq:flux-h}, \eqref{eq:flux-qu} or \eqref{eq:flux-qv} evaluated with
$\bm n=\bm n^K$, i.e., with $(\cdot)^+=(\cdot)^{\mathrm{in}}$ and
$(\cdot)^-=(\cdot)^{\mathrm{ext}}$. Note that the three fluxes change sign
when the orientation of the face is reversed (the averages are unchanged,
while $\bm n$ and the jumps change sign), so that $\fq^{\,K'}=-\fq^{\,K}$ on
a face shared by $K$ and $K'$, and similarly for $\fqu$ and $\fqv$.

\begin{theorem}[Local conservation of the DG--SG scheme]\label{thm:local-conserve}
For every element $K\in\Th$ and each stochastic mode $j=1,\dots,\kpce$, the
semi-discrete DG--SG scheme \eqref{eq:dg-h}, \eqref{eq:dg-qx},
\eqref{eq:dg-qy} with the velocity projection \eqref{eq:velocity-projection}
satisfies
\begin{equation}\label{eq:local-conserve-h}
 \frac{\dd}{\dd t}\int_K\widehat h_{\Delta,j}\dd\bm x
 =-\int_{\partial K}(\fq^{\,K})_j\dd s,
\end{equation}
and
\begin{equation}\label{eq:local-conserve-q}
\begin{aligned}
 \frac{\dd}{\dd t}\int_K\widehat q^{x}_{\Delta,j}\dd\bm x
 ={}&-\int_{\partial K}\Bigl(\fqu^{\,K}
   +\frac12g\,\Pmat(\hD^{\mathrm{ext}})\hD^{\mathrm{in}}\,n^K_x\Bigr)_j\dd s
   -\int_K g\bigl(\Pmat(\hD)(\BD)_x\bigr)_j\dd\bm x\\
 &-\frac12g\int_{\partial K}\bigl(\Pmat(\hD^{\mathrm{in}})\jump{\BD}_K\bigr)_j n^K_x\dd s,\\
 \frac{\dd}{\dd t}\int_K\widehat q^{y}_{\Delta,j}\dd\bm x
 ={}&-\int_{\partial K}\Bigl(\fqv^{\,K}
   +\frac12g\,\Pmat(\hD^{\mathrm{ext}})\hD^{\mathrm{in}}\,n^K_y\Bigr)_j\dd s
   -\int_K g\bigl(\Pmat(\hD)(\BD)_y\bigr)_j\dd\bm x\\
 &-\frac12g\int_{\partial K}\bigl(\Pmat(\hD^{\mathrm{in}})\jump{\BD}_K\bigr)_j n^K_y\dd s.
\end{aligned}
\end{equation}
\end{theorem}

\begin{proof}
Let $\chi_K$ denote the characteristic function of $K$, and let $\bm\delta_j$
be the $j$th unit vector of $\R^{\kpce}$.

\smallskip\noindent\emph{Height.} Take in \eqref{eq:dg-h} the test
function $\bm\psi=\chi_K\bm\delta_j$. The volume terms vanish because
$\chi_K$ is constant on every element, and only the faces $e\subset\partial K$
contribute to the face term. On such a face, $\jump{\chi_K}=1$ if $K=K^+$
and $\jump{\chi_K}=-1$ if $K=K^-$; in the first case $\bm n=\bm n^K$ and
$\fq=\fq^{\,K}$, in the second case $\bm n=-\bm n^K$ and
$\fq=-\fq^{\,K}$. In both cases $\fq\jump{\chi_K}=\fq^{\,K}$, which gives
\eqref{eq:local-conserve-h}.

\smallskip\noindent\emph{Momentum.} Take $\widehat\zeta^{x}_{\Delta,j}=\chi_K$
in \eqref{eq:dg-qx}. The transport volume terms vanish, and the same
orientation argument as above shows that the face term with $\fqu$ equals
$\int_{\partial K}(\fqu^{\,K})_j\dd s$. For the gravity face term we use
$\avg{\chi_K\hD}=\frac12\hD^{\mathrm{in}}$ on $\partial K$ and the fact that
$\jump{\hD+\BD}\,n_x=(\hD^{\mathrm{in}}+\BD^{\mathrm{in}}-\hD^{\mathrm{ext}}-\BD^{\mathrm{ext}})\,n^K_x$
for both orientations of the face, so that
\[
 -\sum_{e\subset\partial K}\int_e g\bigl(\Pmat(\jump{\hD+\BD})\avg{\chi_K\hD}\bigr)_jn_x\dd s
 =\frac12g\int_{\partial K}\bigl(\Pmat(\jump{\hD+\BD}_K)\hD^{\mathrm{in}}\bigr)_jn^K_x\dd s .
\]
Hence \eqref{eq:dg-qx} with this test function reads
\begin{equation}\label{eq:local-conserve-step1}
\begin{aligned}
 \frac{\dd}{\dd t}\int_K\widehat q^{x}_{\Delta,j}\dd\bm x
 ={}&-\int_{\partial K}(\fqu^{\,K})_j\dd s
   -\int_K g\bigl(\Pmat(\hD)(\hD+\BD)_x\bigr)_j\dd\bm x\\
 &-\frac12g\int_{\partial K}\bigl(\Pmat(\jump{\hD+\BD}_K)\hD^{\mathrm{in}}\bigr)_jn^K_x\dd s\\
 &+\frac12\int_K\bigl(\Pmat((\hD)_t)\uD\bigr)_j\dd\bm x
   -\frac12\int_K\bigl(\Pmat(\qxD)(\uD)_x+\Pmat(\qyD)(\uD)_y\bigr)_j\dd\bm x\\
 &+\frac12\sum_{e\subset\partial K}\int_e\bigl(\Pmat(\fq)\jump{\chi_K\uD}\bigr)_j\dd s .
\end{aligned}
\end{equation}
To remove the skew-symmetric terms in the last two lines we use the height
equation \eqref{eq:dg-h} with the vector-valued test function
\begin{equation}\label{eq:local-conserve-testfunction}
 \bm\psi^{(j)}:=\frac12\Pmat(\uD)\bm\delta_j\,\chi_K ,
\end{equation}
which belongs to $[\Vk]^{\kpce}$ because $\Pmat(\uD)\bm\delta_j$ is a polynomial
of degree at most $k$ on $K$. By the triple-product symmetry
\eqref{eq:P-triple},
\[
 (\hD)_t\cdot\Pmat(\uD)\bm\delta_j=\bigl(\Pmat((\hD)_t)\uD\bigr)_j,
 \qquad
 \qxD\cdot\bigl(\Pmat(\uD)\bm\delta_j\bigr)_x=\bigl(\Pmat(\qxD)(\uD)_x\bigr)_j,
 \qquad
 \fq\cdot\Pmat(\uD^{\pm})\bm\delta_j=\bigl(\Pmat(\fq)\uD^{\pm}\bigr)_j,
\]
and analogously for the $y$-derivative, so that \eqref{eq:dg-h} with
$\bm\psi=\bm\psi^{(j)}$ gives
\[
 \frac12\int_K\bigl(\Pmat((\hD)_t)\uD\bigr)_j\dd\bm x
 -\frac12\int_K\bigl(\Pmat(\qxD)(\uD)_x+\Pmat(\qyD)(\uD)_y\bigr)_j\dd\bm x
 +\frac12\sum_{e\subset\partial K}\int_e\bigl(\Pmat(\fq)\jump{\chi_K\uD}\bigr)_j\dd s=0 .
\]
These are exactly the last two lines of \eqref{eq:local-conserve-step1},
which therefore vanish. It remains to split the gravity terms into their
water-height and bottom parts. Since $\Pmat(\widehat a)\widehat b=\Pmat(\widehat b)\widehat a$,
we have $\partial_x\bigl(\Pmat(\hD)\hD\bigr)=2\Pmat(\hD)(\hD)_x$ on $K$, and
the divergence theorem on $K$ gives
\[
 \int_K g\bigl(\Pmat(\hD)(\hD)_x\bigr)_j\dd\bm x
 =\frac12g\int_{\partial K}\bigl(\Pmat(\hD^{\mathrm{in}})\hD^{\mathrm{in}}\bigr)_jn^K_x\dd s .
\]
Combining this boundary term with the part of the jump term involving
$\jump{\hD}_K$ yields
\[
 -\frac12g\int_{\partial K}\Bigl(\Pmat(\hD^{\mathrm{in}})\hD^{\mathrm{in}}
 +\Pmat(\hD^{\mathrm{ext}}-\hD^{\mathrm{in}})\hD^{\mathrm{in}}\Bigr)_jn^K_x\dd s
 =-\frac12g\int_{\partial K}\bigl(\Pmat(\hD^{\mathrm{ext}})\hD^{\mathrm{in}}\bigr)_jn^K_x\dd s,
\]
while the remaining bottom parts are
$-\int_Kg(\Pmat(\hD)(\BD)_x)_j\dd\bm x-\frac12g\int_{\partial K}(\Pmat(\jump{\BD}_K)\hD^{\mathrm{in}})_jn^K_x\dd s$,
and $\Pmat(\jump{\BD}_K)\hD^{\mathrm{in}}=\Pmat(\hD^{\mathrm{in}})\jump{\BD}_K$.
This proves the first identity in \eqref{eq:local-conserve-q}; the second
follows in the same way from \eqref{eq:dg-qy} with
$\widehat\zeta^{y}_{\Delta,j}=\chi_K$.
\end{proof}

\begin{remark}
The momentum flux $\fqu^{\,K}+\frac12g\Pmat(\hD^{\mathrm{ext}})\hD^{\mathrm{in}}n^K_x$
in \eqref{eq:local-conserve-q} is single valued on every interior face:
$\fqu^{\,K'}=-\fqu^{\,K}$ on a face shared by $K$ and $K'$, and
$\Pmat(\hD^{\mathrm{ext}})\hD^{\mathrm{in}}=\Pmat(\hD^{\mathrm{in}})\hD^{\mathrm{ext}}$
by \eqref{eq:P-comm} while $\bm n^{K'}=-\bm n^K$. Its second part is the
discrete counterpart of the SG pressure flux $\frac12g\Pmat(\hhat)\hhat$ in
\eqref{eq:sg-fluxes}, and the remaining terms are the discretization of the
bottom source $-g\Pmat(\hhat)\nabla\Bhat$, including the jump of the discrete
bottom across $\partial K$. For $\kpce=1$,
\eqref{eq:local-conserve-h}--\eqref{eq:local-conserve-q} reduce to the local
conservation identities (10a)--(10b) of \cite{fu2022high}, in which the jump
of the bottom has to be taken with respect to the outward normal $\bm n^K$,
i.e., as $b^{\mathrm{in}}-b^{\mathrm{ext}}=-\jump{b}_K$.
\end{remark}

Many geophysical flows of interest, such as small perturbations of a lake
at rest, are small deviations from an equilibrium of the SWE, and a scheme
that does not preserve the equilibrium exactly generates spurious waves that
can dominate the true dynamics. In the stochastic setting the relevant
equilibrium is the following.

\begin{definition}[Well-balanced DG--SG scheme for the two-dimensional stochastic SWE]\label{def:wb-2d}
The scheme \eqref{eq:dg-h}, \eqref{eq:dg-qx}, \eqref{eq:dg-qy} is said to be
well-balanced if, started from a discrete stochastic ``lake-at-rest'' state,
its solution $(h_{\Delta,\Lambda}, q^x_{\Delta,\Lambda}, q^y_{\Delta,\Lambda})$, with
$h_{\Delta,\Lambda}:=\sum_{j=1}^{\kpce}\widehat h_{\Delta,j}\phi_j$ and similarly
for the discharges, remains in that state,
\begin{align}\label{well-balance-2d}
    q^x_{\Delta,\Lambda}(x,y,t,\bm\xi) &\equiv 0, \qquad
    q^y_{\Delta,\Lambda}(x,y,t,\bm\xi) \equiv 0, \qquad
    h_{\Delta,\Lambda}(x,y,t,\bm\xi) + B_{\Delta,\Lambda}(x,y,\bm\xi) \equiv C(\bm\xi),
\end{align}
where $C(\bm\xi)$ is a random scalar depending only on $\bm\xi$.
Such a state describes still water whose flat surface elevation is
random. In terms of the PCE coefficients, \eqref{well-balance-2d} is
equivalent to
\begin{align}\label{well-balance-vector-2d}
    \widehat{q}^x_{\Delta} &\equiv \bm 0, \qquad
    \widehat{q}^y_{\Delta} \equiv \bm 0, \qquad
    \widehat{h}_{\Delta} + \widehat{B}_{\Delta} \equiv \widehat{C},
    \qquad \forall (x,y,t) \in \Omega \times [0,T],
\end{align}
with spatial domain $\Omega$ and terminal time $T$, where $\widehat C$ is a
constant vector. The vector equations \eqref{well-balance-vector-2d}
represent a steady state of the skew-symmetric SG system \eqref{eq:sg-skew}.
\end{definition}

\begin{proposition}[Well-balanced DG--SG scheme]\label{prop:wb}
Suppose that the bottom topography $B$ is independent of time and that the
discrete hyperbolicity condition $\Pmat(\hD)\succ0$ holds on every element
and face. Then, for any $\alpha_1\ge0$, the semi-discrete DG--SG scheme
\eqref{eq:dg-h}, \eqref{eq:dg-qx}, \eqref{eq:dg-qy} with the velocity
projection \eqref{eq:velocity-projection} is well-balanced in the sense of
Definition~\ref{def:wb-2d}.
\end{proposition}

\begin{proof}
It suffices to show that, for the initial data
\begin{align}\label{wb_initial}
    \uD = \vD &\equiv \bm 0,
    & \hD + \BD &\equiv \widehat C,
    & \forall (x,y) &\in \Omega,\quad t=0,
\end{align}
the time derivatives of all PCE coefficients vanish, i.e.,
\begin{align}
    (\widehat{h}_{\Delta, j})_t = 0,
\qquad
(\widehat{q}^x_{\Delta, j})_t = 0, \qquad (\widehat{q}^y_{\Delta, j})_t = 0,
\qquad \forall j=1,\ldots,\kpce.
\end{align}
Under the condition \eqref{wb_initial} the discharges vanish by
\eqref{eq:velocity-projection}, so that in \eqref{eq:dg-h} all volume terms
and the central part of $\fq$ vanish, and the dissipative part of $\fq$
vanishes because $\jump{\hD+\BD}=0$ on every face; hence
$(\widehat{h}_{\Delta, j})_t = 0$ for all $j=1,\ldots,\kpce$. In
\eqref{eq:dg-qx} the transport terms and the skew-symmetric corrections
vanish because $\uD=\vD=\bm 0$, $\qxD=\qyD=\bm0$ and $(\hD)_t=\bm 0$; the
pressure--source volume term vanishes because $(\hD+\BD)_x=\bm0$ on every
element; and both face terms, which are proportional to
$\jump{\hD+\BD}$ and to $\fqu$, vanish because the water surface is
continuous and the velocities are zero. Thus
$(\widehat{q}^x_{\Delta, j})_t = 0$, and in the same way
$(\widehat{q}^y_{\Delta, j})_t = 0$, for all $j=1,\ldots,\kpce$. The
semi-discrete solution therefore stays in the discrete steady state, which
proves the well-balanced property \eqref{well-balance-vector-2d}.
\end{proof}

\subsection{Entropy conservation and stability}

To ensure the physical admissibility of nonlinear weak solutions, the
scheme should also satisfy a discrete version of the companion entropy law
\eqref{eq:continuum-entropy-law}, in the sense of the entropy stability
theory of \cite{tadmor1987numerical,tadmor2003entropy,fjordholm2012arbitrarily}. Define the discrete energy
\begin{equation}\label{eq:discrete-energy}
  E_\Delta
  :=\sum_{K\in\Th}\int_K
  \left(
    \frac12 g\norml{\hD}^{2}
    +g\hD\transpose\BD
    +\frac12(\qxD)\transpose\uD
    +\frac12(\qyD)\transpose\vD
  \right)\dd\bm x,
\end{equation}
which is the integral over $\Omega$ of the SG entropy density
\eqref{eq:sg-entropy-density} evaluated at the discrete solution.

\begin{definition}[Entropy conservative/stable DG--SG scheme]\label{def:ec-es}
The semi-discrete DG--SG solution defined by \eqref{eq:dg-h},
\eqref{eq:dg-qx}, \eqref{eq:dg-qy} and \eqref{eq:velocity-projection}, with
given numerical fluxes, is said to be entropy conservative (EC) or entropy
stable (ES) if, under periodic boundary conditions,
\[
  \frac{\dd}{\dd t}E_\Delta=0 \quad \text{(EC)},
  \qquad\text{or}\qquad
  \frac{\dd}{\dd t}E_\Delta\le 0 \quad \text{(ES)},
\]
respectively.
\end{definition}

Periodic boundary conditions are assumed in this subsection so that every
face is an interior face and no boundary contributions appear; the
treatment of physical boundaries, and the corresponding boundary energy
flux, are discussed in \cref{sec:augmented-entropy}. The main result of this
section is the following theorem.

\begin{theorem}[Entropy conservative/stable DG--SG scheme]\label{thm:ec-es}
Suppose that the bottom topography $B$ is independent of time and that the
discrete hyperbolicity condition $\Pmat(\hD)\succ0$ holds on every element
and face. Let the dissipation coefficient $\alpha_1\ge0$ in the numerical
fluxes \eqref{eq:flux-h}, \eqref{eq:flux-qu}, \eqref{eq:flux-qv} be a
single-valued function on the faces of $\Th$. The semi-discrete DG--SG
scheme \eqref{eq:dg-h}, \eqref{eq:dg-qx}, \eqref{eq:dg-qy} with the velocity
projection \eqref{eq:velocity-projection} for the skew-symmetric SG system
\eqref{eq:sg-skew} is
\begin{itemize}[leftmargin=2.2em]
\item entropy conservative with the central fluxes, i.e., for
  $\alpha_1\equiv0$;
\item entropy stable with the Lax--Friedrichs-type numerical fluxes, i.e.,
  for any such $\alpha_1\ge 0$, provided that
  $\Pmat(\avg{\hD+\BD})\succeq0$ on every face.
\end{itemize}
\end{theorem}

 Note that only $\nabla B$ enters the continuous
system, while the water surface $h+B$ enters the discretization, so the
vertical datum of the bottom can
always be chosen such that $B_{\Delta,\Lambda}\coloneqq \sum_{j=1}^{\kpce}\widehat B_{\Delta,j}\phi_j$ is
nonnegative for almost every $\bm\xi$, which is sufficient for the condition $\Pmat(\avg{\hD+\BD})\succeq0$ since $\Pmat$ is linear and
$\Pmat(\avg{\hD})\succ0$. The current dissipation depends on the water surface $h+B$, hence replacing $B$ by $B+c$ with a constant $c$ will change the damping of the velocity jumps; a
variant of the momentum fluxes that does not depend on the datum is given in
Remark~\ref{rem:datum-free-flux}.

Before proving the theorem, we record two elementary identities.

\begin{lemma}[Symmetry of the quadratic form]\label{lem:comm-dot}
For any symmetric matrix $Z\in\R^{\kpce\times\kpce}$ and any vectors $\bm b,\bm c\in\R^{\kpce}$,
\begin{equation}\label{eq:lemma1}
  Z\bm b\cdot \bm c=Z\bm c\cdot \bm b.
\end{equation}
\end{lemma}

\begin{proof}
Since $Z=Z\transpose$,
\[
  Z\bm b\cdot \bm c=\bm c\transpose Z\bm b
  =\bm b\transpose Z\transpose \bm c
  =\bm b\transpose Z\bm c
  =Z\bm c\cdot \bm b.
\]
\end{proof}

\begin{lemma}[Interface identity for $\Pmat$]\label{lem:P-jump-average}
For any coefficient vector $\widehat a$ and any traces $\widehat b^+,\widehat b^-\in\R^{\kpce}$,
\begin{equation}\label{eq:lemma2}
  \frac12\left(
    \Pmat(\widehat a)\widehat b^+\cdot \widehat b^+
    -\Pmat(\widehat a)\widehat b^-\cdot \widehat b^-
  \right)
  =\Pmat(\widehat a)\avg{\widehat b}\cdot\jump{\widehat b}.
\end{equation}
\end{lemma}

\begin{proof}
Write $\widehat b^\pm=\avg{\widehat b}\pm\frac12\jump{\widehat b}$ and expand
the two quadratic forms; the terms quadratic in $\avg{\widehat b}$ and in
$\jump{\widehat b}$ cancel, and the mixed terms combine, by
Lemma~\ref{lem:comm-dot} applied to the symmetric matrix $\Pmat(\widehat a)$,
to the right-hand side of \eqref{eq:lemma2}.
\end{proof}

\begin{proof}[Outline of the proof of \cref{thm:ec-es}]
The complete proof is given in \cref{sec:appendix-proof}; here we outline its
steps. Following the modified entropy variable
\eqref{eq:sg-modified-entropy-variable}, we take the test functions
$g(\hD+\BD)$, $\uD$ and $\vD$ in \eqref{eq:dg-h}, \eqref{eq:dg-qx} and
\eqref{eq:dg-qy}, respectively, and add the results. By the velocity
projection \eqref{eq:velocity-projection} and Lemma~\ref{lem:comm-dot}, the
time-derivative terms combine to $\frac{\dd}{\dd t}E_\Delta$. The volume
integrals of the transport terms cancel against the skew-symmetric
corrections by the product rule and Lemma~\ref{lem:comm-dot}, separately for
each flux direction and each velocity component, so that the asymmetry
\eqref{eq:mixed-flux-representations} of the mixed fluxes does not interfere;
the volume integrals of the gravitational terms cancel by
\eqref{eq:velocity-projection} and \eqref{eq:P-comm}. Inserting the numerical
fluxes \eqref{eq:flux-h}, \eqref{eq:flux-qu}, \eqref{eq:flux-qv} and using
Lemma~\ref{lem:P-jump-average}, the remaining face terms reduce to
\begin{equation}\label{eq:entropy-rate}
\begin{aligned}
  \frac{\dd}{\dd t}E_\Delta
  ={}&-\frac12\sum_{e\in\Ehi}\int_e\alpha_1
  \bigl(\Pmat(\avg{\hD+\BD})\jump{\uD}\cdot\jump{\uD}
  +\Pmat(\avg{\hD+\BD})\jump{\vD}\cdot\jump{\vD}\bigr)\dd s\\
  &-\frac g2\sum_{e\in\Ehi}\int_e\alpha_1\norml{\jump{\hD+\BD}}^2\dd s .
\end{aligned}
\end{equation}
The right-hand side vanishes for $\alpha_1\equiv0$, and it is nonpositive for
$\alpha_1\ge0$.
\end{proof}

\begin{remark}[On the proof]
The proof uses only three ingredients: the symmetry of the triple products
through Lemmas~\ref{lem:comm-dot} and~\ref{lem:P-jump-average}, the fact
that the modified entropy variables $g(\hD+\BD)$, $\uD$, $\vD$ are
admissible test functions, and the elementwise velocity projection
\eqref{eq:velocity-projection}. No quadrature-based or nodal
summation-by-parts structure is required, and the argument applies verbatim
to any polynomial degree $k$ and any orthonormal PCE basis, including
tensor-product bases for several random variables.
The quadrature rules used in the implementation, which integrate the
polynomial integrands of the scheme exactly, are stated in
\cref{sec:appendix-discretization}. The well-balanced
property (Proposition~\ref{prop:wb}) holds for the stochastic lake at rest with an arbitrary random
bottom $B_\Lambda$, since the scheme is written in terms of the water surface
$\hD+\BD$ and the initial data are prepared so that $\hD+\BD$ is exactly
constant; see \cref{sec:numerics} for numerical examples showing the well-balanced preservation of stochastic lake-at-rest states.
\end{remark}

\begin{remark}[An alternative datum-independent momentum flux]\label{rem:datum-free-flux}
By \eqref{eq:jump-PH-u-identity} and the linearity of $\Pmat$, the momentum
flux \eqref{eq:flux-qu} equals
$\Pmat(\fq)\avg{\uD}+\frac12\alpha_1\Pmat(\avg{\hD+\BD})\jump{\uD}$. Its
dependence on the vertical datum, and with it the condition on
$\Pmat(\avg{\hD+\BD})$ in \cref{thm:ec-es}, is removed by weighting the
velocity jump with the water height instead of the water surface,
\begin{equation}\label{eq:flux-qu-datum-free}
  \fqu^{\,\ast}
  :=\Pmat(\fq)\avg{\uD}+\frac12\alpha_1\Pmat(\avg{\hD})\jump{\uD},
\end{equation}
and analogously for $\fqv^{\,\ast}$. These fluxes contain the bottom only
through $\jump{\hD+\BD}$, change sign with the orientation of the face, and
vanish for $\uD=\vD=\bm 0$, so that \cref{thm:local-conserve} and
Proposition~\ref{prop:wb} remain valid. Since \eqref{eq:R1-interface} holds for any
single-valued momentum fluxes, inserting \eqref{eq:flux-qu-datum-free} gives
\[
  R_1=-\frac12\sum_{e\in\Ehi}\int_e\alpha_1
  \bigl(\Pmat(\avg{\hD})\jump{\uD}\cdot\jump{\uD}
  +\Pmat(\avg{\hD})\jump{\vD}\cdot\jump{\vD}\bigr)\dd s\le0
\]
by the discrete hyperbolicity condition alone, while $R_2$ is unchanged;
hence \cref{thm:ec-es} holds for $\fqu^{\,\ast}$ and $\fqv^{\,\ast}$ without
a condition on the bottom. All numerical experiments in \cref{sec:numerics}
use \eqref{eq:flux-qu} and \eqref{eq:flux-qv}.
\end{remark}

\section{Algorithmic details and fully discrete scheme}
\label{sec:implementation}

In this section we summarize the components that turn the semi-discrete
scheme of \cref{sec:dg-scheme} into a robust fully discrete algorithm. The
details of each component are given in \cref{sec:appendix-algorithm}. The
complete procedure is listed in \cref{alg:fully-discrete}, where the labels
on the right refer to the parts of \cref{sec:appendix-algorithm} that
describe the corresponding step.
\begin{itemize}[leftmargin=2em,itemsep=1pt,topsep=3pt]
\item \emph{Desingularized velocity projection} (\cref{sec:desingularization}).
  In nearly dry regions the velocity projection \eqref{eq:velocity-projection}
  becomes ill-conditioned. Following \cite{kurganov2007second}, the small
  eigenvalues of $\Pmat(\hD)$ are corrected before the velocities are
  projected, and the discharges are recomputed from the desingularized
  velocities.
\item \emph{Hyperbolicity-preserving criterion} (\cref{sec_hpc}). The
  condition $\Pmat(\hD)\succ0$ is ensured by the positivity of the truncated
  water height at the nodes of a stochastic quadrature rule
  \cite[Theorem~3.2]{dai2021hyperbolicity}, which is preserved in time by a
  time step restriction together with a linear scaling around the cell
  averages.
\item \emph{Troubled-cell indicator and limiter} (\cref{sec_slopelimiting}).
  Troubled cells are identified by the indicator of \cite{fu2017new} applied to
  the mean water surface, and the characteristic TVB limiter of
  \cite{cockburn1998runge} is applied there to the equilibrium variables
  $(\hD+\BD,\qxD,\qyD)$, which leaves the stochastic lake at rest unchanged.
\item \emph{Boundary conditions} (\cref{sec:boundary}). Every boundary face
  carries the same numerical fluxes as an interior face, with the exterior
  trace replaced by a ghost state; periodic, reflecting wall, zeroth-order
  outflow and characteristic outflow conditions are supported.
\item \emph{Time integration} (\cref{sec:time-stepping}). We use the
  third-order SSP Runge--Kutta method with an adaptive time step: the
  admissible step is recomputed at every intermediate stage, and the step is
  restarted with the smaller value if necessary.
\item \emph{Augmented energy} (\cref{sec:augmented-entropy}). The augmented
  energy adds the time-integrated boundary energy flux back to the energy,
  which separates the dissipation of the scheme from the exchange of energy
  through physical boundaries; its monotonicity is the property that we
  monitor numerically.
\end{itemize}

\begin{algorithm}[H]
\caption{Fully discrete two-dimensional DG--SG scheme with SSP-RK3}
\label{alg:fully-discrete}
\begin{algorithmic}[1]
\Require Initial data and bottom topography; mesh and face connectivity;
PCE basis and triple products; final time $T$.
\State \textbf{Initialize (SG):} compute the PCE coefficients
$\hhat,\qxhat,\qyhat,\Bhat$ of the initial data and of the bottom.
\State \textbf{Initialize (DG):} project the coefficients onto
$[\Vk]^{\kpce}$; apply the limiter and the linear scaling to the
initial state (preserving $\hD+\BD$ if the initial state is at rest); project
both velocity components by \eqref{vel_update}--\eqref{eq:desing-recompute}.
\Comment{\ref{sec:desingularization}, \ref{sec_hpc}, \ref{sec_slopelimiting}}
\State Set $t\gets0$.
\While{$t<T$}
  \State Compute $\Delta t$ from the CFL condition and the
  hyperbolicity-preserving bound \eqref{eq:hp-dt-mean} at $U^n$.
  \Comment{\ref{sec:time-stepping}}
  \For{each SSP-RK3 stage $r=1,2,3$}
    \State Assemble the volume and face residuals $\mathcal L(U^{(r-1)})$
    with the projected velocities and the ghost states.
    \Comment{\ref{sec:boundary}}
    \State Form the stage state $U^{(r)}$ by \eqref{eq:ssprk3}.
    \Comment{\ref{sec:time-stepping}}
    \State Identify the troubled cells and apply the characteristic TVB
    limiter to $(\hD+\BD,\qxD,\qyD)$.
    \Comment{\ref{sec_slopelimiting}}
    \State Apply the linear scaling \eqref{eq:hp-linear-scaling} and
    reproject both velocity components.
    \Comment{\ref{sec:desingularization}, \ref{sec_hpc}}
    \State Recompute the admissible step $\Delta t_r$ at $U^{(r)}$.
    \Comment{\ref{sec:time-stepping}}
    \If{$\Delta t_r<\Delta t$}
      \State Set $\Delta t\gets\Delta t_r$ and restart from stage 1.
    \EndIf
  \EndFor
  \State Accept $U^{n+1}=U^{(3)}$, update the energy and the augmented
  energy, and set $t\gets t+\Delta t$.
  \Comment{\ref{sec:augmented-entropy}}
\EndWhile
\end{algorithmic}
\end{algorithm}

\section{Numerical experiments}
\label{sec:numerics}

In this section we assess the performance of the proposed DG--SG scheme
through a sequence of numerical experiments. The examples verify the
properties established in the previous sections, including high-order
accuracy, the well-balanced preservation of stochastic lake-at-rest states,
entropy stability, robustness in nearly dry and dry
regions and on domains with internal walls, and applicability to
multivariate random inputs. We present nine test problems: two accuracy
tests (\cref{sec:num-accuracy}), the stochastic two-dimensional examples
used for the finite-volume SG schemes of \cite{epshteyn2025energy}
(\cref{sec:num-hump,sec:num-perturbation,sec:num-plateau,sec:num-2rv}),
two dam-break problems with uncertain water height, the circular dam break
of \cite{seaid2004non} (\cref{sec:num-circular}) and the curved partial dam
break of \cite{wu2021high} (\cref{sec:num-curved}), and a flow over a
random bottom field in ten random variables (\cref{sec:num-nrv}).
The setup of all test problems is summarized in \cref{tab:num-setup}.

\begin{table}[htbp]
\centering
\captionsetup{font=scriptsize}
\caption{Setup of the test problems: problem data (top) and discretization
(bottom). ``Split'' means that every cell of a structured grid is divided
into two triangles. The limiter is the characteristic TVB limiter of
\cref{sec_slopelimiting} with the listed threshold $\mathrm{tol}$. The
velocity is obtained from the desingularized projection \eqref{vel_update}
with the listed $\epsilon$, or from the exact weak projection
\eqref{eq:velocity-projection}. The last column refers to the
hyperbolicity-preserving scaling of \cref{sec_hpc}.}
\label{tab:num-setup}
\footnotesize
\setlength{\tabcolsep}{3.5pt}
\begin{tabular}{@{}>{\raggedright}p{0.218\textwidth}lc>{\raggedright}p{0.19\textwidth}>{\raggedright}p{0.13\textwidth}ll@{}}
\toprule
Test problem & $\Omega$ & $g$ & Random input & Left and right boundaries & Top and bottom & $T$\\
\midrule
\ref{sec:num-accuracy-cmame} Stochastic hump flow (accuracy) & $[0,2]\times[0,1]$ & $1$ & $\xi\sim\mathcal U(-1,1)$ & zeroth-order outflow & periodic & $0.07$\\
\addlinespace[2pt]
\ref{sec:num-accuracy-fu} Manufactured solution (accuracy) & $[0,1]^2$ & $9.812$ & $\xi\sim\mathcal U(-1,1)$ & periodic & periodic & $0.004$\\
\addlinespace[2pt]
\ref{sec:num-hump} Gaussian-shape hump & $[0,2]\times[0,1]$ & $1$ & $\xi\sim\mathcal U(-1,1)$ & characteristic outflow & wall & $1.8$\\
\addlinespace[2pt]
\ref{sec:num-perturbation} Perturbation of a lake at rest, cases 1 and 2 & $[-1,1]^2$ & $1$ & $\xi\sim\mathcal U(-1,1)$ & characteristic outflow & wall & $0.8$\\
\addlinespace[2pt]
\ref{sec:num-plateau} Submerged flat plateau & $[-0.5,0.5]^2$ & $1$ & $\xi\sim\mathcal U(-1,1)$ & characteristic outflow & wall & $0.65$\\
\addlinespace[2pt]
\ref{sec:num-circular} Circular dam break, wet and dry bed & $[0,50]^2$ & $9.812$ & $\xi\sim\mathcal U(-1,1)$ & wall & wall & $0.69$\\
\addlinespace[2pt]
\ref{sec:num-curved} Curved partial dam break & $[-10,10]^2$, $[-5,5]^2$ & $1$ & $\xi\sim\mathcal U(-1,1)$ & wall (also on both sides of the dam) & wall & $1.5$\\
\addlinespace[2pt]
\ref{sec:num-2rv} Hump, two random variables & $[0,2]\times[0,1]$ & $1$ & $\xi_1\sim\rho_1\propto(1-\xi_1)^2$, $\xi_2\sim\mathcal U(-1,1)$ & zeroth-order outflow & wall & $0.07$\\
\addlinespace[2pt]
\ref{sec:num-nrv} Hump, ten random variables & $[0,2]\times[0,1]$ & $1$ & $\xi_i=0.1\,\eta_i$, $\eta_i\sim\rho_i\propto(1-\eta_i^2)^{i^2}$, $i=1,\dots,10$ & characteristic, reference state \eqref{eq:num-hump2rv-ic} & periodic & $0.07$\\
\bottomrule
\end{tabular}

\medskip
\begin{tabular}{@{}l>{\raggedright}p{0.312\textwidth}ccl>{\raggedright}p{0.17\textwidth}ll@{}}
\toprule
 & Mesh & $k$ & $\kpce$ & $\Ccfl$ & Limiter & Velocity & Scaling\\
\midrule
\ref{sec:num-accuracy-cmame} & $N_x\times N_y$ rectangles, split; $N_x=N_y=16,32,64$; reference: $128\times128$ & $1,2,3$ & $4$ & $0.05$ ($k=3$: $0.02$) & none & exact & off\\
\addlinespace[2pt]
\ref{sec:num-accuracy-fu} & $N_x\times N_y$ squares, split; $N_x=N_y=25,50,100$ & $1,2,3$ & $4$ & $0.05$ & none & exact & off\\
\addlinespace[2pt]
\ref{sec:num-hump} & DistMesh, $\Delta_0=0.02$ (11345 triangles) & $2$ & $2$ & $0.05$ & $\mathrm{tol}=0.02$ & $\epsilon=\hK$ & on\\
\addlinespace[2pt]
\ref{sec:num-perturbation} & DistMesh, $\Delta_0=0.025$ (14632 triangles) & $2$ & $2$ & $0.05$ & $\mathrm{tol}=0.02$ & $\epsilon=\hK$ & on\\
\addlinespace[2pt]
\ref{sec:num-plateau} & DistMesh, $\Delta_0=0.025$ (3638 triangles) & $2$ & $2$ & $0.02$ & $\mathrm{tol}=0.02$ & $\epsilon=\hK$ & on\\
\addlinespace[2pt]
\ref{sec:num-circular} & $50\times50$ squares, split ($\Delta x=1$, 5\,000 triangles) & $2$ & $2$ & $0.05$ & $\mathrm{tol}=0.02$; dry-cell limiter, $\varepsilon_{\mathrm{dry}}=5\times10^{-3}$ & $\epsilon=10^{-2}$ & on\\
\addlinespace[2pt]
\ref{sec:num-curved} & DistMesh, $\Delta_0=1$ (938 triangles); $20\times20$ quadrilaterals, split (800 triangles) & $3$ & $2$ & $0.05$ & none & $\epsilon=10^{-2}$ & on\\
\addlinespace[2pt]
\ref{sec:num-2rv} & DistMesh, $\Delta_0=0.06$ & $2$ & $9$ & $0.05$ & none & exact & off\\
\addlinespace[2pt]
\ref{sec:num-nrv} & $64\times64$ rectangles, split ($8192$ triangles, $\max_K\hK=0.014$) & $2$ & $21$ & $0.05$ & none & exact & off\\
\bottomrule
\end{tabular}
\end{table}

Unless stated otherwise, the following setup is used. The random input is a
scalar Beta-distributed random variable supported on $[-1,1]$,
\[
 \xi\sim\operatorname{Beta}_{[-1,1]}(\beta+1,\alpha+1),
 \qquad
 \rho(\xi)\propto(1-\xi)^\alpha(1+\xi)^\beta,
 \qquad \alpha,\beta>-1,
\]
and the associated orthonormal basis consists of Jacobi polynomials; the
special case $\alpha=\beta=0$ is the uniform distribution on $[-1,1]$ with
Legendre polynomials as the basis. In every example with one random
variable, $\xi\sim\mathcal U(-1,1)$; \cref{sec:num-2rv} uses two independent
random variables, one of them with a Jacobi weight, and \cref{sec:num-nrv}
ten independent variables with Jacobi weights. We plot the mean and the standard deviation of
the water surface $w:=h+B$, or of the water height $h$ in the dam-break
examples, where $B=0$.

\paragraph{Discretization.} All examples use an orthonormal PCE basis in
the random variables and the DG--SG scheme of \cref{sec:dg-scheme} with
polynomial degree $k$ and the entropy stable (ES) fluxes, on either an
unstructured DistMesh triangulation \cite{persson2004simple} with target
edge length $\Delta_0$ or a structured grid of split quadrilaterals. Time
stepping is SSP-RK3 with the CFL number $\Ccfl$ and the local mesh size
$\hK$ of \eqref{eq:num-cfl}. Where listed in \cref{tab:num-setup}, the
desingularized velocity projection (\cref{sec:desingularization}), the
hyperbolicity-preserving scaling (\cref{sec_hpc}) and the characteristic TVB
limiter with threshold $\mathrm{tol}$ (\cref{sec_slopelimiting}) are used,
and the initial data are projected so that the lake at rest is represented
exactly. The details common to all examples are given in
\cref{sec:appendix-discretization}.
Compared with \cite{epshteyn2025energy}, which uses finite volumes on a
$200\times200$ Cartesian grid, $\kpce=4$ and periodic top and bottom
boundaries, the runs of \cref{sec:num-hump,sec:num-perturbation,sec:num-plateau}
use unstructured DG meshes, two PCE terms and walls on the top and bottom
boundaries.

\paragraph{Energy diagnostics.} We record the discrete energy $E_\Delta(t)$
of \eqref{eq:discrete-energy} and plot its relative change
$(E_\Delta(t)-E_\Delta(0))/E_\Delta(0)$. In the examples with outflow boundaries
$E_\Delta$ also changes when waves cross the boundary, so $E_\Delta$ alone
does not show the dissipation of the scheme. For the runs of
\cref{sec:num-hump,sec:num-perturbation,sec:num-plateau,sec:num-2rv,sec:num-nrv} we
therefore also record the augmented energy $\widetilde E_\Delta(t)$ of
\eqref{eq:discrete-augmented-entropy}, with the boundary functional
\eqref{eq:numerical-boundary-entropy-flux} evaluated from the same ES
fluxes as on interior faces and the ghost state, and the time integral
accumulated with the SSP-RK3 stage weights $(\frac16,\frac16,\frac23)$. By
\eqref{eq:augmented-entropy-stability}, $\widetilde E_\Delta$ is
non-increasing in the semi-discrete scheme; $E_\Delta$ itself is
non-increasing when all boundaries are periodic or reflecting.

\subsection{Accuracy tests}
\label{sec:num-accuracy}

Two smooth problems are used to measure the order of convergence: the
stochastic hump flow of \cite[Section~6.1]{epshteyn2025energy}, measured
against a fine-mesh reference solution, and a stochastic extension of the
manufactured solution of \cite[Example~4.7(b)]{fu2022high}, measured against
the exact solution. The settings of both tests are listed in
\cref{tab:num-setup}. The meshes
are structured: the domain is divided into $N_x\times N_y$ congruent
rectangles, each split
into two triangles, with the diagonal alternating from cell to cell.

\paragraph{Error norms.} For the PCE coefficient vector
$\widehat{\bm r}(x,y)=(\widehat r_1,\dots,\widehat r_{\kpce})$ of the error $\err$ of a
scalar field, the orthonormality of the basis gives
$\|\err(x,y,\cdot)\|_{L^2_\rho}^2=\sum_{j=1}^{\kpce}\widehat r_j(x,y)^2$, and we report
\begin{equation}\label{eq:num-errnorms}
  \|\err\|_{L^1(\Omega;L^2_\rho)}=\int_\Omega\Big(\sum_{j=1}^{\kpce}\widehat r_j^2\Big)^{1/2}\dd\bm x,
  \qquad
  \|\err\|_{L^2(\Omega;L^2_\rho)}=\Big(\int_\Omega\sum_{j=1}^{\kpce}\widehat r_j^2\dd\bm x\Big)^{1/2}.
\end{equation}
The first norm is the one used in \cite[Eq.~(6.6)]{epshteyn2025energy}; the
second reduces to the $L^2$ error of \cite{fu2022high} in the deterministic
case. For the discharge vector $\bm m=(q^x,q^y)$ and the velocity vector
$\bm u=(u,v)$ the errors of the two components are combined as
\begin{equation}\label{eq:num-errnorms-vector}
  \begin{aligned}
  \|\err_{\bm m}\|_{L^1(\Omega;L^2_\rho)}&=\|\err_{q^x}\|_{L^1(\Omega;L^2_\rho)}+\|\err_{q^y}\|_{L^1(\Omega;L^2_\rho)},\\
  \|\err_{\bm m}\|_{L^2(\Omega;L^2_\rho)}&=\Big(\|\err_{q^x}\|^2_{L^2(\Omega;L^2_\rho)}+\|\err_{q^y}\|^2_{L^2(\Omega;L^2_\rho)}\Big)^{1/2},
  \end{aligned}
\end{equation}
and likewise for $\bm u$ from $\err_u$ and $\err_v$; the $L^2$ combination is the
$L^2(\Omega;L^2_\rho)$ norm of the vector error, as in
\cite[Tables~3--5]{fu2022high}. The integrals over $\Omega$ are evaluated with
the volume quadrature of the test mesh. The reference or exact solution is
evaluated at the quadrature points of the test mesh: from its own DG
polynomial in the triangle that contains the point for a reference solution,
and from the exact PCE coefficients for the manufactured solution. The
velocity error of the manufactured solution is measured against the
projected velocity of the $L^2$-projected exact state. Orders are $\log_2$ of
the ratio of consecutive errors. Besides the three fields we also list the
sum $\err_h+\err_{\bm m}$ of the water-height and discharge errors and its
order.

\subsubsection{Stochastic hump flow}
\label{sec:num-accuracy-cmame}

The initial data are the deterministic uniform flow
\cite[Eq.~(6.4)]{epshteyn2025energy}
\begin{equation}\label{eq:num-acc-cmame-ic}
  w(x,y,0,\xi)=1,\qquad u(x,y,0,\xi)=0.3,\qquad v(x,y,0,\xi)=0,
\end{equation}
over the stochastic elliptic hump \cite[Eq.~(6.5)]{epshteyn2025energy}
\begin{equation}\label{eq:num-acc-cmame-bottom}
  B(x,y,\xi)=0.5\exp\!\big(-25(x-1)^2-50(y-0.5)^2\big)+0.1(\xi+1),
\end{equation}
with the boundary conditions of \cite{epshteyn2025energy}. For every
degree $k$, the solution is compared with a reference solution computed with
the same scheme and the same polynomial degree on the $128\times128$ mesh, so that the
tables measure the self-convergence of each discretization.

\Cref{tab:num-acc-cmame-L1,tab:num-acc-cmame-L2} report the errors and
orders of convergence in the two norms. The water height converges with the expected order $k+1$ for all
three degrees.
The orders of the discharge and of the velocity are somewhat below $k+1$
on these meshes; their errors are smaller than that of the water height, so
that the combined error $\err_h+\err_{\bm m}$ converges with an order close
to $k+1$ for all three degrees. For reference, the
second-order finite-volume schemes of \cite[Table~1]{epshteyn2025energy}
reach $L^1$ errors of about $10^{-5}$ for $h$ on a $400\times400$ grid,
whereas the $\mathbb P^2$ solution on the $64\times64$ mesh and the
$\mathbb P^3$ solution on the $32\times32$ mesh have errors of
$2\times10^{-6}$ and $10^{-6}$; since the two sets of errors are measured
against different reference solutions, this comparison is only indicative.

\begin{table}[htbp]
\centering
\captionsetup{font=scriptsize}
\caption{Stochastic hump accuracy test: $L^1(\Omega;L^2_\rho)$ errors
\eqref{eq:num-errnorms}--\eqref{eq:num-errnorms-vector} at $T=0.07$ against the
$128\times128$ reference solution, and orders. $\kpce=4$.}
\label{tab:num-acc-cmame-L1}
\small\setlength{\tabcolsep}{3.5pt}
\begin{tabular}{@{}lrrcrcrcrc@{}}
\toprule
& & \multicolumn{2}{c}{$h$} & \multicolumn{2}{c}{$\bm m=(q^x,q^y)$} & \multicolumn{2}{c}{$\bm u=(u,v)$} & \multicolumn{2}{c}{$\err_h+\err_{\bm m}$}\\
\cmidrule(lr){3-4}\cmidrule(lr){5-6}\cmidrule(lr){7-8}\cmidrule(lr){9-10}
& $N_x$ & error & order & error & order & error & order & error & order\\
\midrule
$\mathbb P^1$ & 16 & 1.150e-3 &  & 3.621e-4 &  & 2.514e-4 &  & 1.512e-3 & \\
 & 32 & 2.931e-4 & 1.97 & 9.950e-5 & 1.86 & 8.715e-5 & 1.53 & 3.926e-4 & 1.95\\
 & 64 & 7.300e-5 & 2.01 & 2.677e-5 & 1.89 & 2.519e-5 & 1.79 & 9.977e-5 & 1.98\\
\midrule
$\mathbb P^2$ & 16 & 1.350e-4 &  & 6.613e-5 &  & 7.393e-5 &  & 2.011e-4 & \\
 & 32 & 1.703e-5 & 2.99 & 8.916e-6 & 2.89 & 1.073e-5 & 2.78 & 2.595e-5 & 2.95\\
 & 64 & 2.116e-6 & 3.01 & 1.202e-6 & 2.89 & 1.496e-6 & 2.84 & 3.318e-6 & 2.97\\
\midrule
$\mathbb P^3$ & 16 & 1.463e-5 &  & 6.978e-6 &  & 8.329e-6 &  & 2.161e-5 & \\
 & 32 & 9.417e-7 & 3.96 & 4.965e-7 & 3.81 & 6.136e-7 & 3.76 & 1.438e-6 & 3.91\\
 & 64 & 5.789e-8 & 4.02 & 3.404e-8 & 3.87 & 4.353e-8 & 3.82 & 9.193e-8 & 3.97\\
\bottomrule
\end{tabular}
\end{table}

\begin{table}[htbp]
\centering
\captionsetup{font=scriptsize}
\caption{Stochastic hump accuracy test: $L^2(\Omega;L^2_\rho)$ errors
\eqref{eq:num-errnorms}--\eqref{eq:num-errnorms-vector} at $T=0.07$ against the
$128\times128$ reference solution, and orders. $\kpce=4$.}
\label{tab:num-acc-cmame-L2}
\small\setlength{\tabcolsep}{3.5pt}
\begin{tabular}{@{}lrrcrcrcrc@{}}
\toprule
& & \multicolumn{2}{c}{$h$} & \multicolumn{2}{c}{$\bm m=(q^x,q^y)$} & \multicolumn{2}{c}{$\bm u=(u,v)$} & \multicolumn{2}{c}{$\err_h+\err_{\bm m}$}\\
\cmidrule(lr){3-4}\cmidrule(lr){5-6}\cmidrule(lr){7-8}\cmidrule(lr){9-10}
& $N_x$ & error & order & error & order & error & order & error & order\\
\midrule
$\mathbb P^1$ & 16 & 2.679e-3 &  & 7.619e-4 &  & 4.639e-4 &  & 3.440e-3 & \\
 & 32 & 6.943e-4 & 1.95 & 2.123e-4 & 1.84 & 1.770e-4 & 1.39 & 9.066e-4 & 1.92\\
 & 64 & 1.730e-4 & 2.00 & 5.575e-5 & 1.93 & 5.410e-5 & 1.71 & 2.288e-4 & 1.99\\
\midrule
$\mathbb P^2$ & 16 & 3.352e-4 &  & 1.149e-4 &  & 1.282e-4 &  & 4.501e-4 & \\
 & 32 & 4.253e-5 & 2.98 & 1.525e-5 & 2.91 & 1.976e-5 & 2.70 & 5.778e-5 & 2.96\\
 & 64 & 5.316e-6 & 3.00 & 1.995e-6 & 2.93 & 2.758e-6 & 2.84 & 7.311e-6 & 2.98\\
\midrule
$\mathbb P^3$ & 16 & 3.742e-5 &  & 1.190e-5 &  & 1.308e-5 &  & 4.932e-5 & \\
 & 32 & 2.435e-6 & 3.94 & 8.715e-7 & 3.77 & 1.042e-6 & 3.65 & 3.307e-6 & 3.90\\
 & 64 & 1.515e-7 & 4.01 & 5.915e-8 & 3.88 & 7.656e-8 & 3.77 & 2.106e-7 & 3.97\\
\bottomrule
\end{tabular}
\end{table}

\subsubsection{Manufactured solution with a random water height}
\label{sec:num-accuracy-fu}

The second test uses the method of manufactured solutions with the exact
fields of \cite[Example~4.7(b)]{fu2022high}, to which a random perturbation
of the water height is added:
\begin{equation}\label{eq:num-acc-fu-exact}
  \begin{aligned}
    h(x,y,t,\xi)&=10+\mathrm e^{\sin(2\pi x)}\cos(2\pi y)\cos(2\pi t)+0.1\,\xi,\\
    q^x(x,y,t)&=\sin(\cos(2\pi x))\sin(2\pi y)\cos(2\pi t),\\
    q^y(x,y,t)&=\cos(2\pi x)\cos(\sin(2\pi y))\cos(2\pi t),
  \end{aligned}
  \qquad B(x,y)=\sin(2\pi x)+\sin(2\pi y),
\end{equation}
on the periodic unit square. The source term is the residual of
the continuous SG system,
\begin{equation}\label{eq:num-acc-fu-source}
  S_h=\partial_t\hhat+\nabla\cdot\widehat{\bm q},\qquad
  \bm S_q=\partial_t\widehat{\bm q}+\nabla\cdot\big(\Pmat(\widehat{\bm q})\uhat\big)
  +g\,\Pmat(\hhat)\nabla(\hhat+\Bhat),\qquad
  \Pmat(\hhat)\uhat=\widehat{\bm q},
\end{equation}
evaluated on the exact PCE coefficients \eqref{eq:num-acc-fu-exact} and
$L^2$-projected onto the DG space at every Runge--Kutta stage time, so that
\eqref{eq:num-acc-fu-exact} solves the forced SG system exactly for every
truncation $\kpce$. The initial data are the $L^2$ projection of the exact
state at $t=0$, and the solution is compared with the exact
coefficients at the final time, which is that of \cite{fu2022high}.

\Cref{tab:num-acc-fu-L1,tab:num-acc-fu-L2} report the errors against
the exact coefficients. The DG--SG scheme attains the optimal order $k+1$
in both norms for the water height, the discharge and the velocity and for
all three degrees, in agreement with the deterministic results of
\cite[Tables~4 and~5]{fu2022high}. The discharge errors are about ten
times the velocity errors, since the discharge is the product of the
velocity and the water height, which varies around the depth $10$.

\begin{table}[htbp]
\centering
\captionsetup{font=scriptsize}
\caption{Manufactured solution: $L^1(\Omega;L^2_\rho)$ errors
\eqref{eq:num-errnorms}--\eqref{eq:num-errnorms-vector} at $T=0.004$ against the
exact PCE coefficients, and orders. $\kpce=4$.}
\label{tab:num-acc-fu-L1}
\small\setlength{\tabcolsep}{3.5pt}
\begin{tabular}{@{}lrrcrcrcrc@{}}
\toprule
& & \multicolumn{2}{c}{$h$} & \multicolumn{2}{c}{$\bm m=(q^x,q^y)$} & \multicolumn{2}{c}{$\bm u=(u,v)$} & \multicolumn{2}{c}{$\err_h+\err_{\bm m}$}\\
\cmidrule(lr){3-4}\cmidrule(lr){5-6}\cmidrule(lr){7-8}\cmidrule(lr){9-10}
& $N_x$ & error & order & error & order & error & order & error & order\\
\midrule
$\mathbb P^1$ & 25 & 2.160e-3 &  & 5.036e-2 &  & 4.988e-3 &  & 5.252e-2 & \\
 & 50 & 5.388e-4 & 2.00 & 1.302e-2 & 1.95 & 1.287e-3 & 1.95 & 1.356e-2 & 1.95\\
 & 100 & 1.348e-4 & 2.00 & 3.272e-3 & 1.99 & 3.237e-4 & 1.99 & 3.407e-3 & 1.99\\
\midrule
$\mathbb P^2$ & 25 & 7.215e-5 &  & 1.428e-3 &  & 1.426e-4 &  & 1.500e-3 & \\
 & 50 & 8.981e-6 & 3.01 & 1.818e-4 & 2.97 & 1.798e-5 & 2.99 & 1.908e-4 & 2.98\\
 & 100 & 1.076e-6 & 3.06 & 2.176e-5 & 3.06 & 2.171e-6 & 3.05 & 2.284e-5 & 3.06\\
\midrule
$\mathbb P^3$ & 25 & 2.363e-6 &  & 5.288e-5 &  & 5.314e-6 &  & 5.524e-5 & \\
 & 50 & 1.495e-7 & 3.98 & 3.507e-6 & 3.91 & 3.508e-7 & 3.92 & 3.657e-6 & 3.92\\
 & 100 & 9.128e-9 & 4.03 & 2.219e-7 & 3.98 & 2.232e-8 & 3.97 & 2.310e-7 & 3.98\\
\bottomrule
\end{tabular}
\end{table}

\begin{table}[htbp]
\centering
\captionsetup{font=scriptsize}
\caption{Manufactured solution: $L^2(\Omega;L^2_\rho)$ errors
\eqref{eq:num-errnorms}--\eqref{eq:num-errnorms-vector} at $T=0.004$ against the
exact PCE coefficients, and orders. $\kpce=4$.}
\label{tab:num-acc-fu-L2}
\small\setlength{\tabcolsep}{3.5pt}
\begin{tabular}{@{}lrrcrcrcrc@{}}
\toprule
& & \multicolumn{2}{c}{$h$} & \multicolumn{2}{c}{$\bm m=(q^x,q^y)$} & \multicolumn{2}{c}{$\bm u=(u,v)$} & \multicolumn{2}{c}{$\err_h+\err_{\bm m}$}\\
\cmidrule(lr){3-4}\cmidrule(lr){5-6}\cmidrule(lr){7-8}\cmidrule(lr){9-10}
& $N_x$ & error & order & error & order & error & order & error & order\\
\midrule
$\mathbb P^1$ & 25 & 3.248e-3 &  & 4.964e-2 &  & 4.806e-3 &  & 5.289e-2 & \\
 & 50 & 8.096e-4 & 2.00 & 1.271e-2 & 1.97 & 1.226e-3 & 1.97 & 1.352e-2 & 1.97\\
 & 100 & 2.025e-4 & 2.00 & 3.189e-3 & 1.99 & 3.079e-4 & 1.99 & 3.391e-3 & 2.00\\
\midrule
$\mathbb P^2$ & 25 & 1.228e-4 &  & 1.476e-3 &  & 1.469e-4 &  & 1.599e-3 & \\
 & 50 & 1.510e-5 & 3.02 & 1.867e-4 & 2.98 & 1.820e-5 & 3.01 & 2.018e-4 & 2.99\\
 & 100 & 1.863e-6 & 3.02 & 2.209e-5 & 3.08 & 2.201e-6 & 3.05 & 2.395e-5 & 3.07\\
\midrule
$\mathbb P^3$ & 25 & 4.356e-6 &  & 6.080e-5 &  & 6.126e-6 &  & 6.516e-5 & \\
 & 50 & 2.720e-7 & 4.00 & 4.026e-6 & 3.92 & 4.021e-7 & 3.93 & 4.298e-6 & 3.92\\
 & 100 & 1.682e-8 & 4.01 & 2.545e-7 & 3.98 & 2.562e-8 & 3.97 & 2.713e-7 & 3.99\\
\bottomrule
\end{tabular}
\end{table}

\subsection{Gaussian-shape hump with stochastic bottom}
\label{sec:num-hump}

Following \cite[Section~6.2]{epshteyn2025energy}, the randomness is in the
position of the Gaussian hump,
\begin{equation}\label{eq:num-hump-bottom}
  B(x,y,\xi)=0.8\exp\!\big(-5(x-0.9+0.1\xi)^2-50(y-0.5)^2\big),
\end{equation}
and the initial water surface is the discontinuous, deterministic pulse
\begin{equation}\label{eq:num-hump-ic}
  w(x,y,0,\xi)=\begin{cases}1.01, & 0.05<x<0.15,\\ 1, & \text{otherwise},\end{cases}
  \qquad u(x,y,0)=v(x,y,0)=0.
\end{equation}
The test is a modification of the tests in
\cite{dai2022hyperbolicity,xing2010positivity}; it examines how the
uncertainty located in the topography is transferred to the waves that
interact with it. The parameters are listed in
\cref{tab:num-setup}; the mean
water surface with the standard deviation disks is plotted at the output
times of \cite[Section~6.2]{epshteyn2025energy}.

\begin{figure}[htbp]
\centering
\includegraphics[width=\textwidth,trim={0 0 0 28},clip]{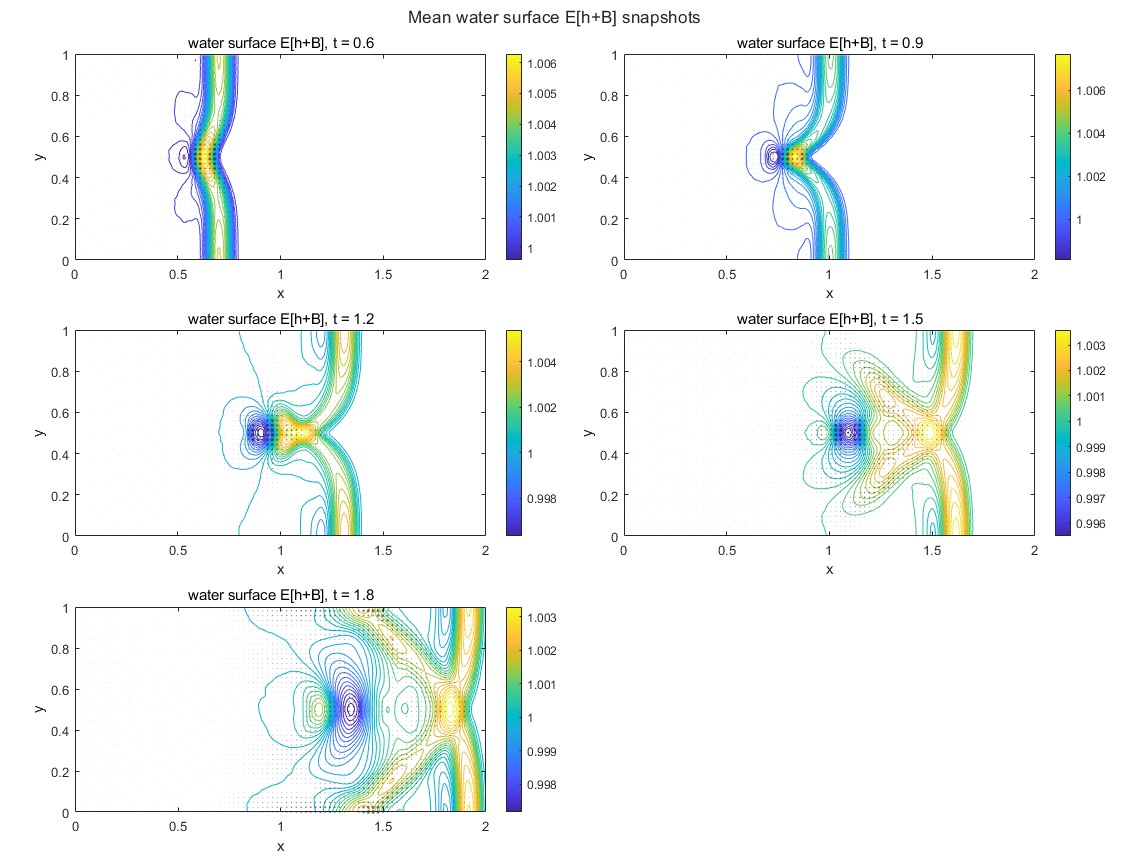}
\captionsetup{font=scriptsize}
\caption{Gaussian-shape hump with stochastic bottom: contours of the mean
water surface $\E[w]$ at $t=0.6$, $0.9$, $1.2$, $1.5$, $1.8$ (left to right,
top to bottom). Disk area is proportional to $\Std[w]$, normalized by its
maximum in each panel; the maximum standard deviation is
$2.9\times10^{-3}$, $5.2\times10^{-3}$, $3.1\times10^{-3}$,
$1.8\times10^{-3}$ and $1.2\times10^{-3}$, respectively. $\mathbb P^2$ DG,
$\kpce=2$, $\Delta_0=0.02$, $\mathrm{tol}=0.02$.}
\label{fig:num-hump-snapshots}
\end{figure}

In \cref{fig:num-hump-snapshots}, the initial pulse splits into a
left-going wave, which leaves the domain through the outflow boundary
$x=0$, and a right-going wave that interacts with the hump: the front slows
down over the shallower water above the hump, the transmitted part bends,
and the reflected waves interact with the walls $y=0,1$ to produce the
patterns at $t=1.5$ and $1.8$. The uncertainty of the water surface is generated at the
hump, where the random position of the bottom enters, when the wave
arrives there; afterwards it spreads along the wave structures and its maximum
decays from $5.2\times10^{-3}$ at $t=0.9$ to $1.2\times10^{-3}$ at $t=1.8$.
The wave structures and the distribution of the uncertainty agree with the
second-order finite-volume results of \cite[Figure~3]{epshteyn2025energy}
on a $200\times200$ grid; the DG--SG solution resolves the fronts more
sharply and retains a larger maximum of $\Std[w]$, consistent with its lower
numerical dissipation (the same trend is seen in \cite{epshteyn2025energy}
between the first- and second-order schemes). Ahead of the right-going
wave, where the bottom is random, but the water is still undisturbed, the
deviation of the mean water surface from the rest state $w=1$ and the
standard deviation $\Std[w]$ remain at round-off level: the
stochastic lake at rest over the random topography is preserved and no
spurious waves are generated by the uncertainty of the bottom, as guaranteed
by Proposition~\ref{prop:wb}.

\begin{figure}[htbp]
\centering
\includegraphics[width=0.92\textwidth]{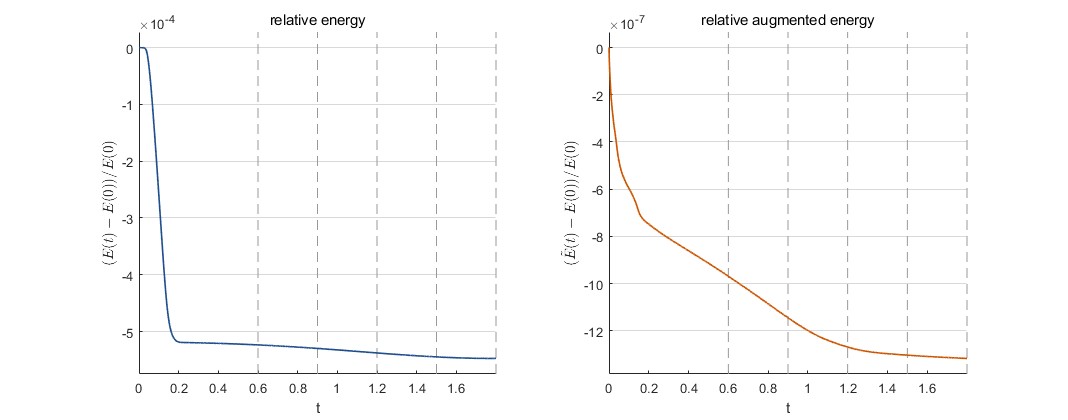}
\captionsetup{font=scriptsize}
\caption{Gaussian-shape hump with stochastic bottom: relative change in the
energy $E_\Delta$ (left) and in the augmented energy $\widetilde E_\Delta$
(right). Dashed lines mark the output times. Note the different vertical
scales, $10^{-4}$ and $10^{-7}$.}
\label{fig:num-hump-energy}
\end{figure}

The relative energy in \cref{fig:num-hump-energy} drops sharply until
$t\approx0.2$, when the left-going wave leaves the domain, and then
decreases slowly as the reflected waves reach the outflow boundaries. The
augmented energy, which accounts for the boundary energy flux, decreases
monotonically, and its change is three orders of magnitude smaller than
that of the energy; its decay is
fastest while the limiter acts on the initial discontinuity and when the
wave reaches the hump. This is consistent with \cref{thm:ec-es}, with
the augmented energy identity \eqref{eq:augmented-entropy-stability} and
with the energy plots of \cite[Figure~4]{epshteyn2025energy}.

\subsection{Modified two-dimensional perturbation of a lake at rest}
\label{sec:num-perturbation}

This example also serves as a test of the well-balanced property of the
scheme (Proposition~\ref{prop:wb}): the initial state is a small
perturbation of a lake at rest over a non-flat bottom.
Following \cite[Section~6.5]{epshteyn2025energy}, a stochastic perturbation is
placed on the lake at rest,
\begin{equation}\label{eq:num-pert-ic}
  w(x,y,0,\xi)=\begin{cases}1+0.001(\xi+1), & |xy|\le0.05,\\ 1, & \text{otherwise},\end{cases}
  \qquad u(x,y,0)=v(x,y,0)=0,
\end{equation}
over the deterministic bottom
\begin{equation}\label{eq:num-pert-bottom}
  B_a(x,y)=\begin{cases}
    a\,\big(\cos(5\pi(xy+0.35))+1\big), & -0.55<xy<-0.15,\\
    0.125\big(\cos(10\pi(xy-0.35))+1\big), & 0.25<xy<0.45,\\
    0, & \text{otherwise}.
  \end{cases}
\end{equation}
This is a two-dimensional modification of the classical
lake-at-rest perturbation example of \cite{LEVEQUE1998346}; see also
\cite{chertock2015well,dai2024energy}. Two
bottoms are considered, as in \cite[Eqs.~(6.15a) and~(6.15b)]{epshteyn2025energy}:
case~1 with $a=0.25$, where the hump in the quadrants $xy<0$ rises to $0.5$,
and case~2 with $a=0.45$, where it rises to $0.9$ and leaves a water depth of
only $0.1$ above its crest. The hump in the quadrants $xy>0$ (height $0.25$)
is the same in both cases. The bottom is symmetric under
$(x,y)\mapsto(-x,-y)$ and $x\leftrightarrow y$, and so is the initial
perturbation. The waves have amplitudes comparable to the imposed
perturbation of $10^{-3}$, so the example is sensitive to any loss of the
well-balanced property. Both cases are computed with the augmented energy
recorded and with the same discretization, listed in
\cref{tab:num-setup}.

\begin{figure}[htbp]
\centering
\includegraphics[width=0.9\textwidth]{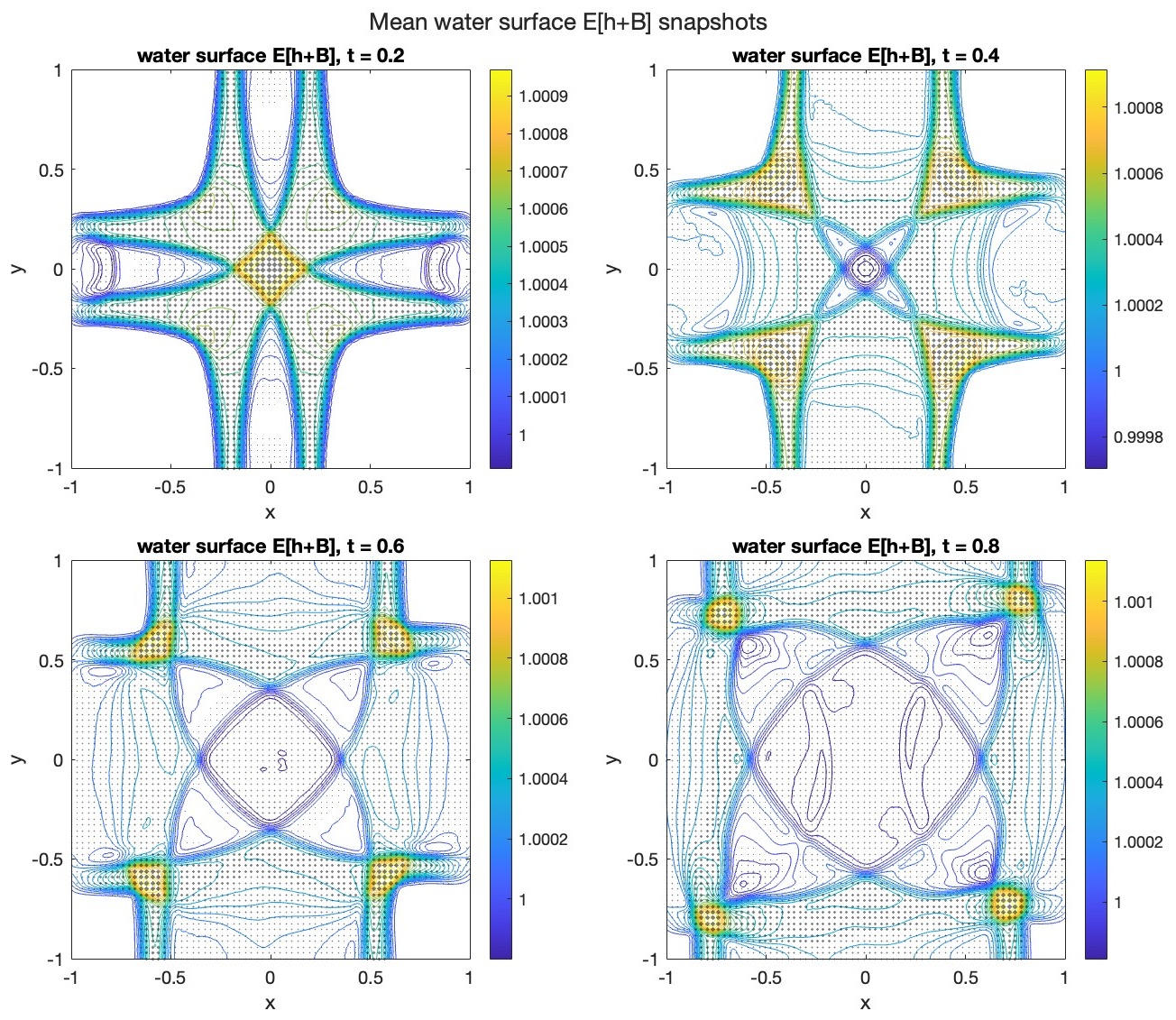}
\captionsetup{font=scriptsize}
\caption{Modified perturbation of a lake at rest, case~1 ($a=0.25$): contours
of the mean water surface $\E[w]$ at $t=0.2$, $0.4$, $0.6$, $0.8$ (left to
right, top to bottom). Disk area is proportional to $\Std[w]$, normalized by
its maximum in each panel. $\mathbb P^2$ DG, $\kpce=2$, $\Delta_0=0.025$,
$\mathrm{tol}=0.02$.}
\label{fig:num-perturbation-case1}
\end{figure}

\begin{figure}[htbp]
\centering
\includegraphics[width=0.9\textwidth]{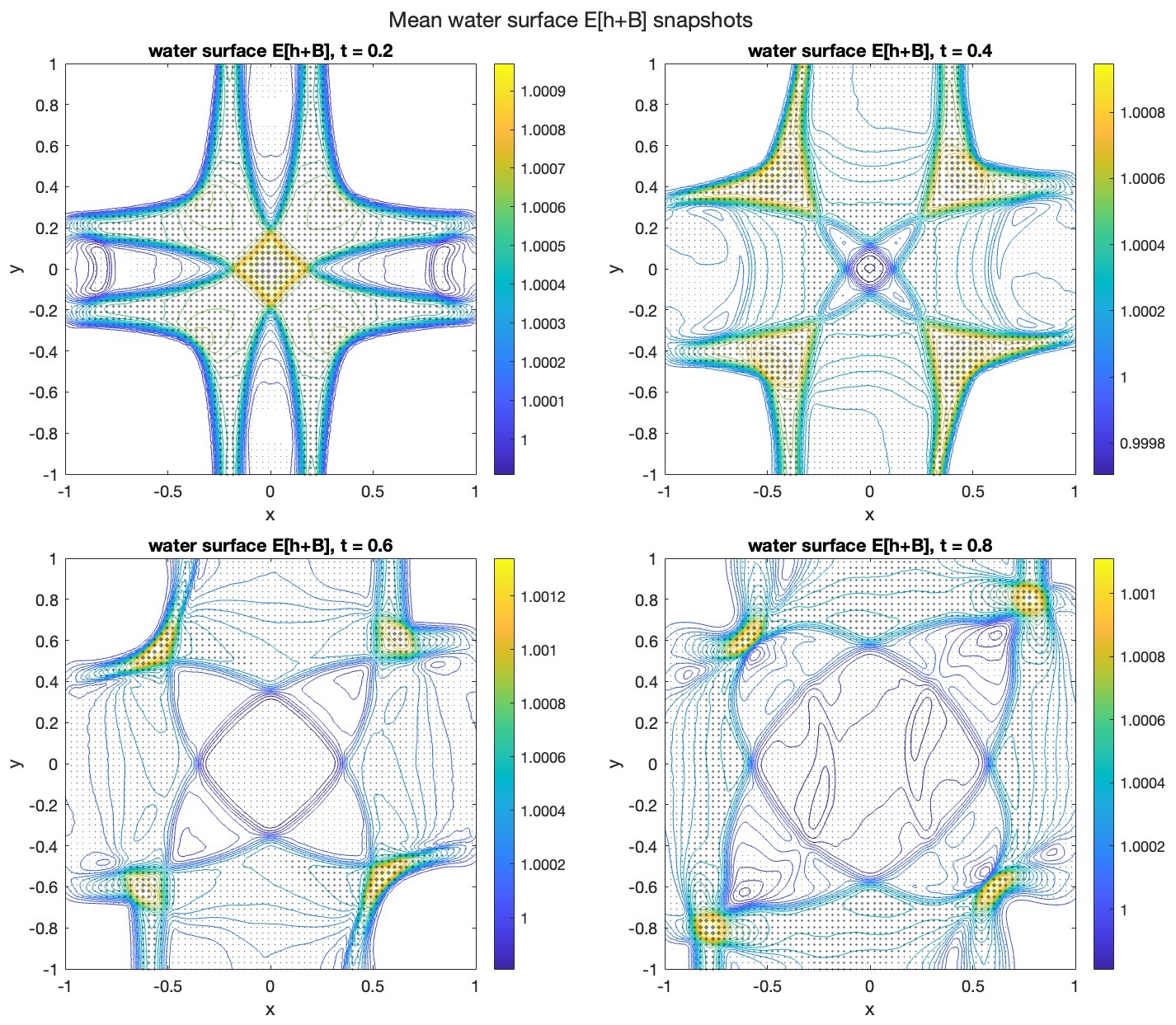}
\captionsetup{font=scriptsize}
\caption{Modified perturbation of a lake at rest, case~2 ($a=0.45$): contours
of the mean water surface $\E[w]$ at $t=0.2$, $0.4$, $0.6$, $0.8$. Same
discretization and plotting as \cref{fig:num-perturbation-case1}.}
\label{fig:num-perturbation-case2}
\end{figure}

\Cref{fig:num-perturbation-case1,fig:num-perturbation-case2} show the
mean water surface and its standard deviation for the two bottoms. The
cross-shaped perturbation propagates along the diagonals toward the
corners of the domain, while its parts touching $x=\pm1$ leave through the
outflow boundaries; the wave speed $\sqrt{gh}$ decreases over the humps, so
the fronts bend over them. In the quadrants $xy>0$, where the bottom is
identical, the two solutions coincide; in the quadrants $xy<0$ the crests
of case~2 at $t=0.6$ and $0.8$ lie closer to the center than those of
case~1 and follow the shape of the higher hump. Since the uncertainty
originates from the initial water surface, $\Std[w]$ is distributed in
proportion to the mean perturbation and travels with the waves. The mean
perturbation remains of the size of the initial perturbation $10^{-3}$
throughout; a scheme that is not well-balanced produces in this example a
different wave structure and larger perturbations
\cite[Figure~14]{epshteyn2025energy}. The wave patterns agree with the
well-balanced finite-volume results of \cite[Figures~12 and~13]{epshteyn2025energy},
with sharper fronts on the unstructured DG mesh.

\begin{figure}[htbp]
\centering
\begin{subfigure}{\textwidth}\centering
  \includegraphics[width=0.9\linewidth]{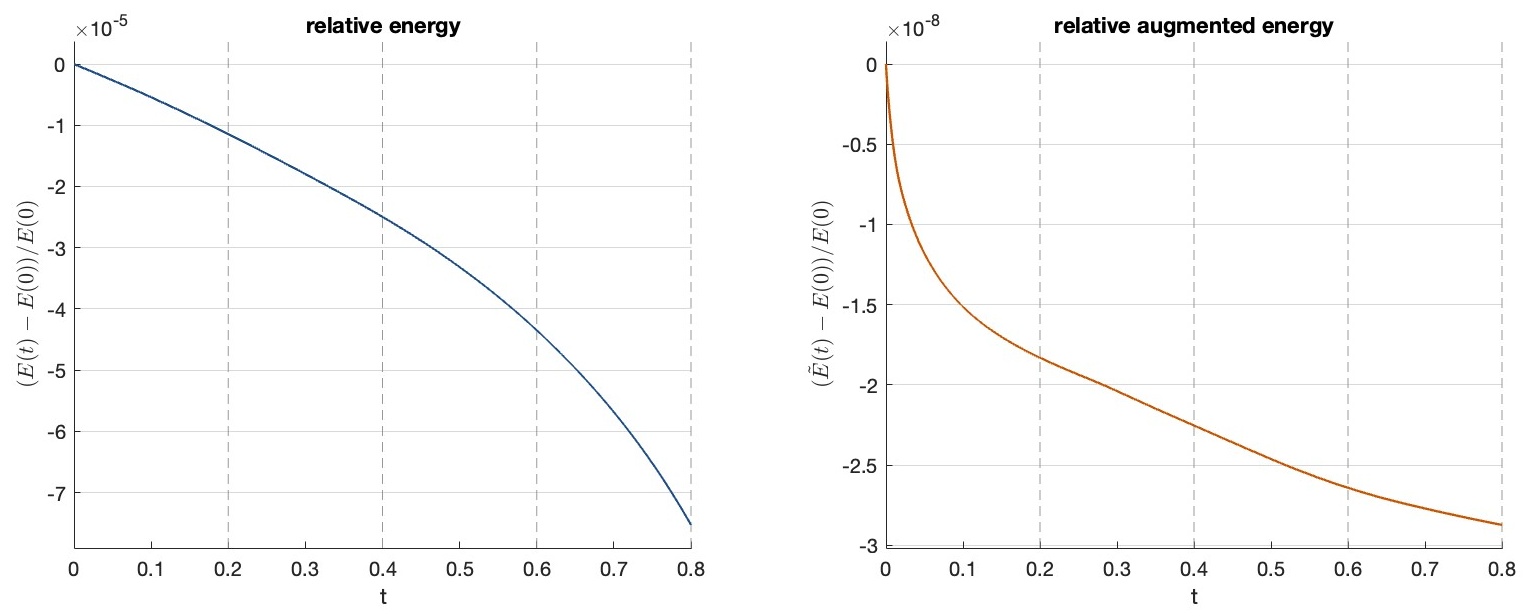}
  \caption{Case~1, $a=0.25$}
\end{subfigure}

\medskip
\begin{subfigure}{\textwidth}\centering
  \includegraphics[width=0.9\linewidth]{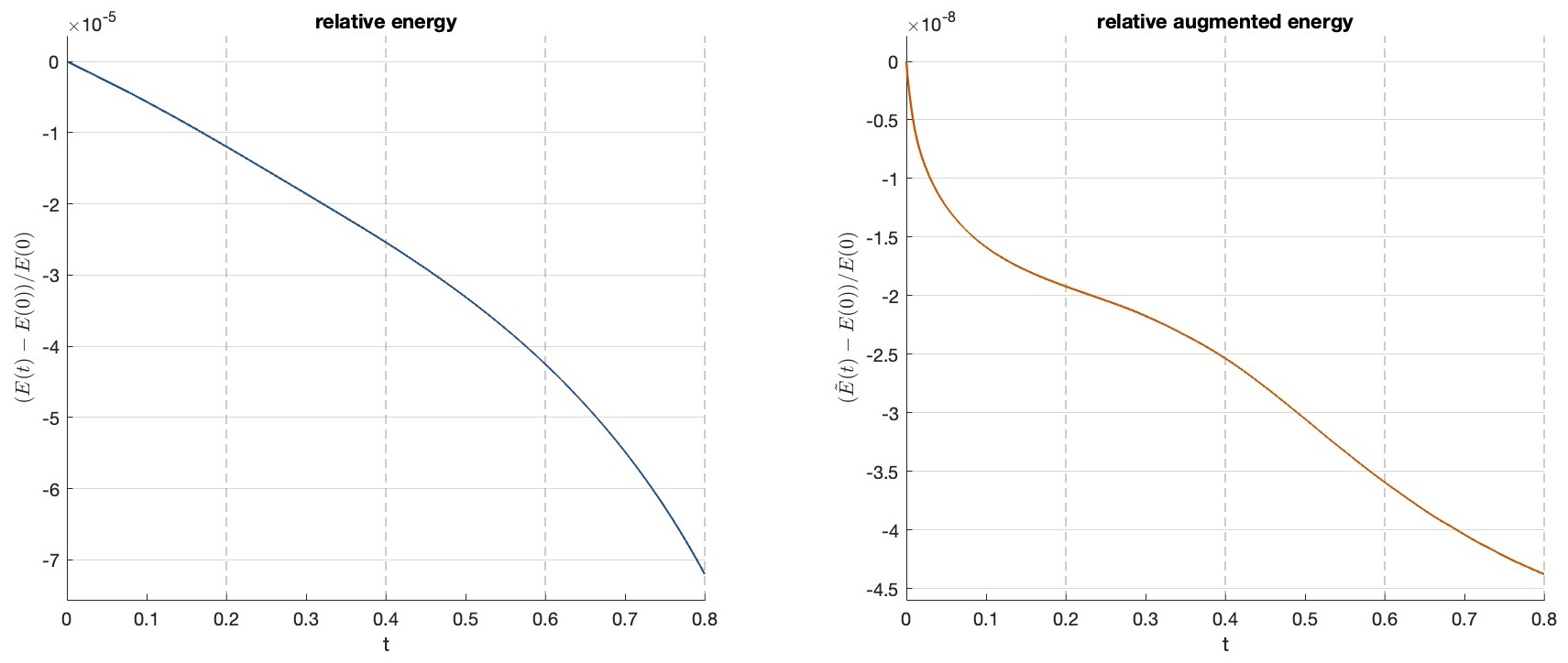}
  \caption{Case~2, $a=0.45$}
\end{subfigure}
\captionsetup{font=scriptsize}
\caption{Modified perturbation of a lake at rest: relative change in the
energy $E_\Delta$ (left) and in the augmented energy $\widetilde E_\Delta$
(right). Dashed lines mark the output times. Note the different vertical
scales, $10^{-5}$ and $10^{-8}$.}
\label{fig:num-perturbation-energy}
\end{figure}

In \cref{fig:num-perturbation-energy} the relative energy decreases
from the start, because the initial perturbation touches the outflow
boundaries, and almost identically for the two bottoms. The augmented
energy decreases monotonically in both cases, and its change is three
orders of magnitude smaller than that of the energy; its decay is fastest
while the limiter acts on the edges
of the initial perturbation and, in case~2, after $t\approx0.4$, when the
waves steepen over the higher hump, consistent with \cref{thm:ec-es}.

\subsection{A submerged flat plateau}
\label{sec:num-plateau}

This example also serves as a test of the well-balanced property
(Proposition~\ref{prop:wb}), here for a lake at rest over a stochastic
bottom.
Following \cite[Section~6.4]{epshteyn2025energy}, a deterministic perturbation
of size $10^{-4}$ is placed on the lake at rest,
\begin{equation}\label{eq:num-plateau-ic}
  w(x,y,0,\xi)=\begin{cases}1.0001, & -0.4<x<-0.3,\\ 1, & \text{otherwise},\end{cases}
  \qquad u(x,y,0)=v(x,y,0)=0,
\end{equation}
over the stochastic bottom
\begin{equation}\label{eq:num-plateau-bottom}
  B(x,y,\xi)=\begin{cases}
    0.9998, & r\le0.1,\\
    9.997(0.2-r)+0.0001(\xi+1), & 0.1<r\le0.2,\\
    0.0001, & \text{otherwise},
  \end{cases}
  \qquad r:=\sqrt{x^2+y^2}+0.0001.
\end{equation}
The plateau lies only $2\times10^{-4}$ below the
water surface, and the randomness enters only on its slope $0.1<r\le0.2$. The
water layer above the plateau is very thin, so this example exercises the
velocity desingularization, the hyperbolicity-preserving scaling and the
stage restart together with the well-balanced property. The case-specific
parameters are listed in \cref{tab:num-setup}; both the energy and the
augmented energy are recorded.

\begin{figure}[htbp]
\centering
\includegraphics[width=0.8\textwidth,trim={0 0 0 0},clip]{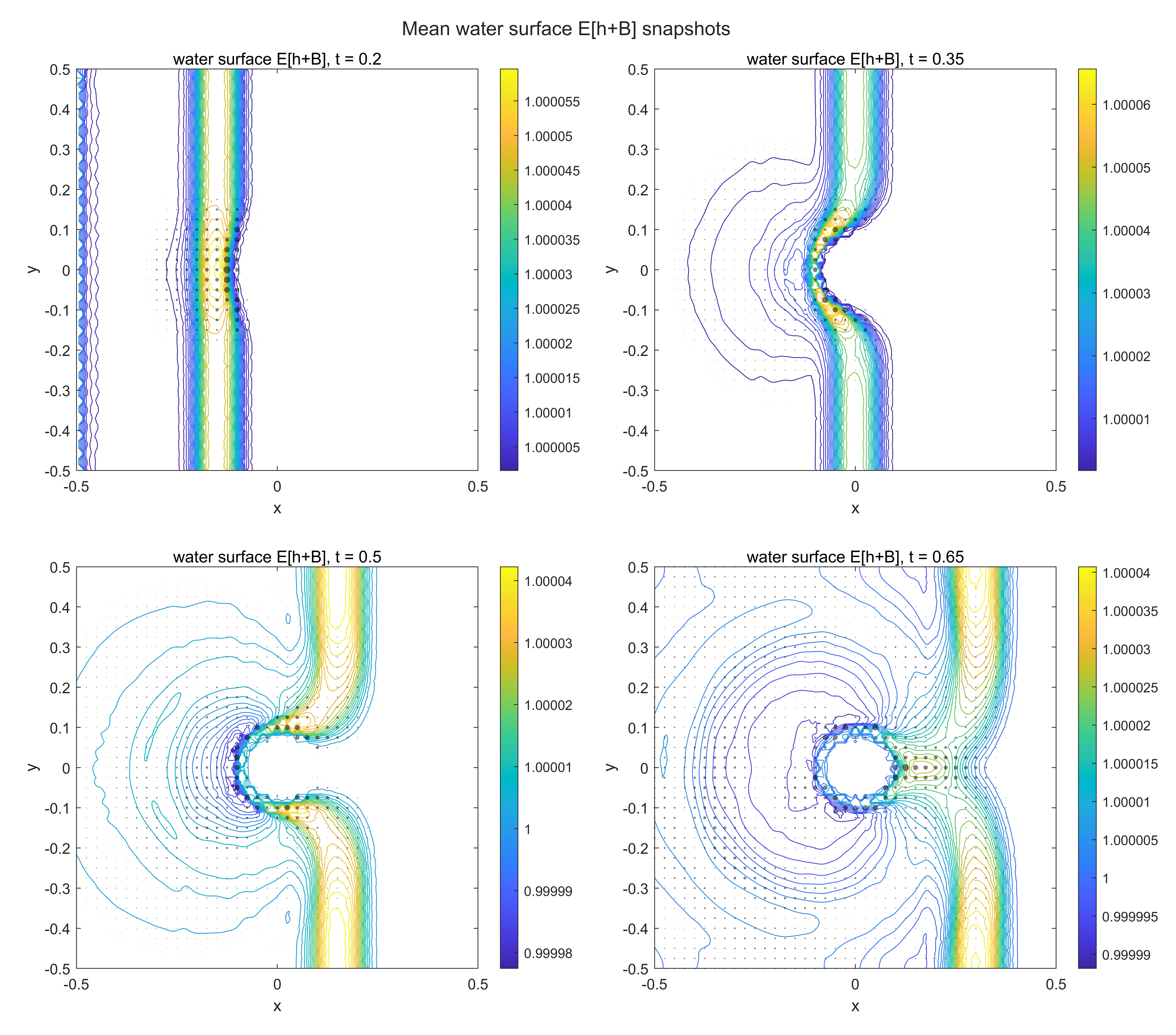}
\captionsetup{font=scriptsize}
\caption{Submerged flat plateau: contours of the mean water surface $\E[w]$
at $t=0.2$, $0.35$, $0.5$, $0.65$ (left to right, top to bottom). Disk area
is proportional to $\Std[w]$, normalized by its maximum in each panel; the
maximum standard deviation is $6.8\times10^{-9}$, $2.6\times10^{-8}$,
$2.2\times10^{-8}$ and $1.4\times10^{-8}$, respectively. $\mathbb P^2$ DG,
$\kpce=2$, $\Delta_0=0.025$, $\mathrm{tol}=0.02$.}
\label{fig:num-plateau-snapshots}
\end{figure}

The initial perturbation generates two waves. The left-going wave
leaves the domain through $x=-0.5$; the right-going wave slows down
strongly over the plateau, where the water depth is only $2\times10^{-4}$,
wraps around it and sends reflected waves back, producing the patterns at
$t=0.5$ and $0.65$ in \cref{fig:num-plateau-snapshots}, in agreement with
the finite-volume results of \cite[Figure~7]{epshteyn2025energy}. The
velocity desingularization and the hyperbolicity-preserving scaling are
active in the nearly dry region above the plateau. Away from the waves, in
particular on the random slope ahead of the right-going wave and on the
plateau itself, the deviation of the mean water surface from the rest state
and the standard deviation $\Std[w]$ remain at
round-off level although $\Std[B]$ there is of order $10^{-4}$; throughout
the evolution the displacement of the mean surface does not exceed the
initial perturbation $10^{-4}$, and $\Std[w]$ stays below $3\times10^{-8}$
and concentrated on the random slope. The scheme thus preserves the
stochastic lake at rest under a small perturbation in the presence of a
nearly dry region. The small structures along $x=-0.5$ at $t=0.2$ are
produced by the limiter acting on the boundary cells while the left-going
wave leaves the domain.

\begin{figure}[htbp]
\centering
\includegraphics[width=0.8\textwidth]{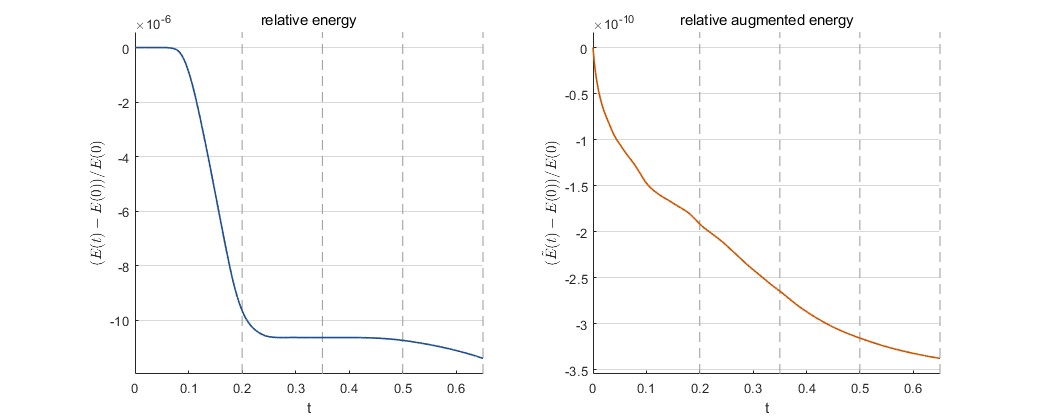}
\captionsetup{font=scriptsize}
\caption{Submerged flat plateau: relative change in the energy $E_\Delta$
(left) and in the augmented energy $\widetilde E_\Delta$ (right). Dashed
lines mark the output times. Note the different vertical scales, $10^{-6}$
and $10^{-10}$.}
\label{fig:num-plateau-energy}
\end{figure}

The relative energy in \cref{fig:num-plateau-energy} decreases while
the left-going wave leaves the domain, is nearly constant during the
interaction with the plateau and decreases again when the reflected waves
reach the left boundary; its small isolated increases are due to the
boundary exchange. The augmented energy decreases monotonically throughout,
and its change is more than four orders of magnitude smaller than that of
the energy, with the fastest decay
while the limiter acts on the initial perturbation, as expected from the
entropy stability of the scheme.

\subsection{Circular dam break}
\label{sec:num-circular}

Following \cite[Section~3.2.1]{seaid2004non} (see also
\cite{toro2001shock}), the domain is $\Omega=[0,50]^2$
with a flat bottom $B=0$, and a cylindrical dam of radius $11$ centered at
$(25,25)$ holds the water at rest. The dam is removed at $t=0$. We make the
inner water height uncertain,
\begin{equation}\label{eq:num-circ-ic}
  h(x,y,0,\xi)=\begin{cases}10+0.1\,\xi, & r\le11,\\ h_{\mathrm{out}}, & r>11,\end{cases}
  \qquad u(x,y,0)=v(x,y,0)=0,
\end{equation}
where $r$ is the distance from $(25,25)$, with $h_{\mathrm{out}}=1$ (wet bed)
or $h_{\mathrm{out}}=10^{-12}$ (dry bed). The parameters are
listed in \cref{tab:num-setup}.
The initial data are obtained by the $L^2$ projection on the cut cells,
followed by the limiter and the hyperbolicity-preserving scaling before the
first step; the floor of the scaling is $10^{-6}$. On the dry bed the cells ahead of the
front hold a water height of $10^{-12}$, so the run relies on the
desingularization, the hyperbolicity-preserving scaling, a dry-cell limiter
\cite{fu2022high} that flags cells whose minimal sampled height is below
$\varepsilon_{\mathrm{dry}}$ and reduces them to their cell averages, and a
cap on the velocity magnitude, which is needed when
$\varepsilon_{\mathrm{dry}}$ is small.

\begin{figure}[htbp]
\centering
\begin{subfigure}[t]{0.45\textwidth}\centering
  \includegraphics[width=\linewidth,trim={74 40 74 52},clip]{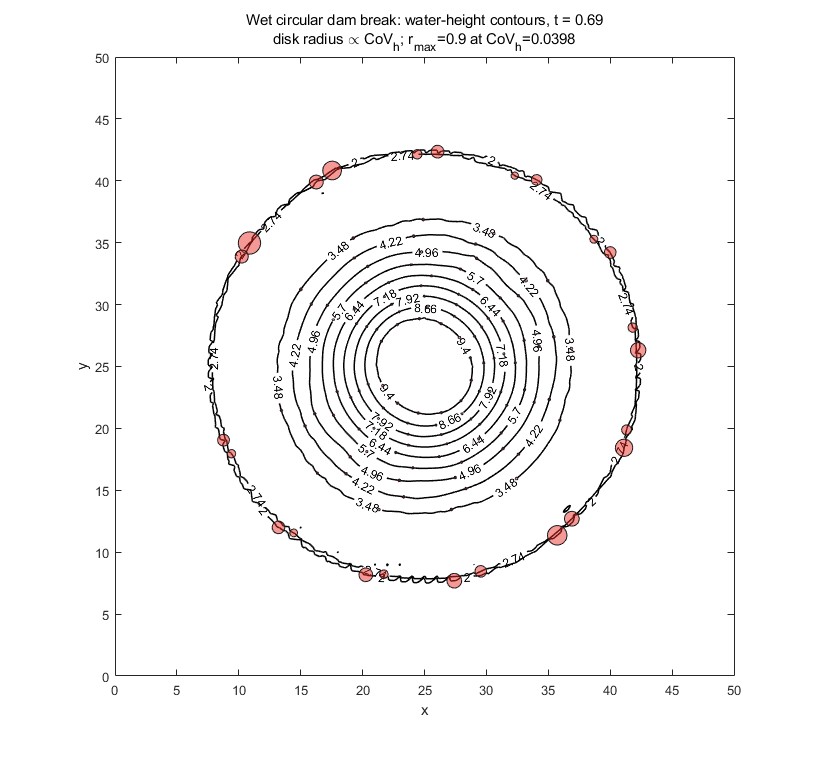}
  \caption{Wet bed, $h_{\mathrm{out}}=1$}
\end{subfigure}\hfill
\begin{subfigure}[t]{0.45\textwidth}\centering
  \includegraphics[width=\linewidth,trim={74 40 74 52},clip]{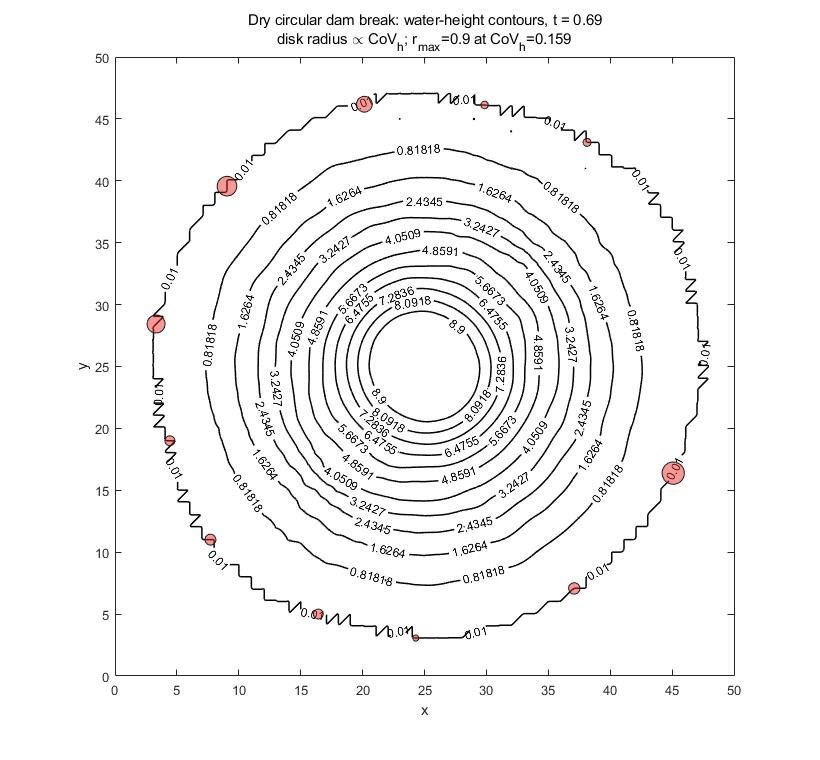}
  \caption{Dry bed, $h_{\mathrm{out}}=10^{-12}$}
\end{subfigure}
\captionsetup{font=scriptsize}              
\caption{Circular dam break at $t=0.69$: contours of $\E[h]$ at the levels
$h=2,\,2.74,\,3.48,\dots,9.40$ on the wet bed and
$h=0.01,\,0.82,\,1.63,\dots,8.90$ on the dry bed, with disks whose radius
is proportional to $\mathrm{CoV}_h=\Std[h]/|\E[h]|$ at twelve points sampled
along each level; the largest disk corresponds to $\mathrm{CoV}_h=0.040$ on
the wet bed and $0.159$ on the dry bed. Representative values of
$\mathrm{CoV}_h$ on the successive levels (one sampling point per level,
outermost level first) are $0.040$, $0.023$, $0.0049$, $0.0046$, $0.0044$,
$0.0041$, $0.0039$, $0.0038$, $0.0037$, $0.0037$, $0.0040$ on the wet
bed and $0.012$, $0.010$, $0.0081$, $0.0063$, $0.0055$, $0.0049$, $0.0045$,
$0.0042$, $0.0039$, $0.0038$, $0.0037$, $0.0037$ on the dry bed.
$\mathbb P^2$ DG, $\kpce=2$, $5\,000$ triangles, $\mathrm{tol}=0.02$.}
\label{fig:num-circular}
\end{figure}

\Cref{fig:num-circular} shows the contours of $\E[h]$ at $t=0.69$, at
the levels of the deterministic solution of
\cite[Figures~9 and~10]{seaid2004non}, with disks whose radius is
proportional to the coefficient of variation
$\mathrm{CoV}_h=\Std[h]/|\E[h]|$. On the wet bed the removal of the dam
produces an outward circular bore and an inward rarefaction wave; on the
dry bed no bore forms and the rarefaction wave extends to the wet--dry
front. In both cases the circular symmetry is well preserved on the
split-square mesh, and the scheme remains robust with the water height
$10^{-12}$ ahead of the front, where the hyperbolicity-preserving scaling,
the desingularization and the dry-cell limiter are active;
on an unstructured triangular mesh (not shown) the solution looks
essentially the same, which indicates the robustness of the proposed scheme
with respect to the mesh.
The relative uncertainty is small near the center, where $\E[h]$ is close
to its random initial value, and grows toward the front, most strongly on
the dry bed, where $\E[h]$ tends to zero while $\Std[h]$ does not (see the
largest values of $\mathrm{CoV}_h$ in the caption).

The positions of the contour levels of $\E[h]$, including the bore on
the wet bed, are visually consistent with those of the deterministic
solution of \cite[Figures~9 and~10]{seaid2004non}, which was computed with
a relaxation scheme on a $50\times50$ grid; compared with the reference
plots, the DG--SG bore appears steeper and the wet--dry front on the dry
bed slightly more advanced, as expected from the lower numerical
dissipation of the high-order scheme on the same number of cells. The
mesh-scale kinks of the $h=0.01$ contour reflect the elementwise treatment
of the dry cells.

\subsection{Curved partial dam break}
\label{sec:num-curved}

\begin{figure}[htbp]
\centering
\includegraphics[width=0.8\textwidth]{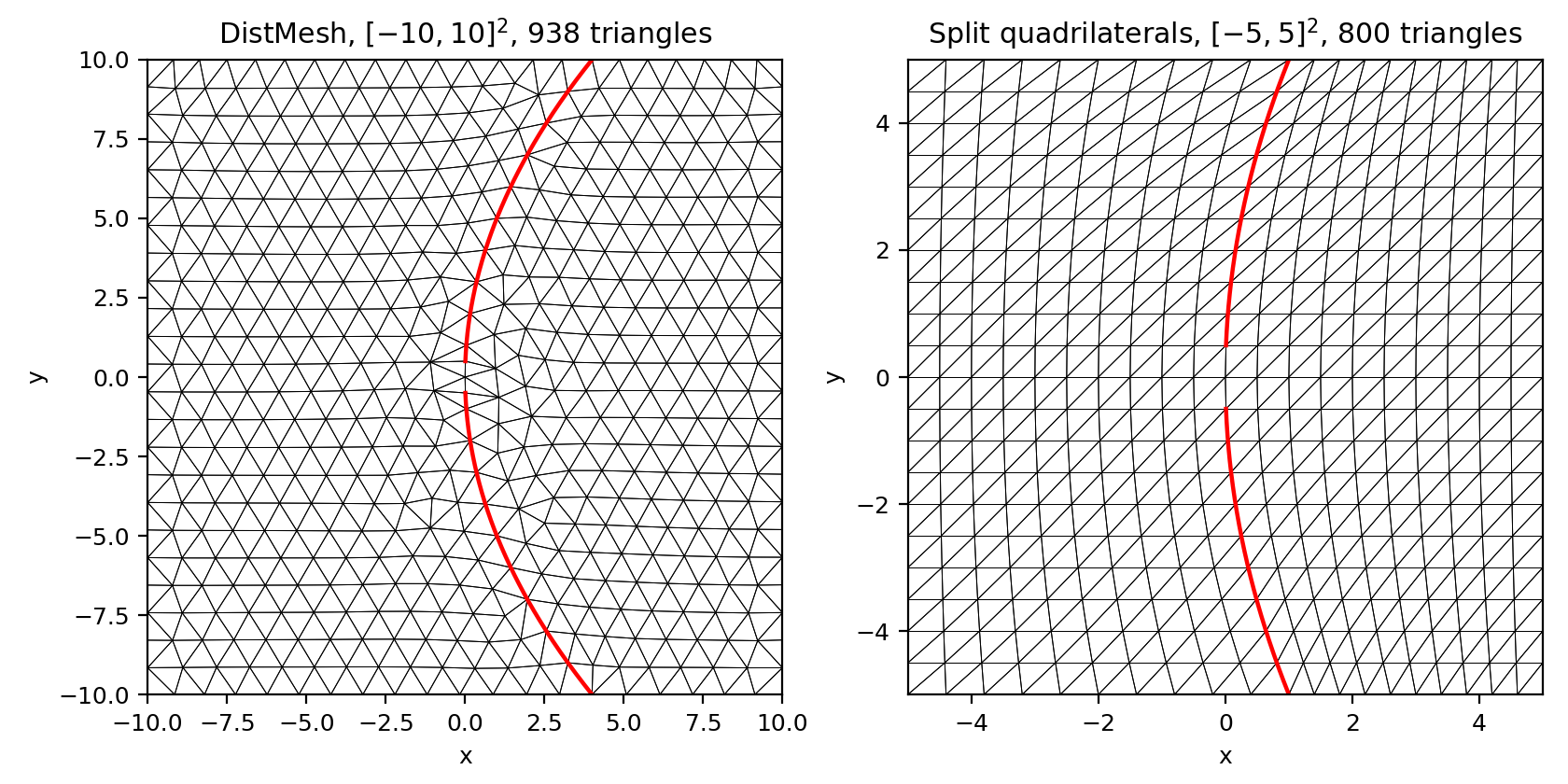}
\captionsetup{font=scriptsize}
\caption{Curved partial dam break: the two DG--SG meshes, DistMesh on
$[-10,10]^2$ (left) and split quadrilaterals on $[-5,5]^2$ (right), with the
dam \eqref{eq:num-curved-dam} for $|y|\ge0.5$ in red.}
\label{fig:num-curved-meshes}
\end{figure}

Following \cite[Section~4.3]{wu2021high}, the bottom is flat ($B=0$), and a
dam along the parabola
\begin{equation}\label{eq:num-curved-dam}
  x=\frac{y^2}{25}
\end{equation}
with an opening $-0.5<y<0.5$ separates still water of height $10$ on the left
from height $5$ on the right. The dam is a reflecting wall on both sides. We
make the left height uncertain,
\begin{equation}\label{eq:num-curved-ic}
  h(x,y,0,\xi)=\begin{cases}10+0.1\,\xi, & x\le y^2/25,\\ 5, & x>y^2/25,\end{cases}
  \qquad u(x,y,0)=v(x,y,0)=0.
\end{equation}
As in \cite{wu2021high}, the polynomial degree is $k=3$ and no
limiter is used. To test the robustness
of the scheme with respect to the geometry and to the mesh, the problem is
solved in two configurations (\cref{tab:num-setup}): on the domain
$[-10,10]^2$ with a DistMesh triangulation whose element edges are fitted to
\eqref{eq:num-curved-dam}, and on the domain $[-5,5]^2$ with a structured
mesh, a $20\times20$ grid of quadrilaterals mapped to fit the same parabola
and opening and split into triangles. In both meshes, the dam is represented
by straight element edges (\cref{fig:num-curved-meshes}). This example tests
the scheme on a domain with an internal wall and, since no limiter is used,
the robustness of the ES flux and of the hyperbolicity-preserving scaling at a strong
discontinuity.

\begin{figure}[htbp]
\centering
\begin{subfigure}{\textwidth}\centering
  \includegraphics[width=\linewidth,trim={200 40 55 55},clip]{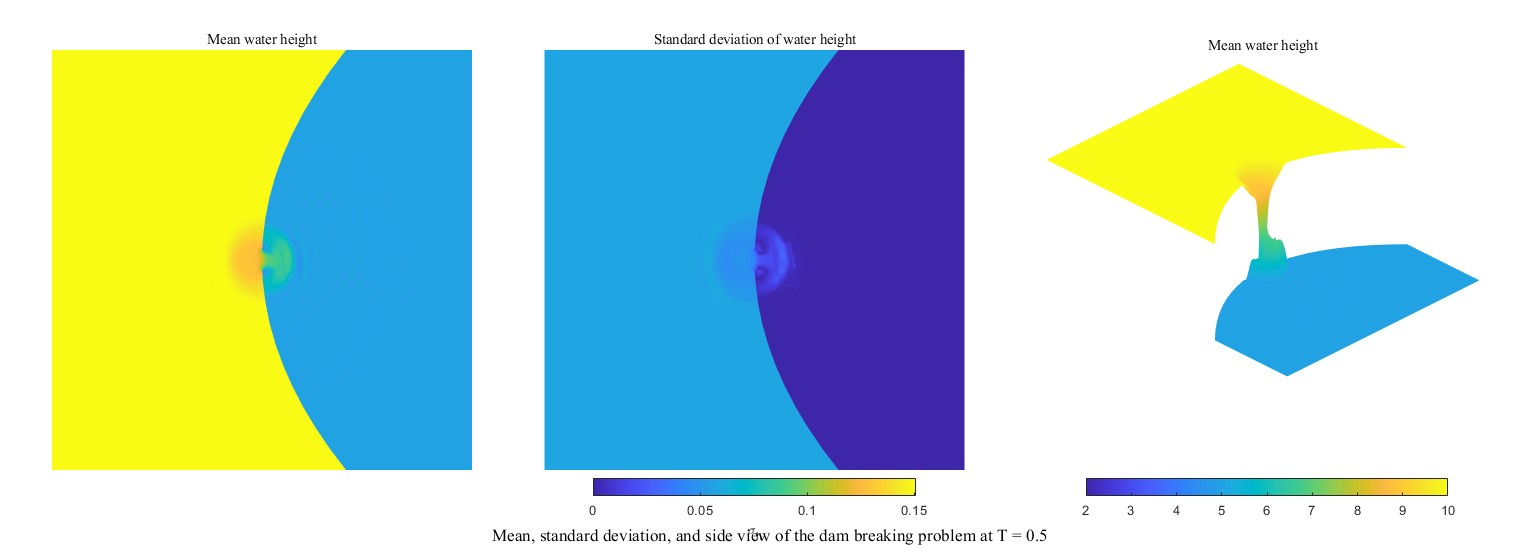}
  \captionsetup{font=scriptsize}
  \caption{$t=0.5$, $\max\Std[h]=0.061$}
\end{subfigure}

\medskip
\begin{subfigure}{\textwidth}\centering
  \includegraphics[width=\linewidth,trim={200 40 55 55},clip]{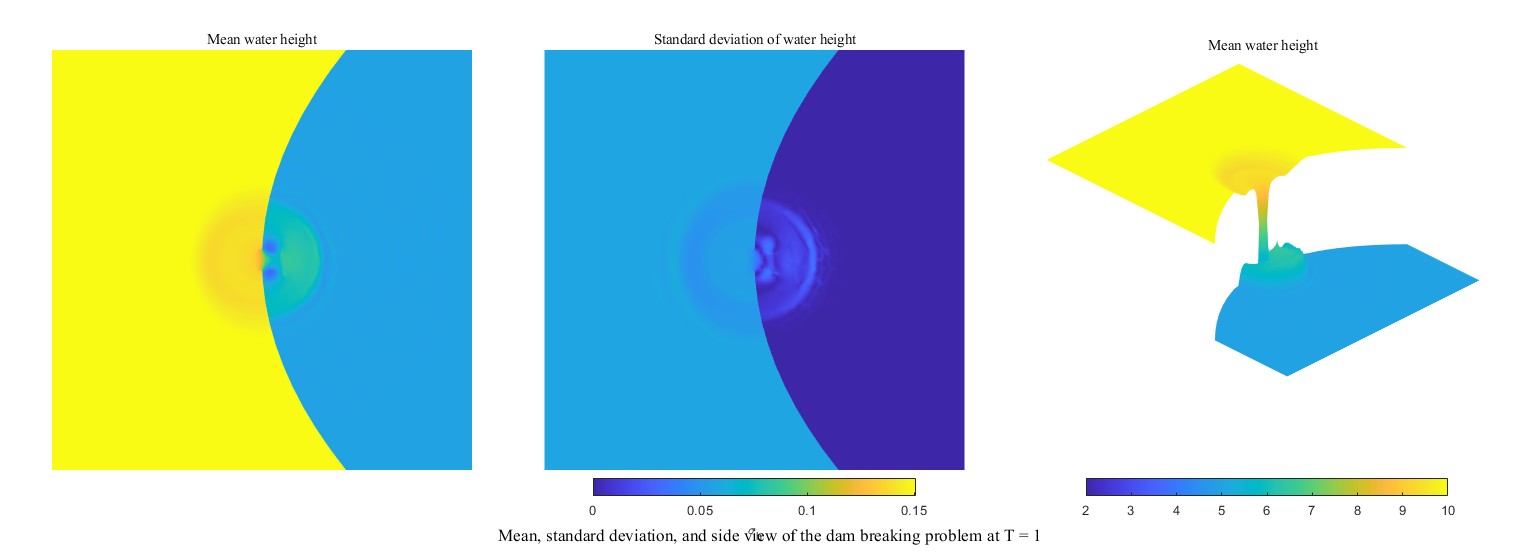}
  \captionsetup{font=scriptsize}
  \caption{$t=1$, $\max\Std[h]=0.062$}
\end{subfigure}

\medskip
\begin{subfigure}{\textwidth}\centering
  \includegraphics[width=\linewidth,trim={200 40 55 55},clip]{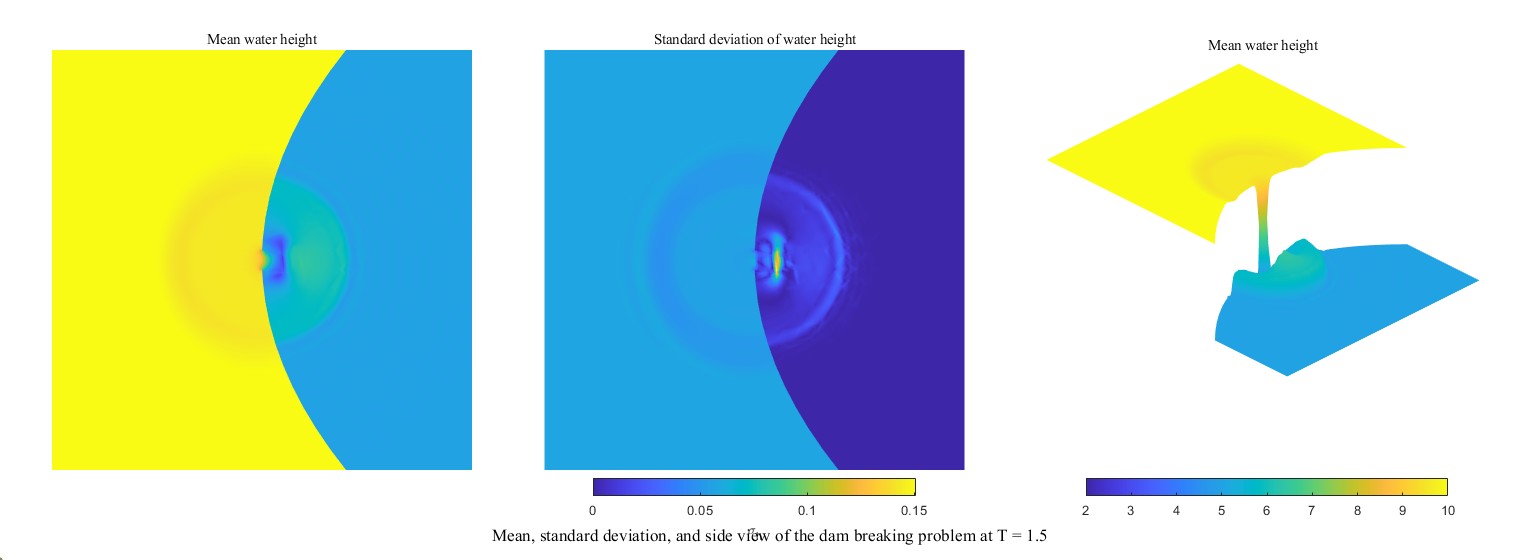}
  \captionsetup{font=scriptsize}
  \caption{$t=1.5$, $\max\Std[h]=0.151$}
\end{subfigure}
\captionsetup{font=scriptsize}
\caption{Curved partial dam break, DistMesh mesh on $[-10,10]^2$. Left:
$\E[h]$; middle: $\Std[h]$; right: oblique view of $\E[h]$. $\mathbb P^3$ DG,
$\kpce=2$, no limiter.}
\label{fig:num-curved-distmesh}
\end{figure}

\begin{figure}[htbp]
\centering
\begin{subfigure}{\textwidth}\centering
  \includegraphics[width=\linewidth,trim={200 40 55 55},clip]{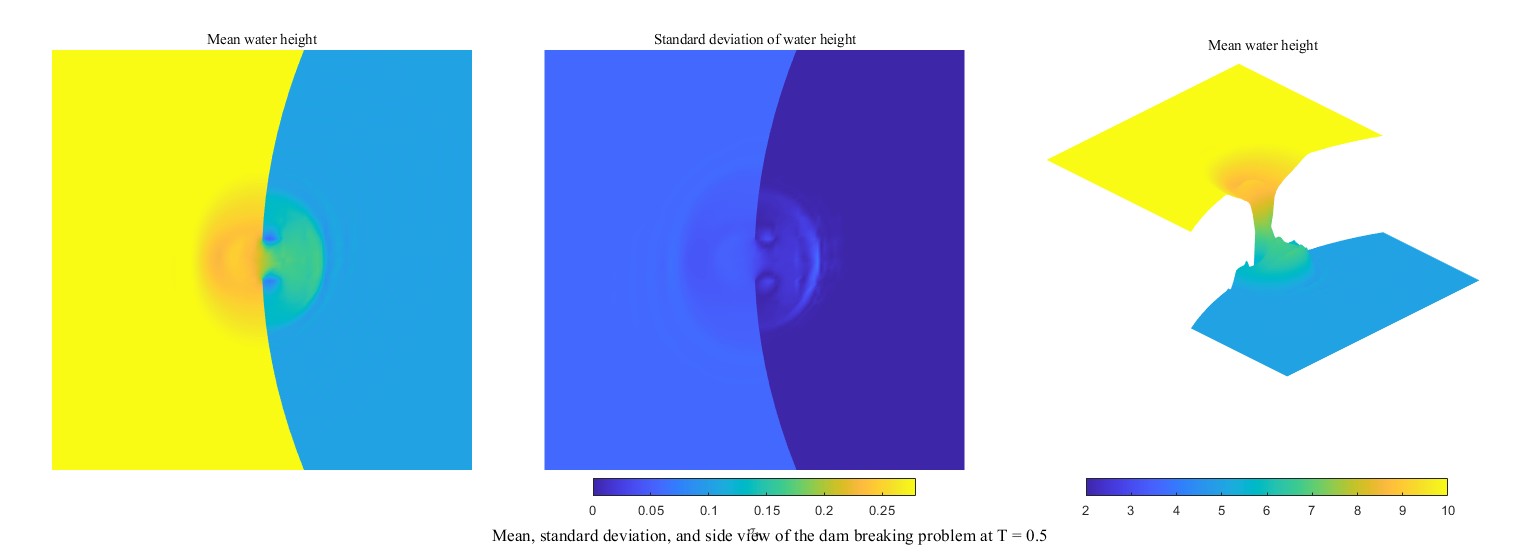}
  \captionsetup{font=scriptsize}
  \caption{$t=0.5$, $\max\Std[h]=0.065$}
\end{subfigure}

\medskip
\begin{subfigure}{\textwidth}\centering
  \includegraphics[width=\linewidth,trim={200 40 55 55},clip]{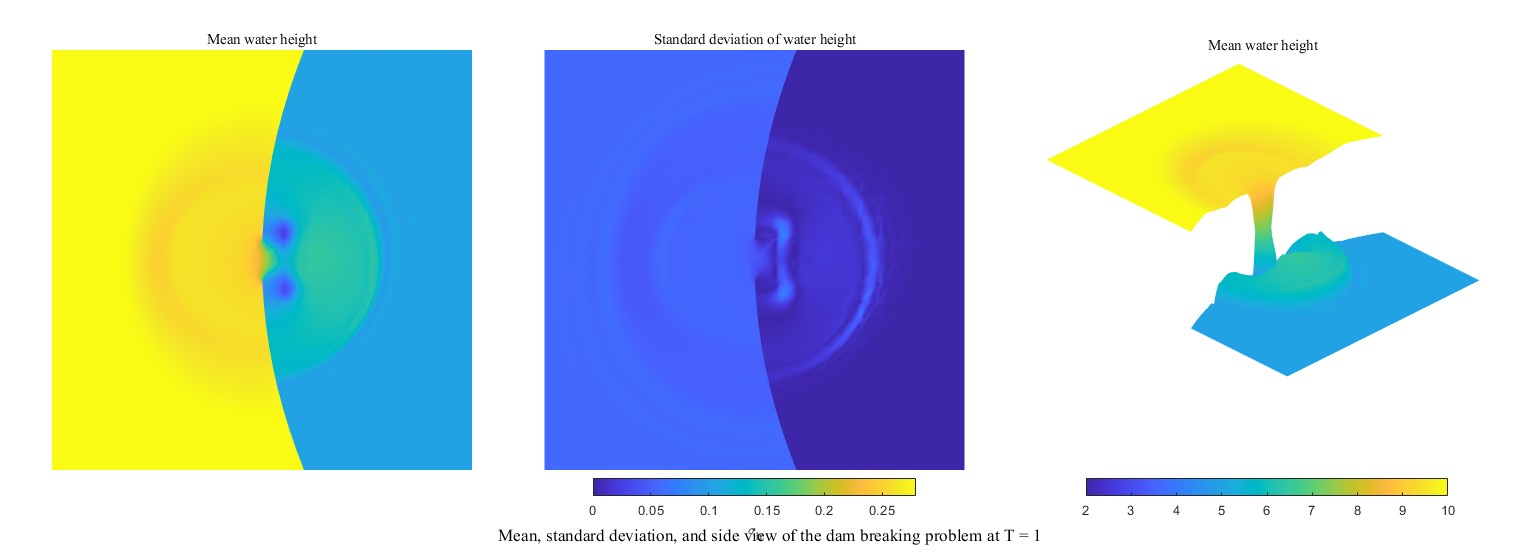}
  \captionsetup{font=scriptsize}
  \caption{$t=1$, $\max\Std[h]=0.083$}
\end{subfigure}

\medskip
\begin{subfigure}{\textwidth}\centering
  \includegraphics[width=\linewidth,trim={200 40 55 55},clip]{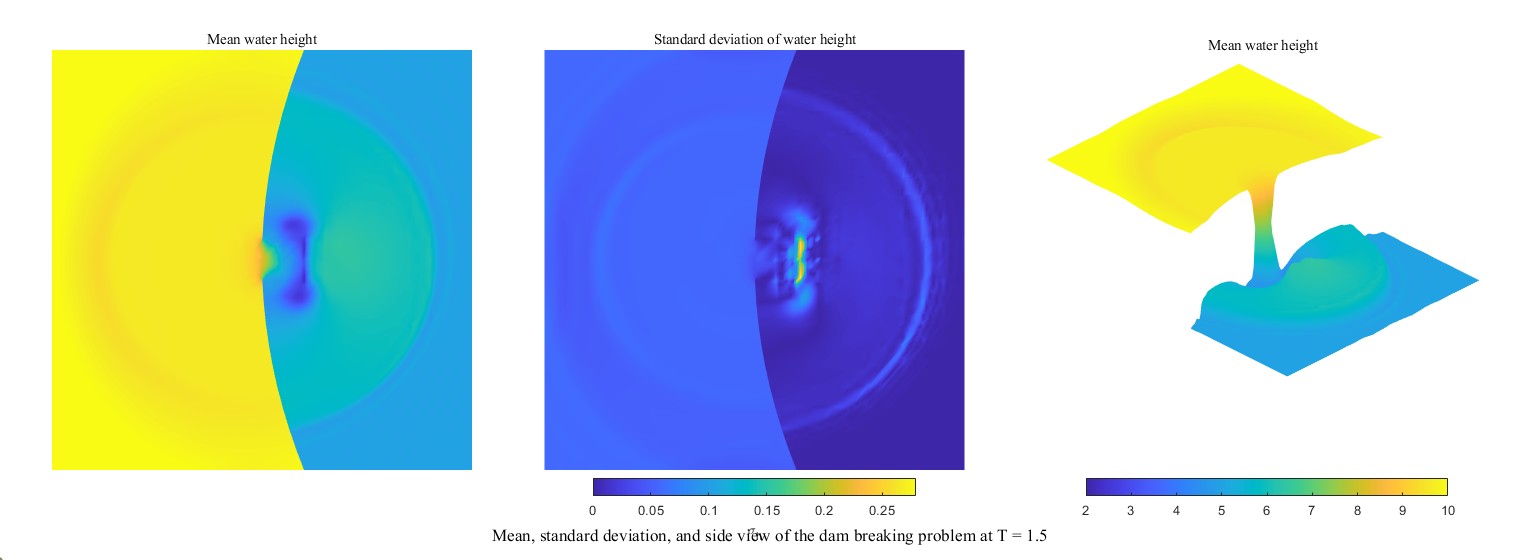}
  \captionsetup{font=scriptsize}
  \caption{$t=1.5$, $\max\Std[h]=0.279$}
\end{subfigure}
\captionsetup{font=scriptsize}
\caption{Curved partial dam break, split-quadrilateral mesh on $[-5,5]^2$.
Left: $\E[h]$; middle: $\Std[h]$; right: oblique view of $\E[h]$.
$\mathbb P^3$ DG, $\kpce=2$, no limiter.}
\label{fig:num-curved-quad}
\end{figure}

\Cref{fig:num-curved-distmesh,fig:num-curved-quad} show, at each
output time, the mean water height, its standard deviation and an oblique
view of the mean, with the color range $[2,10]$ of \cite{wu2021high}. The
water flows through the opening as a jet and spreads into the right basin
as a roughly circular bore, a rarefaction wave propagates into the left
basin, and two shallow spots form at the ends of the opening. These are the
features of the deterministic solution of \cite[Figures~9--11]{wu2021high},
computed with curved triangles and the same degree $k=3$, and the positions
of the bore front and of the rarefaction edge at $t=0.5$, $1$, $1.5$ are
visually consistent with the reference on both meshes. The front
positions agree between the two DG--SG runs, and the straight element edges
representing the curved dam do not visibly affect the solution at this
resolution. The uncertainty, initially located in the left basin, is
carried through the opening with the flow:
$\Std[h]$ keeps its initial
value $0.1/\sqrt3\approx0.058$ in the undisturbed part of the left basin and
vanishes in the undisturbed part of the right basin, and at $t=1.5$ its
maximum is attained in a thin layer across the jet downstream of the
opening. This maximum is larger on the finer split-quadrilateral mesh. Both
computations remain stable up to $T=1.5$ without any limiting, which
illustrates the robustness of the ES flux combined with the
hyperbolicity-preserving scaling at a strong discontinuity.

\subsection{Gaussian-shape hump with two independent random variables}
\label{sec:num-2rv}

Following \cite[Section~6.6, Eqs.~(6.16)--(6.17)]{epshteyn2025energy}, the
uncertainties are in the two widths of the hump,
\begin{equation}\label{eq:num-hump2rv-bottom}
  B(x,y,\bm\xi)=0.5\exp\!\big(-12.5(\xi_1+1)(x-1)^2-25(\xi_2+1)(y-0.5)^2\big),
\end{equation}
with the deterministic initial uniform flow of the accuracy test,
\begin{equation}\label{eq:num-hump2rv-ic}
  w(x,y,0,\bm\xi)=1,\qquad u(x,y,0,\bm\xi)=0.3,\qquad v(x,y,0,\bm\xi)=0.
\end{equation}
The random vector
$\bm\xi=(\xi_1,\xi_2)$ has independent components on $[-1,1]$:
$\xi_1\sim\operatorname{Beta}_{[-1,1]}(1,3)$, that is, $(\alpha,\beta)=(2,0)$ in
the parameterization of \cref{sec:numerics}, and $\xi_2$ is uniform.
The PCE basis is the tensor product of the orthonormal Jacobi polynomials of
degree $\le2$ in $\xi_1$ and the Legendre polynomials of degree $\le2$ in
$\xi_2$, three terms each and $\kpce=9$ in total; the triple products are
computed with a $60\times60$ Gauss--Jacobi rule. The parameters are listed in
\cref{tab:num-setup}. This example tests the product-matrix construction
for non-identical input distributions and the coupled velocity projection
with a larger stochastic system.

\begin{figure}[htbp]
\centering
\includegraphics[width=\textwidth]{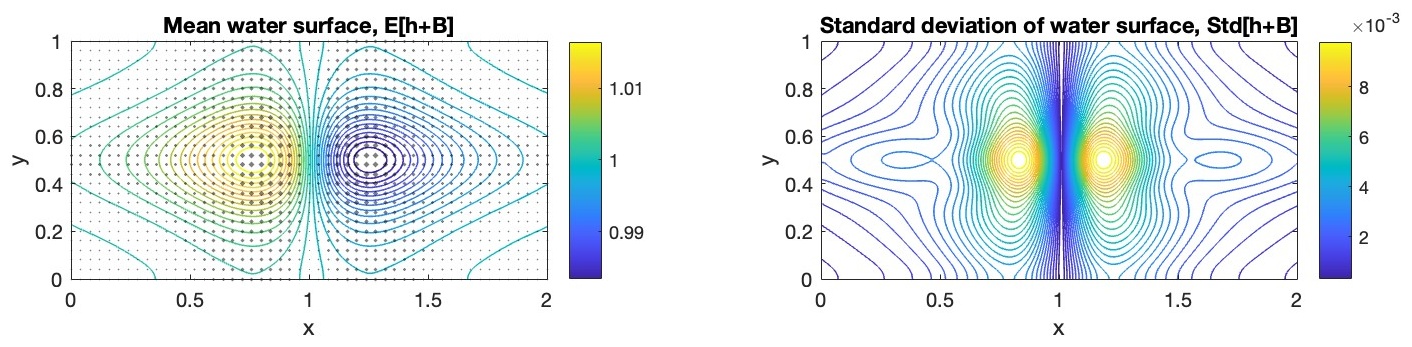}
\captionsetup{font=scriptsize}
\caption{Hump with two independent random variables at $T=0.07$. Left:
contours of the mean water surface $\E[w]$ with black disks at every fourth
plotting-grid point, area proportional to $\Std[w]$. Right: contours of
$\Std[w]$. $\mathbb P^2$ DG, $\kpce=9$, $\Delta_0=0.06$, no limiter.}
\label{fig:num-hump2rv}
\end{figure}

\Cref{fig:num-hump2rv} shows the mean and the standard deviation of
the water surface at $T=0.07$. The mean water surface shows the
elevation--depression pattern of the flow over the hump, raised upstream of
the crest and lowered downstream of it, as in
\cite[Figure~2]{epshteyn2025energy}. The uncertainty, now
controlled by the two random widths, is largest near the maximum and the
minimum of the mean surface and smallest along the line $x=1$, where the
bottom slope vanishes for every realization; this distribution agrees with
the finite-volume results of \cite[Figure~15]{epshteyn2025energy}.
The maximum of $\Std[w]$ is $1.01\times10^{-2}$, attained at
$(x,y)\approx(1.19,0.5)$. The
example demonstrates the applicability of the scheme to multivariate,
non-identically distributed inputs, with the product matrix built from the
tensor-product Jacobi--Legendre basis.

\begin{figure}[htbp]
\centering
\includegraphics[width=0.92\textwidth]{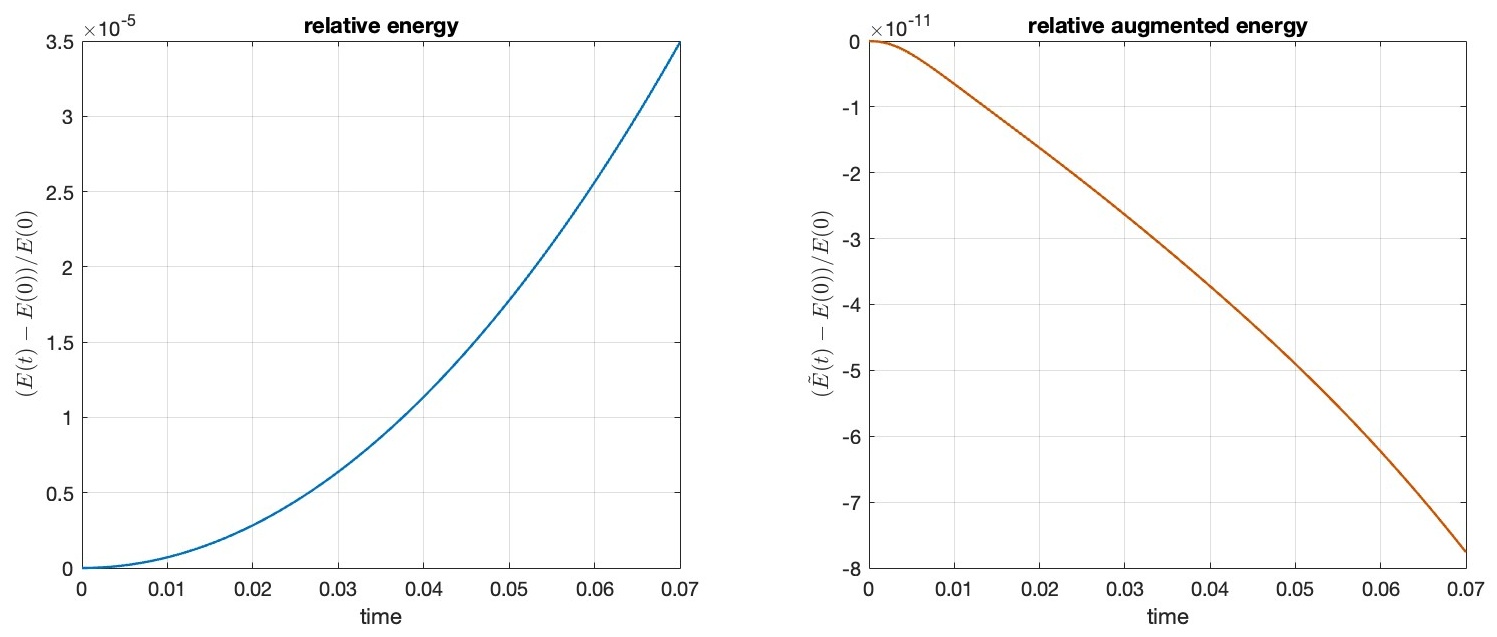}
\captionsetup{font=scriptsize}
\caption{Hump with two independent random variables: relative change in the
energy $E_\Delta$ (left) and in the augmented energy $\widetilde E_\Delta$
(right). Note the vertical scales, $10^{-5}$ and $10^{-11}$.}
\label{fig:num-hump2rv-energy}
\end{figure}

The relative energy in \cref{fig:num-hump2rv-energy} increases during
the computation, because the net energy flux through the zeroth-order
outflow boundaries is not zero once the water surface deforms over the
hump; the dissipation of the scheme is therefore visible only in the
augmented energy, which decreases monotonically and whose change is six
orders of magnitude smaller than that of the energy. Since no limiter is used, this
decrease is the dissipation of the ES flux alone, consistent
with the energy plots of \cite[Figure~16]{epshteyn2025energy}.

\subsection{Gaussian-shape hump with a random bottom field in ten random variables}
\label{sec:num-nrv}

The last example tests the scheme with a higher-dimensional random input.
The bottom is the deterministic hump of the accuracy test of
\cref{sec:num-accuracy-cmame}, perturbed by the truncation of a
Karhunen--Lo\`eve-type random field \cite{ghanem1991stochastic,le2010spectral}
(see \cite{mishra2012multilevel} for shallow water flows over random
topography of this type),
\begin{equation}\label{eq:num-humpnrv-bottom}
  B(x,y,\bm\xi)=0.5\exp\!\big(-25(x-1)^2-50(y-0.5)^2\big)+0.1
  +\sum_{i=1}^{d} f_i(x,y)\,\xi_i,
  \qquad
  f_i(x,y)=\cos\Big(\frac{i\pi x}{2}\Big)\cos(2i\pi y),
\end{equation}
with the deterministic uniform flow \eqref{eq:num-hump2rv-ic} as initial
data. The modes $f_i$ are
$1$-periodic in $y$, so that the random bottom is compatible with the
periodic boundaries, and orthogonal on $\Omega$ with
$\|f_i\|^2_{L^2(\Omega)}=\tfrac12$. The random variables are independent,
$\xi_i=0.1\,\eta_i$, where $\eta_i$ is supported on $[-1,1]$ with the
symmetric Beta density $\rho_i(\eta_i)\propto(1-\eta_i^2)^{i^2}$, i.e.,
$\alpha=\beta=i^2$ in the parameterization of \cref{sec:numerics}, so that
\begin{align*}
 \E[\xi_i]&=0,&
 \Var[\xi_i]&=\frac{10^{-2}}{2i^2+3},\\
 \E[B]&=0.5\exp\!\big(-25(x-1)^2-50(y-0.5)^2\big)+0.1,&
 \Var[B]&=10^{-2}\sum_{i=1}^{d}\frac{f_i(x,y)^2}{2i^2+3}.
\end{align*}
The variance of the modes decreases with their frequency, and the
concentration parameters are chosen such that the variance of the full
field converges: $\int_\Omega\Var[B]\dd\bm x=\frac{10^{-2}}2\sum_{i\ge1}(2i^2+3)^{-1}$
is a convergent series, and we take $d$ large enough that the truncated sum
reproduces it to within $10\%$, which gives $d=10$ (the first ten modes carry
$90.0\%$ of the variance of the full field; the pointwise maximum of
$\Std[B]$, attained where every $\cos^2$ factor equals one, is $0.0654$
against $0.0689$ for the full field). The PCE basis is the hyperbolic-cross
set \cite{shen2010sparse,dung2018hyperbolic} of degree two in the ten variables,
$\Lambda=\{\bm a\in\mathbb N_0^{10}:\prod_{i=1}^{10}(a_i+1)\le3\}$, which contains
the constant, the ten linear modes $\eta_i$ and the ten squares $\eta_i^2$,
$\kpce=21$; the basis functions are products of univariate orthonormal
Jacobi polynomials and the triple products in \eqref{eq:P-definition} are
products of the univariate ones. Since \eqref{eq:num-humpnrv-bottom} is affine in
$\bm\xi$, the initial data are represented exactly in this basis. The mesh is
chosen so that the shortest retained mode, $f_{10}$, has at least six cells
per period in each direction. The characteristic boundaries use the
initial data \eqref{eq:num-hump2rv-ic}
over the random bottom as the exterior reference state, so that at the
inflow boundary $x=0$ the incoming characteristic quantities are those of
the uniform flow. The parameters are
listed in \cref{tab:num-setup}.

\begin{figure}[htbp]
\centering
\includegraphics[width=\textwidth,trim={93 23 89 31},clip]{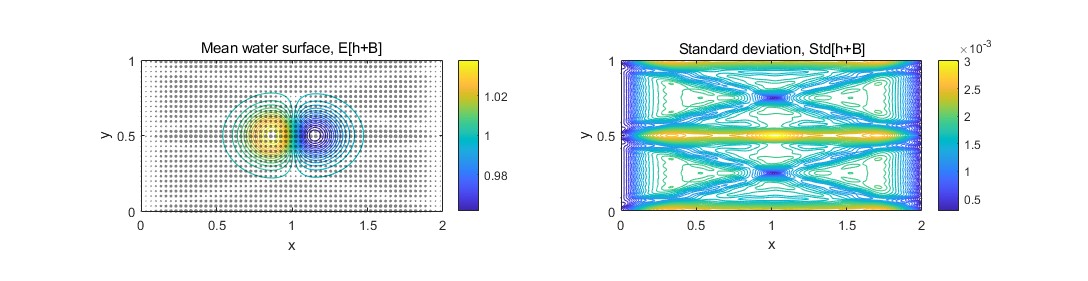}
\captionsetup{font=scriptsize}
\caption{Hump with the random bottom field in ten random variables at
$T=0.07$. Left: contours of the mean water surface $\E[w]$ with black disks
at every sixth plotting-grid point, area proportional to $\Std[w]$. Right:
contours of $\Std[w]$ (maximum $3.1\times10^{-3}$). $\mathbb P^2$ DG,
$\kpce=21$, $64\times64\times2$ triangles, no limiter.}
\label{fig:num-humpnrv}
\end{figure}

\Cref{fig:num-humpnrv} shows the mean and the standard deviation of
the water surface at $T=0.07$. The mean shows the same
elevation--depression pattern over the hump as in the one- and two-variable
versions of this flow
(\cref{sec:num-accuracy-cmame,sec:num-2rv}), since the random part of the
bottom has zero mean. The standard deviation, in contrast, is no longer
concentrated around the hump but inherits the spatial structure of the
random field: the water surface is flat at $t=0$, so at early times
$\partial_tw=-\partial_xq^x\approx0.3\,\partial_xB$ and $\Std[w]$ is
generated by the slope of the random bottom in the flow direction.
Its largest values lie along the lines $y=0$, $y=0.5$ and $y=1$, where all
modes $\cos(2i\pi y)$ attain their extrema, with the maximum
$3.1\times10^{-3}$ on $y=0.5$; the X-shaped lines of small
uncertainty are the common zeros of $\partial_xf_i$ of the leading modes,
and the nonlinear coupling with the hump deforms this pattern around $x=1$.
At the boundaries $x=0,2$, where $\partial_xf_i=0$ for every mode,
$\Std[w]\le7.7\times10^{-4}$. The response is dominated by the low modes:
the first five variables carry most of $\Var[w]$, and the same computation
with $d=5$ changes $\Std[w]$ by at most $12\%$ of its maximum ($5\%$ in
$L^2(\Omega)$). The quadratic PCE modes carry less than $10^{-4}$ of the
variance, and the linear basis ($\kpce=11$) reproduces $\Std[w]$ to a
relative difference of $10^{-4}$, reflecting the affine dependence of the
input on $\bm\xi$ and the small size of the uncertainty.

\begin{figure}[htbp]
\centering
\includegraphics[width=\textwidth]{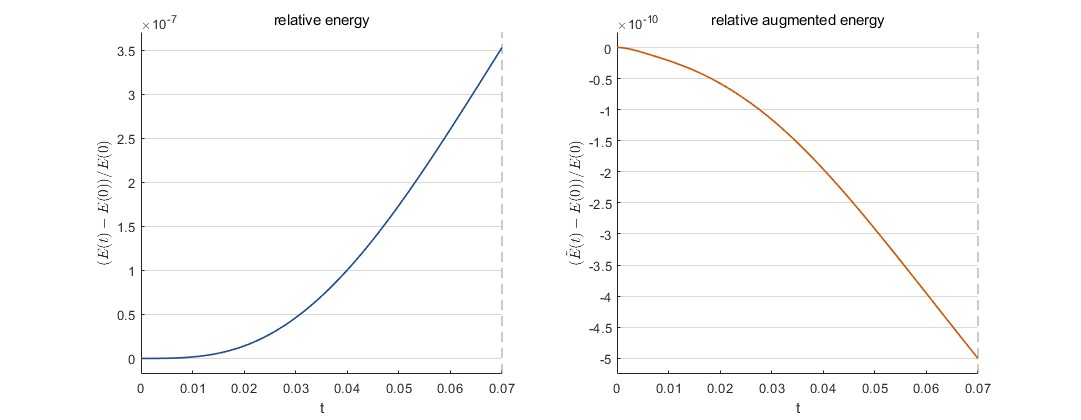}
\captionsetup{font=scriptsize}
\caption{Hump with the random bottom field in ten random variables: relative
change in the energy $E_\Delta$ (left) and in the augmented energy
$\widetilde E_\Delta$ (right).}
\label{fig:num-humpnrv-energy}
\end{figure}

As in \cref{sec:num-2rv}, the relative energy in
\cref{fig:num-humpnrv-energy} increases, by $3.5\times10^{-7}$ at $T$,
through the net energy flux across the characteristic boundaries, whereas
the augmented energy decreases monotonically to $-5.0\times10^{-10}$;
without a limiter this is the dissipation of the ES flux alone.
The example thus demonstrates that the product-matrix construction, the
coupled velocity projection and the ES fluxes apply without
modification to a $21$-mode system in ten non-identically distributed
random variables.

\section{Conclusion}\label{sec:conclusion}

In this paper we have developed a high-order discontinuous
Galerkin--stochastic Galerkin method for the two-dimensional shallow water
equations with random inputs. Starting from the hyperbolicity-preserving SG
formulation of \cite{dai2022hyperbolicity}, we derived a skew-symmetric SG
system that is hyperbolic under the same positive-definiteness condition on
the SG product matrix of the water height and that admits a companion
entropy balance law with respect to a modified entropy variable that is
linear in the water surface and in the projected velocities. This structure
was carried over to a DG discretization on triangular meshes with
elementwise projected velocity components. With compatible height, momentum,
and gravity face terms, the semi-discrete scheme is locally conservative for
the water height and the discharge, exactly well-balanced for stochastic lake-at-rest states,
and entropy conservative with central fluxes or entropy stable with
Lax--Friedrichs-type fluxes; the proof relies only on the symmetry of the
SG triple products and on the velocity projection, and the kinetic volume
terms cancel separately in both coordinate directions, so that the
directional asymmetry of the SG mixed fluxes does not interfere with the
entropy balance.

The fully discrete algorithm combines SSP Runge--Kutta time stepping with a
hyperbolicity-preserving time step restriction and linear scaling enforced
at every stage, a desingularized velocity projection, a characteristic TVB
limiter applied to equilibrium variables, and reflecting, periodic and
characteristic outflow boundary conditions implemented through ghost traces.
The augmented energy of \cite{epshteyn2025energy}, formulated here for the
DG weak form, retains the line integral of the
discrete boundary energy flux and thereby distinguishes the dissipation of
the scheme from the physical energy exchange through outflow boundaries.

The numerical experiments confirm the designed order of accuracy $k+1$ up
to degree three for the water height, the discharge and the velocity, the
well-balanced property to round-off over random topography, including in a
nearly dry region above a submerged plateau, the robustness of the scheme in
dry-bed and internal-wall dam-break configurations, the applicability to
several non-identically distributed random variables, including a truncated
Karhunen--Lo\`eve-type bottom field in ten variables, and the monotone decay
of the augmented energy in all examples with outflow boundaries, while the
energy itself reflects the exchange with the exterior. Future work includes
the treatment of moving-water steady states of the SWE, an efficient and
accurate treatment of dry and nearly dry regions, the extension of the
structure-preserving SG system to other nonlinear stochastic hyperbolic
systems, and data-driven methods that learn the physical flux of the SG
system from data with uncertainties.

\section*{Acknowledgements}
Funding information - Akil Narayan and Yinqian Yu were partially supported by NSF DMS-2207207. Akil Narayan was partially supported by AFOSR FA9550-23-1-0749. This research of Yekaterina Epshteyn is based upon work supported by
and while serving at the National Science Foundation. Any opinion,
findings, and conclusions or recommendations expressed in this material
are those of the author(s) and do not necessarily reflect the views
of the National Science Foundation.

\section*{Declaration of the use of generative SI}
During the preparation of this work the authors used generative SI tools
(ChatGPT and Claude) as writing and programming assistants. The tools were
used to polish the
language and the organization of the manuscript, to consolidate the
bibliography, and to assist in extending and refactoring parts of the
authors' MATLAB implementation. The
mathematical formulation, the derivation and proofs of the results, the
design of the numerical scheme, the numerical experiments and the
interpretation of the results are the authors' own. All text and code
produced with the assistance of the tools were reviewed, tested, and edited
by the authors, who take full responsibility for the content of the
publication.

\appendix
\section{Algorithmic details}
\label{sec:appendix-algorithm}

This appendix gives the details of the components of the fully discrete
algorithm summarized in \cref{sec:implementation}. For an
element $K\in\Th$ we denote by $\hK$ its local mesh size, defined in
\eqref{eq:num-cfl} as the ratio of twice the area to the longest edge; the
letter $h$ is reserved for the water height.

\subsection{Desingularized velocity projection}
\label{sec:desingularization}

The velocity projection \eqref{eq:velocity-projection} requires, on every
element, the solution of a linear system whose matrix is the mass matrix of
$\Pmat(\hD)$. When the water height is very small, as in the nearly dry
regions of the examples in \cref{sec:num-plateau,sec:num-circular}, this
matrix becomes ill-conditioned and the projected velocities may take
unphysically large values. To address this issue we adopt the
desingularization procedure introduced in \cite{kurganov2007second} for the
deterministic SWE, whose stochastic variant was used in
\cite{dai2021hyperbolicity,dai2022hyperbolicity,dai2024energy,epshteyn2025energy}
and in the one-dimensional DG--SG scheme \cite{epshteyn2026structure}.

At every volume quadrature point of an element $K$, let
$\Pmat(\hD)=W\transpose DW$ be the eigenvalue decomposition of the
symmetric matrix $\Pmat(\hD)$, with $W$ orthogonal and $D=\mathrm{diag}(\lambda_1,\dots,\lambda_{\kpce})$.
For a given positive constant $\epsilon$ we define
\begin{align}\label{eq:desing-eigs}
    D^{\mathrm{cor}} \coloneqq \mathrm{diag}(\lambda_1^{\mathrm{cor}},\dots,\lambda_{\kpce}^{\mathrm{cor}}),
    & &
    \lambda_i^{\mathrm{cor}} \coloneqq \frac{\sqrt{\lambda_i^4+\max(\lambda_i^4,\epsilon^4)}}{\sqrt{2}\,\lambda_i},
\end{align}
so that $\lambda_i^{\mathrm{cor}}=\lambda_i$ when $\lambda_i\ge\epsilon$
and $\lambda_i^{\mathrm{cor}}\ge\epsilon/\sqrt2$ otherwise, and we set
$\Pcor(\hD)\coloneqq W\transpose D^{\mathrm{cor}}W$. The desingularized
velocities are obtained by replacing $\Pmat(\hD)$ with $\Pcor(\hD)$ in the
velocity projection: on every $K\in\Th$,
\begin{subequations}\label{vel_update}
\begin{align}
 \int_K[\Pcor(\hD)\uD-\qxD]\cdot\bm\zeta^x\dd\bm x&=0,
 &&\forall\bm\zeta^x\in[\Vk]^{\kpce},\\
 \int_K[\Pcor(\hD)\vD-\qyD]\cdot\bm\zeta^y\dd\bm x&=0,
 &&\forall\bm\zeta^y\in[\Vk]^{\kpce}.
\end{align}
\end{subequations}
To maintain the consistency relation \eqref{eq:velocity-projection} between
the discharges and the velocities that enter the scheme, the discharges are
subsequently recomputed from the desingularized velocities with the original
product matrix,
\begin{subequations}\label{eq:desing-recompute}
\begin{align}
 \int_K[\Pmat(\hD)\uD-\qxD]\cdot\bm\zeta^x\dd\bm x&=0,
 &&\forall\bm\zeta^x\in[\Vk]^{\kpce},\\
 \int_K[\Pmat(\hD)\vD-\qyD]\cdot\bm\zeta^y\dd\bm x&=0,
 &&\forall\bm\zeta^y\in[\Vk]^{\kpce}.
\end{align}
\end{subequations}
Away from nearly dry regions, where all eigenvalues satisfy
$\lambda_i\ge\epsilon$, the two steps reduce to the exact projection
\eqref{eq:velocity-projection} and the discharges are unchanged. Unless
stated otherwise, we take $\epsilon=\hK$ on each element $K$. The same
eigenvalue correction is applied to $\Pmat(\hD)$ when the normal SG Jacobian
\eqref{eq:sg-normal-jacobian} is evaluated on the faces for the wave speed
$\alpha_1$.

\subsection{Hyperbolicity-preserving criterion}\label{sec_hpc}

The hyperbolicity of the skew-symmetric SG system \eqref{eq:sg-skew}
requires that $\Pmat(\hhat)$ be positive definite
(Proposition~\ref{prop:hyperbolicity}). For the semi-discrete scheme this corresponds to
the condition $\Pmat(\hD)\succ0$ at every quadrature point, in the volume and
on the faces of every element, at which the scheme is evaluated. Even if the
condition holds at the current time, additional conditions must be imposed
to ensure that it holds at the next time step. A sufficient and computable
condition is provided by \cite[Theorem~3.2]{dai2021hyperbolicity}: if the
truncated water height $h_{\Delta,\Lambda}\coloneqq \sum_{j=1}^{\kpce}\widehat h_{\Delta,j}\phi_j$ is
positive at the nodes of an $\Nxi$-point positive quadrature rule in the
stochastic space that integrates exactly all products of three
polynomials of the PCE space $\mathbb P_\Lambda$, then $\Pmat(\hD)$ is positive definite.
Based on this, we define the hyperbolicity-preserving property of a fully
discrete scheme as follows.

\begin{definition}[Hyperbolicity-preserving fully discrete DG--SG scheme]
\label{def:hp-scheme}
Let $\{\bm x_i\}_{i=1}^{N_q} \subset \mathbb{R}^2$ be the set of physical quadrature points,
in the volume and on the faces, involved in the DG--SG scheme, and let
$\{\bm\xi^{(m)}\}_{m=1}^{\Nxi}$ be the nodes of a positive stochastic quadrature
rule that integrates exactly all products of three polynomials of
$\mathbb P_\Lambda$. Suppose that at the current time level $t^n$
\begin{equation}\label{pp_creterion_assumption}
    h_{\Delta,\Lambda}(\bm x_i,t^n,\bm\xi^{(m)}) \coloneqq \sum_{j=1}^{\kpce} \widehat h_{\Delta,j}(\bm x_i,t^n)\phi_j(\bm\xi^{(m)}) >0,
    \quad \forall\, i=1,\dots,N_q,\ m=1,\dots,\Nxi.
\end{equation}
The fully discrete DG--SG scheme is said to preserve the hyperbolicity
condition if
\begin{equation}\label{pp_creterion_condition}
    h_{\Delta,\Lambda}(\bm x_i,t^{n+1},\bm\xi^{(m)})>0, \quad \forall\, i=1,\dots,N_q,\ m=1,\dots,\Nxi.
\end{equation}
\end{definition}

Suppose that \eqref{pp_creterion_assumption} holds and that a forward Euler
step is applied. To enforce \eqref{pp_creterion_condition} at the next time
level we use one of the following two strategies, which are the
two-dimensional analogs of the strategies of the one-dimensional scheme
\cite{epshteyn2026structure}.

\begin{remark}[Mean-positivity-preserving time step and linear scaling]
\label{rem:linear-scaling}
Let $\Phi(\bm\xi)\coloneqq(\phi_1(\bm\xi),\dots,\phi_{\kpce}(\bm\xi))\transpose$ and
let $\overline{\hD}_K$ denote the cell average of $\hD$ on $K$. By the local
conservation property (\cref{thm:local-conserve}), a forward Euler step
updates the cell average as
$\overline{\hD}_K^{\,n+1}=\overline{\hD}_K^{\,n}-\frac{\Delta t^n}{|K|}\oint_{\partial K}\fq^{\,K}\dd s$.
Hence the cell-average water height stays positive at every stochastic node
if, on each element $K$,
\begin{equation}\label{eq:hp-dt-mean}
\Delta t^n < \Delta t_K^n \coloneqq \min_{m=1,\dots,\Nxi}
 \left|\frac{|K|\,(\overline{\hD}_K^{\,n})\transpose\Phi(\bm\xi^{(m)})}
 {\bigl(\oint_{\partial K}\fq^{\,K}\dd s\bigr)\transpose\Phi(\bm\xi^{(m)})}\right| .
\end{equation}
Positivity at all quadrature points $\mathcal S_K$ of the element is then restored
by the linear scaling
\begin{equation}\label{eq:hp-linear-scaling}
\begin{aligned}
 \tilde h^{n}_{\Delta,\Lambda}(\bm x,\bm\xi^{(m)})
 &=\theta_K\bigl(h^{n}_{\Delta,\Lambda}(\bm x,\bm\xi^{(m)})-\overline{h}^{\,n}_{\Delta,\Lambda,K}(\bm\xi^{(m)})\bigr)
 +\overline{h}^{\,n}_{\Delta,\Lambda,K}(\bm\xi^{(m)}),\\
 \theta_K&=\min_{m=1,\dots,\Nxi}\min\Bigl\{1,\ \frac{\overline{h}^{\,n}_{\Delta,\Lambda,K}(\bm\xi^{(m)})}
 {\overline{h}^{\,n}_{\Delta,\Lambda,K}(\bm\xi^{(m)})-\mu_K(\bm\xi^{(m)})}\Bigr\},
\end{aligned}
\end{equation}
where $\overline{h}^{\,n}_{\Delta,\Lambda,K}(\bm\xi^{(m)})=(\overline{\hD}_K^{\,n})\transpose\Phi(\bm\xi^{(m)})$
and $\mu_K(\bm\xi^{(m)})=\min_{\bm x\in\mathcal S_K}h^{n}_{\Delta,\Lambda}(\bm x,\bm\xi^{(m)})$.
In the implementation the numerator of $\theta_K$ is reduced by a small
positive floor, so that the scaled water height stays strictly positive.
The scaling
acts on the PCE coefficient vector as a whole, so it preserves the cell
averages of all stochastic modes. This linear scaling extends the
positivity-preserving limiter of \cite{xing2010positivity} (see
\cite{zhang2012maximum,xing2013positivity} for triangular meshes) to the
DG--SG formulation; we also refer to it as the hyperbolicity-preserving
scaling. In the numerical experiments we use \eqref{eq:hp-dt-mean} with
a safety factor $0.9$ together with \eqref{eq:hp-linear-scaling}.
\end{remark}

\begin{remark}[Pointwise restriction]\label{rem:pointwise-restriction}
Alternatively, one can impose a more direct restriction on $\Delta t$ at
every quadrature point. Writing the forward Euler update of the height
coefficients as $\hD^{n+1}=\hD^{n}+\Delta t^n\mathcal L_h(\widehat U^n)$,
positivity at every quadrature point and every stochastic node holds if
\begin{equation}\label{eq:hp-dt-pointwise}
\Delta t^n < \Delta t_K^n \coloneqq \min_{i=1,\dots,|\mathcal S_K|}\ \min_{m=1,\dots,\Nxi}
 \left|\frac{\hD^{n}(\bm x_i)\transpose\Phi(\bm\xi^{(m)})}
 {\mathcal L_h(\widehat U^n)(\bm x_i)\transpose\Phi(\bm\xi^{(m)})}\right| ,
\end{equation}
where $|\mathcal S_K|$ denotes the number of physical quadrature points of the
element. This restriction needs no scaling, but it is typically more
restrictive than \eqref{eq:hp-dt-mean}.
\end{remark}

The strategies of Remarks~\ref{rem:linear-scaling}
and~\ref{rem:pointwise-restriction} are not the only ways to enforce \eqref{pp_creterion_condition}; we simply adopt
them here. In our numerical experiments we did not observe extremely small
time steps caused by the hyperbolicity-preserving restriction.

\subsection{Troubled-cell indicator and characteristic TVB limiter}
\label{sec_slopelimiting}

In order to suppress spurious oscillations near discontinuities, a slope
limiting procedure is applied after every Runge--Kutta stage. We follow the
standard procedure of RKDG methods \cite{cockburn2001runge,qiu2005comparison},
which consists of two steps: first, the troubled cells, i.e., the elements
in which limiting is needed, are identified; second, a slope limiter is
applied to the solution in those cells.

For the first step we extend the troubled-cell indicator of Fu and Shu
\cite{fu2017new}, in the scaled form used in \cite{fu2022high}, to the
stochastic setting. For a piecewise polynomial $p\in \Vk$ the
indicator on the element $K$ is
\begin{equation}\label{eq:troubled-indicator}
    \mathcal I_K(p) = \frac{\sum_{K' \in \omega(K)}\bigl| \overline{\overline{p}}_{K'} - \overline{p}_K\bigr|}{\overline{p}_{\max} - \overline{p}_{\min}},
\end{equation}
where $\omega(K)$ is the set of elements sharing an edge with $K$,
$\overline{p}_K$ is the cell average of $p$ on $K$,
$\overline{\overline{p}}_{K'}$ is the cell average over $K$ of the
polynomial $p|_{K'}$ extended from the neighbor $K'$ to $K$, and
$\overline{p}_{\max}$ and $\overline{p}_{\min}$ are the largest and smallest
cell averages over the domain. For a given tolerance $\mathrm{tol}>0$, the
element $K$ is marked as a troubled cell if $\mathcal I_K(p)>\mathrm{tol}$.
If all cell averages of $p$ coincide, as for a lake at rest, no element is
marked. In the
DG--SG scheme we take the spatial linear ($\mathbb P^1$) projection of the
stochastic mean of the water surface,
$\E[h_{\Delta,\Lambda}+B_{\Delta,\Lambda}]=\widehat h_{\Delta,1}+\widehat B_{\Delta,1}$,
as the indicator function $p$; on boundary elements the neighbors are the
ghost elements defined by the boundary conditions of \cref{sec:boundary}.

In the second step we apply the characteristic-wise TVB limiter for
triangular meshes of Cockburn and Shu \cite[Section~3.4.2]{cockburn1998runge}
to the equilibrium variables
$(\hD\transpose+\BD\transpose,(\qxD)\transpose,(\qyD)\transpose)\transpose$
in the troubled cells. The complete modal polynomial is first projected onto
$\mathbb P^1$; the differences between the values of the $\mathbb P^1$
polynomial at the three edge midpoints and the cell average are limited,
with the parameter $\nu>1$ of \cite{cockburn1998runge} and the TVB constant
$\Mtvb\ge0$, in the local characteristic variables of the full SG normal Jacobian
\eqref{eq:sg-normal-jacobian} evaluated at the cell average, and the limited
midpoint values are converted back to modal coefficients. Following
\cite{cockburn1998runge}, the original $\mathbb P^k$ polynomial is retained
whenever its $\mathbb P^1$ projection is not modified by the limiter, so
that the limiter does not reduce the accuracy in smooth regions. The use of
the water surface $\hD+\BD$, rather than $\hD$ alone, is compatible with the
well-balanced structure: the limiter leaves the stochastic lake at rest
unchanged, whereas limiting $\hD$ over a variable bottom would generate
spurious waves. After limiting, the velocities are recomputed by the
projection of \cref{sec:desingularization}.

\subsection{Boundary conditions}
\label{sec:boundary}

Every boundary face carries the same numerical fluxes
\eqref{eq:flux-h}, \eqref{eq:flux-qu}, \eqref{eq:flux-qv} as an interior
face, with the exterior trace $(\cdot)^+$ replaced by a ghost state that
encodes the boundary condition. This face-trace mechanism supports the
following boundary conditions.
\begin{itemize}[leftmargin=2em]
\item \emph{Periodic:} the exterior trace is the trace of the paired
  element on the opposite side of the domain; matching boundary partitions
  are required.
\item \emph{Reflecting wall:} the water height and the bottom are copied,
  and the normal components of the discharge and of the projected velocity
  are reflected,
  \[
   \bm q^+=\bm q^- -2(\bm q^-\cdot\bm n)\bm n,\qquad
   \bm u^+=\bm u^- -2(\bm u^-\cdot\bm n)\bm n,
  \]
  componentwise in SG space, so that the normal mass flux through the wall
  vanishes.
\item \emph{Zeroth-order outflow:} all exterior traces equal the interior
  traces.
\item \emph{Characteristic outflow boundary:} following the
  characteristic treatment of \cite{atkins1997continued}, the SG normal Jacobian
  \eqref{eq:sg-normal-jacobian} is evaluated at the interior trace and
  diagonalized; the amplitudes of the incoming characteristic fields are
  prescribed from a far-field state (in the examples, the lake at rest
  $w=1$, $q^x=q^y=0$), and the outgoing amplitudes are extrapolated from the
  interior. The ghost state is assembled from the resulting characteristic
  amplitudes.
\end{itemize}
On domains with several boundary components, such as the basin with an
internal dam of \cref{sec:num-curved}, the boundary faces are tagged, so
that different conditions can be imposed on different components. The ghost
states are also used as neighbors in the troubled-cell indicator and in the
limiter, and in the boundary energy functional of \cref{sec:augmented-entropy}.

\subsection{SSP Runge--Kutta time integration with adaptive time step}
\label{sec:time-stepping}

Let $\mathcal L(\widehat U)$ denote the spatial DG--SG operator of
\cref{sec:dg-scheme}, including the velocity projection.
Here the time derivative $(\hD)_t$ in the momentum equations
\eqref{eq:dg-qx} and \eqref{eq:dg-qy} is first obtained from the height
equation \eqref{eq:dg-h} and then substituted, so that $\mathcal L$ is an
explicit operator. We integrate in
time with the third-order strong-stability-preserving Runge--Kutta method
(SSP-RK3) \cite{shu1988efficient,gottlieb2001strong},
\begin{subequations}\label{eq:ssprk3}
\begin{align}
 U^{(1)}&=U^n+\Delta t\,\mathcal L(U^n),\\
 U^{(2)}&=\tfrac34U^n+\tfrac14
   [U^{(1)}+\Delta t\,\mathcal L(U^{(1)})],\\
 U^{n+1}&=\tfrac13U^n+\tfrac23
   [U^{(2)}+\Delta t\,\mathcal L(U^{(2)})].
\end{align}
\end{subequations}
The time step restriction of \cref{sec_hpc} is derived for a forward Euler
step. Since \eqref{eq:ssprk3} is a convex combination of forward Euler
steps, the restriction carries over to each stage, provided that it is
satisfied with the stage state \cite{gottlieb2001strong}. We therefore
employ the adaptive strategy of \cite{chertock2015well-adaptiveRK}: the
admissible time step is recomputed at every intermediate stage from the
CFL condition with the spectral radius of the normal SG Jacobian and from
the hyperbolicity-preserving bound; if a stage requires a smaller step than
the one in use, the step is rejected and restarted from $U^n$ with the
smaller step. The limiter of \cref{sec_slopelimiting}, the linear scaling
of \cref{sec_hpc} and the velocity projection are applied after every
stage. We note that \cref{thm:ec-es} concerns the semi-discrete scheme; a
fully discrete entropy inequality for explicit Runge--Kutta time stepping
can be enforced, for instance, with relaxation Runge--Kutta methods
\cite{ketcheson2019relaxation,ranocha2020relaxation}, which we do not pursue
here. The complete procedure is summarized in \cref{alg:fully-discrete}.

\subsection{Augmented energy}
\label{def_aug_energy}
\label{sec:augmented-entropy}

\Cref{thm:ec-es} is stated for periodic domains, where every face is
interior. In the presence of physical boundaries, the energy
inside the domain also changes through the energy flux across the boundary.
Such a change is not a deficiency of the numerical scheme but reflects the
dynamics of the model: for instance, with outflow boundaries energy may leave
the domain with outgoing waves, or enter it with an inflow, and the total
energy may then increase; see \cite{epshteyn2025energy} for a discussion in
the finite-volume setting.
For periodic or impermeable boundaries no additional energy enters the
domain, and $E_\Delta$ itself is non-increasing.

To separate the dissipation of the scheme from the boundary exchange, we use
the augmented energy introduced for finite-volume schemes in
\cite{epshteyn2025energy} (see also \cite{epshteyn2026structure} for the
one-dimensional DG--SG scheme), which adds the boundary energy flux back to
the energy; what is specific to the present setting is the boundary
functional generated by the two-dimensional DG weak form. At the continuous
level,
integrating the companion law \eqref{eq:continuum-entropy-law} over
$\Omega$ and applying the divergence theorem gives
\begin{align}
 \frac{\dd}{\dd t}\int_\Omega E\dd\bm x
 &=-\int_\Omega
   [(H^x)_x+(H^y)_y]\dd\bm x\notag\\
 &=-\int_{\partial\Omega}
   (n_xH^x+n_yH^y)\dd s,
\label{eq:2d-integration-by-parts}
\end{align}
with $\bm n$ the outward normal of $\Omega$. The continuous augmented
energy
\begin{equation}\label{eq:continuum-augmented-entropy}
 \widetilde E(t)
 =\int_\Omega E(t)\dd\bm x
 +\int_0^t\int_{\partial\Omega}
   \bm H\cdot\bm n\dd s\,\dd\tau,
 \qquad
 \widetilde E(0)=E(0),
\end{equation}
is therefore constant for smooth solutions even when the energy inside the
domain is not.

The same calculation can be carried out for the DG weak form. Repeating the
proof of \cref{thm:ec-es} with physical boundary faces, all interior-face
contributions cancel or produce the nonnegative dissipation in
\eqref{eq:entropy-rate}, and the boundary faces
$e\in\Ehb$ contribute the face terms of \eqref{eq:dg-h}, \eqref{eq:dg-qx}
and \eqref{eq:dg-qy} tested with the interior traces of the entropy
variables.
Denoting the interior trace by the superscript $\mathrm{in}$, the boundary
entropy functional generated by the numerical fluxes is
\begin{equation}
\begin{aligned}
 H^*_{\bm n}
 ={}&
 \fqu^{\,K}\cdot\uD^{\mathrm{in}}
 +\fqv^{\,K}\cdot\vD^{\mathrm{in}}
 -\frac12\Pmat(\fq^{\,K})\uD^{\mathrm{in}}\cdot\uD^{\mathrm{in}}
 -\frac12\Pmat(\fq^{\,K})\vD^{\mathrm{in}}\cdot\vD^{\mathrm{in}}
 +g\fq^{\,K}\cdot(\hD^{\mathrm{in}}+\BD^{\mathrm{in}})
 +\mathcal C^*_{g,\bm n} ,\\
 \mathcal C^*_{g,\bm n}
 ={}&\frac g2\Big(
 \Pmat(\hD^{\mathrm{in}})\jump{\hD+\BD}_K \Big)\cdot
 \bigl(n_x\uD^{\mathrm{in}}+n_y\vD^{\mathrm{in}}\bigr),
\end{aligned}
\label{eq:numerical-boundary-entropy-flux}
\end{equation}
where the fluxes are evaluated from the interior trace and the ghost state
of \cref{sec:boundary} with the outward normal, and $\mathcal C^*_{g,\bm n}$ is
the contribution of the gravity face corrections in
\eqref{eq:dg-qx}--\eqref{eq:dg-qy} on the boundary face, tested with the
interior velocity traces $\uD^{\mathrm{in}}$ and $\vD^{\mathrm{in}}$. This term is proportional to
$\Pmat(\jump{\hD+\BD})$ and therefore vanishes whenever the exterior and
interior water surfaces agree.

For the zeroth-order outflow condition, $U^+=U^-$ and
$\mathcal C^*_{g,\bm n}=0$; the numerical fluxes are consistent, so that
\begin{equation}\label{eq:outflow-entropy-flux}
\begin{aligned}
 H^*_{\bm n}
 ={}&n_x\left[
 \frac12\Pmat(\qxD)\uD\cdot\uD
 +\frac12\Pmat(\qxD)\vD\cdot\vD
 +g\qxD\cdot(\hD+\BD)\right]
 \\
 &+n_y\left[
 \frac12\Pmat(\qyD)\uD\cdot\uD
 +\frac12\Pmat(\qyD)\vD\cdot\vD
 +g\qyD\cdot(\hD+\BD)\right],
\end{aligned}
\end{equation}
which is exactly the normal component of the SG entropy flux
\eqref{eq:sg-entropy-flux}. The discrete augmented energy is defined as
\begin{equation}\label{eq:discrete-augmented-entropy}
 \widetilde E_\Delta(t)
 =E_\Delta(t)
 +\int_0^t\sum_{e\in\Ehb}\int_e
   H^*_{\bm n}(\tau)\dd s\,\dd\tau ,
 \qquad \widetilde E_\Delta(0)=E_\Delta(0).
\end{equation}
By construction, for the semi-discrete EC scheme
$\dd\widetilde E_\Delta/\dd t=0$, and for the ES scheme
\begin{equation}\label{eq:augmented-entropy-stability}
 \frac{\dd\widetilde E_\Delta}{\dd t}
 =-\mathcal D_\Delta(t)\le0,
\end{equation}
where $\mathcal D_\Delta\ge0$ is the sum of the three interior-face
quadratic forms in \eqref{eq:entropy-rate}. Since
$H^*_{\bm n}$ is built from the ghost states themselves, this identity holds
for any boundary closure: it restates the interior estimate of
\cref{thm:ec-es} and is not a stability statement for the boundary
treatment. The augmented energy is thus a diagnostic that isolates the
interior dissipation of the scheme from the boundary exchange, and its
monotonicity is the property that we monitor numerically on domains with
outflow boundaries. In the fully
discrete scheme the time integral in \eqref{eq:discrete-augmented-entropy}
is accumulated with the same Runge--Kutta stage weights as the solution
update.

\begin{remark}
For characteristic outflow boundaries, walls, or prescribed inflow, the
general functional \eqref{eq:numerical-boundary-entropy-flux} must be used
rather than the physical outflow flux \eqref{eq:outflow-entropy-flux}, so
that the diagnostic is exactly consistent with the ghost state and with the
DG residual. In particular, at a reflecting wall the ghost state of
\cref{sec:boundary} produces a vanishing normal mass flux, so that no
physical energy is exchanged through the wall.
\end{remark}

\section{Discretization common to all numerical examples}
\label{sec:appendix-discretization}

This appendix lists the discretization used in all examples of
\cref{sec:numerics}; the parameters that differ between the examples are
given in \cref{tab:num-setup}.
\begin{itemize}[leftmargin=2em]
\item \emph{Stochastic space:} orthonormal PCE with the number of terms
  $\kpce$ listed in \cref{tab:num-setup}; for one random variable these
  are the polynomials of degree at most $\kpce-1$, and the bases of
  \cref{sec:num-2rv,sec:num-nrv} are described there. For one random
  variable the triple products in \eqref{eq:P-definition} are computed
  with a 120-point Gauss rule.
\item \emph{Physical space:} modal DG with the total-degree $\mathbb P^k$
  monomial basis on affine triangles, on either an unstructured DistMesh
  triangulation \cite{persson2004simple} with target edge length $\Delta_0$ or
  a structured grid of split quadrilaterals. The degree $k$ and the mesh are
  listed in \cref{tab:num-setup}.
  Volume integrals are evaluated with a tensor-product Gauss rule mapped to
  the triangle, and face integrals with a Gauss--Legendre rule; both
  integrate all polynomial integrands of the scheme exactly (degree up to
  $3k$ in the volume and $4k$ on the faces).
\item \emph{Numerical flux:} the ES fluxes
  \eqref{eq:flux-h}, \eqref{eq:flux-qu}, \eqref{eq:flux-qv} with $\alpha_1$
  the spectral radius of the normal SG Jacobian
  \eqref{eq:sg-normal-jacobian} evaluated from the two face traces.
  The entropy conservative scheme with the central fluxes
  ($\alpha_1\equiv0$) is an intermediate result in \cref{thm:ec-es} and is not tested separately; its conservation
  property could only be examined in a smooth periodic problem designed for
  that purpose.
\item \emph{Time stepping:} SSP-RK3 \eqref{eq:ssprk3} with
  \begin{equation}\label{eq:num-cfl}
    \Delta t=\frac{\Ccfl}{\max_{K\in\Th}\max_{\text{faces, nodes}}(\alpha_1/\hK)},
    \qquad \hK=\frac{2|K|}{\text{longest edge of }K}.
  \end{equation}
  Where the hyperbolicity-preserving scaling of
  \cref{sec_hpc} is used (\cref{tab:num-setup}), the step is further restricted by its criterion
  (linear scaling, safety factor $0.9$), and the stages are restarted as in
  \cref{alg:fully-discrete}. Steps are shortened so that the output times
  are hit exactly.
\item \emph{Limiter:} the characteristic TVB limiter of
  \cref{sec_slopelimiting} ($\nu=1.5$, $\Mtvb=0$) in the troubled cells
  detected by the indicator applied to $\E[w]$ with the threshold
  $\mathrm{tol}$ listed in \cref{tab:num-setup}. The full $\mathbb P^k$ polynomial
  is kept whenever its $\mathbb P^1$ projection is not modified.
\item \emph{Velocity:} the coupled projection \eqref{eq:velocity-projection},
  either exact or with the desingularization of \eqref{vel_update}; the
  choice and the value of $\epsilon$ are listed in \cref{tab:num-setup}.
\item \emph{Boundary conditions:} listed in \cref{tab:num-setup}. The
  characteristic outflow boundaries of
  \cref{sec:num-hump,sec:num-perturbation,sec:num-plateau} use the far-field
  state $w=1$, $q^x=q^y=0$. Outflow boundaries are implemented through ghost
  states (\cref{sec:boundary}), so that every boundary face carries the same
  numerical flux as an interior face.
\item \emph{Initial data:} $w$, $B$ and the discharges are $L^2$-projected
  onto the DG--PCE space, with $h$ obtained from the samples of $w-B$.
  Denoting the projection by $\Pi_\Delta$, this gives
  $\Pi_\Delta h+\Pi_\Delta B=\Pi_\Delta w$, so that the lake at rest is
  represented exactly.
\end{itemize}

\section{Proof of \texorpdfstring{Theorem~\ref{thm:ec-es}}{the entropy conservation and stability theorem}}
\label{sec:appendix-proof}

This appendix gives the complete proof of \cref{thm:ec-es}, which is
outlined in \cref{sec:dg-scheme}.

\begin{proof}[Proof of \cref{thm:ec-es}]
The modified entropy variable
\eqref{eq:sg-modified-entropy-variable} suggests the choice of test
functions
\[
  \widehat\psi_{\Delta,j}=g(\widehat h_{\Delta,j}+\widehat B_{\Delta,j}),
  \qquad
  \widehat\zeta^{x}_{\Delta,j}=\widehat u_{\Delta,j},
  \qquad
  \widehat\zeta^{y}_{\Delta,j}=\widehat v_{\Delta,j},
  \qquad j=1,\dots,\kpce,
\]
all of which belong to $\Vk$. We insert them into \eqref{eq:dg-h},
\eqref{eq:dg-qx}, and \eqref{eq:dg-qy}, sum over $j$, and add the three
resulting identities. Collecting all time-derivative terms on the left-hand
side and the spatial terms on the right-hand side, we write the result as
\[
  \mathrm{LHS}=R_1+R_2,
\]
where $R_1$ collects the spatial terms without gravitational contributions
and $R_2$ the spatial terms with gravitational contributions.

\smallskip\noindent\textbf{Time-derivative terms.}
Using \eqref{eq:velocity-projection}, Lemma~\ref{lem:comm-dot}, and the
product rule, the time-derivative terms read
\begin{align}\label{eq:LHS-first}
  \mathrm{LHS}
  &=\sum_{K\in\Th}\int_Kg(\hD)_t\cdot(\hD+\BD)\dd\bm x
  +\sum_{K\in\Th}\int_K(\qxD)_t\cdot \uD\dd\bm x
  +\sum_{K\in\Th}\int_K(\qyD)_t\cdot \vD\dd\bm x\notag\\
  &\quad-\frac12\sum_{K\in\Th}\int_K\Pmat((\hD)_t)\uD\cdot\uD\dd\bm x
  -\frac12\sum_{K\in\Th}\int_K\Pmat((\hD)_t)\vD\cdot \vD\dd\bm x .
\end{align}
Here the velocity projection \eqref{eq:velocity-projection} gives, in the
weak sense (i.e., after multiplication by a test function in
$[\Vk]^{\kpce}$ and integration over an element),
\[
  \qxD=\Pmat(\hD)\uD,
  \qquad
  \qyD=\Pmat(\hD)\vD.
\]
Differentiating these relations in time and testing with $\uD$ and $\vD$,
which are admissible test functions, gives, on every element $K$,
\begin{align*}
  \int_K(\qxD)_t\cdot\uD\dd\bm x
  &=\int_K\Pmat((\hD)_t)\uD\cdot\uD\dd\bm x
    +\int_K\Pmat(\hD)(\uD)_t\cdot\uD\dd\bm x,\\
  \int_K(\qyD)_t\cdot \vD\dd\bm x
  &=\int_K\Pmat((\hD)_t)\vD\cdot \vD\dd\bm x
    +\int_K\Pmat(\hD)(\vD)_t\cdot \vD\dd\bm x.
\end{align*}
By Lemma~\ref{lem:comm-dot}, in the same weak sense (i.e., after integration
over $K$),
\[
  \Pmat(\hD)(\uD)_t\cdot\uD=\qxD\cdot(\uD)_t,
  \qquad
  \Pmat(\hD)(\vD)_t\cdot \vD=\qyD\cdot(\vD)_t.
\]
Therefore \eqref{eq:LHS-first} becomes
\begin{equation}
\begin{aligned}
  \mathrm{LHS}
  &=\sum_{K\in\Th}\int_Kg(\hD)_t\cdot \hD\dd\bm x
  +\sum_{K\in\Th}\int_Kg(\hD)_t\cdot \BD\dd\bm x\\
  &\quad+\sum_{K\in\Th}\int_K\frac12(\qxD)_t\cdot \uD\dd\bm x
  +\sum_{K\in\Th}\int_K\frac12\qxD\cdot (\uD)_t\dd\bm x\\
  &\quad+\sum_{K\in\Th}\int_K\frac12(\qyD)_t\cdot \vD\dd\bm x
  +\sum_{K\in\Th}\int_K\frac12\qyD\cdot (\vD)_t\dd\bm x\\
  &=\frac{\dd}{\dd t}E_\Delta,
\end{aligned}
\end{equation}
where the last equality is the time derivative of \eqref{eq:discrete-energy}
and uses that $\BD$ is independent of time.

\smallskip\noindent\textbf{Non-gravitational spatial terms.}
The terms of $R_1$ are the transport terms of the two momentum equations and
the skew-symmetric corrections:
\begin{equation}
\begin{aligned}
  R_1
  &=\sum_{K\in\Th}\sum_{j=1}^{\kpce}\int_K
  \bigl(\Pmat(\qxD)\uD\bigr)_j(\widehat u_{\Delta,j})_x\dd\bm x
  +\sum_{K\in\Th}\sum_{j=1}^{\kpce}\int_K
  \bigl(\Pmat(\qyD)\uD\bigr)_j(\widehat u_{\Delta,j})_y\dd\bm x\\
  &\quad+
  \sum_{K\in\Th}\sum_{j=1}^{\kpce}\int_K
  \bigl(\Pmat(\qxD)\vD\bigr)_j(\widehat v_{\Delta,j})_x\dd\bm x
  +\sum_{K\in\Th}\sum_{j=1}^{\kpce}\int_K
  \bigl(\Pmat(\qyD)\vD\bigr)_j(\widehat v_{\Delta,j})_y\dd\bm x\\
  &\quad-
  \frac12\sum_{K\in\Th}\sum_{j=1}^{\kpce}\int_K
  \bigl(\Pmat(\qxD)(\widehat u_{\Delta,j}\uD)_x\bigr)_j\dd\bm x
  -\frac12\sum_{K\in\Th}\sum_{j=1}^{\kpce}\int_K
  \bigl(\Pmat(\qyD)(\widehat u_{\Delta,j}\uD)_y\bigr)_j\dd\bm x\\
  &\quad-
  \frac12\sum_{K\in\Th}\sum_{j=1}^{\kpce}\int_K
  \bigl(\Pmat(\qxD)(\widehat v_{\Delta,j}\vD)_x\bigr)_j\dd\bm x
  -\frac12\sum_{K\in\Th}\sum_{j=1}^{\kpce}\int_K
  \bigl(\Pmat(\qyD)(\widehat v_{\Delta,j}\vD)_y\bigr)_j\dd\bm x\\
  &\quad-
  \left(
  \sum_{e\in\Ehi}\sum_{j=1}^{\kpce}\int_e(\fqu)_j\jump{\widehat u_{\Delta,j}}\dd s
  -\sum_{e\in\Ehi}\sum_{j=1}^{\kpce}\int_e\frac12\bigl(\Pmat(\fq)\jump{\widehat u_{\Delta,j}\uD}\bigr)_j\dd s
  \right)\\
  &\quad-
  \left(
  \sum_{e\in\Ehi}\sum_{j=1}^{\kpce}\int_e(\fqv)_j\jump{\widehat v_{\Delta,j}}\dd s
  -\sum_{e\in\Ehi}\sum_{j=1}^{\kpce}\int_e\frac12\bigl(\Pmat(\fq)\jump{\widehat v_{\Delta,j}\vD}\bigr)_j\dd s
  \right).
\end{aligned}
\end{equation}
The volume terms in $R_1$ cancel pairwise by the product rule and
Lemma~\ref{lem:comm-dot}; this is the discrete counterpart of the
cancellation produced by the skew-symmetric terms in
\eqref{eq:skew-parameterized}. For instance, summing over $j$,
\begin{align*}
  &\bigl(\Pmat(\qxD)\uD\bigr)\cdot(\uD)_x
  -\frac12\sum_{j=1}^{\kpce}\bigl(\Pmat(\qxD)(\widehat u_{\Delta,j}\uD)_x\bigr)_j \\
  &\qquad
  =\Pmat(\qxD)\uD\cdot(\uD)_x
  -\frac12\Pmat(\qxD)\uD\cdot(\uD)_x
  -\frac12\Pmat(\qxD)(\uD)_x\cdot\uD=0.
\end{align*}
The same cancellation holds for the $y$-derivative of $\uD$ and for both
derivative terms involving $\vD$. We emphasize that the cancellation takes
place separately for each of the four combinations of flux direction and
velocity component, so that the asymmetry
\eqref{eq:mixed-flux-representations} of the mixed fluxes, which is part of
the hyperbolicity-preserving construction, does not interfere with
it. Only the face terms remain. On an interior face,
Lemma~\ref{lem:P-jump-average} gives
\begin{align*}
  \frac12\sum_{j=1}^{\kpce}\bigl(\Pmat(\fq)\jump{\widehat u_{\Delta,j}\uD}\bigr)_j
  &=\frac12\bigl(\Pmat(\fq)\uD^+\cdot\uD^+
    -\Pmat(\fq)\uD^-\cdot\uD^-\bigr)
  =\Pmat(\fq)\avg{\uD}\cdot\jump{\uD},\\
  \frac12\sum_{j=1}^{\kpce}\bigl(\Pmat(\fq)\jump{\widehat v_{\Delta,j}\vD}\bigr)_j
  &=\Pmat(\fq)\avg{\vD}\cdot\jump{\vD}.
\end{align*}
Therefore
\begin{equation}\label{eq:R1-interface}
  R_1
  =-\sum_{e\in\Ehi}\int_e
  \bigl(\fqu-\Pmat(\fq)\avg{\uD}\bigr)\cdot\jump{\uD}\dd s
  -\sum_{e\in\Ehi}\int_e
  \bigl(\fqv-\Pmat(\fq)\avg{\vD}\bigr)\cdot\jump{\vD}\dd s.
\end{equation}
Substituting the fluxes \eqref{eq:flux-h}, \eqref{eq:flux-qu}, and \eqref{eq:flux-qv} into \eqref{eq:R1-interface}, and using
\begin{equation}\label{eq:jump-PH-u-identity}
  \jump{\Pmat(\hD+\BD)\uD}-\Pmat(\jump{\hD+\BD})\avg{\uD}
  =\Pmat(\avg{\hD+\BD})\jump{\uD},
\end{equation}
with the analogous identity for $\vD$, yields
\begin{align}\label{eq:R1-diss}
  R_1
  &=-\frac12\sum_{e\in\Ehi}\int_e\alpha_1\,
  \Pmat\bigl(\avg{\hD+\BD}\bigr)\jump{\uD}\cdot\jump{\uD}\dd s\notag\\
  &\quad-
  \frac12\sum_{e\in\Ehi}\int_e\alpha_1\,
  \Pmat\bigl(\avg{\hD+\BD}\bigr)\jump{\vD}\cdot\jump{\vD}\dd s\le 0.
\end{align}
The inequality uses that the symmetric matrix $\Pmat(\avg{\hD+\BD})$ is
positive semidefinite on every face, which is the hypothesis of the theorem
in the entropy stable case. If $\lambda_{\min}\ge0$ denotes the minimum over all faces of
the smallest eigenvalue of $\Pmat(\avg{\hD+\BD})$, then
$R_1\le-\frac12\lambda_{\min}\sum_{e}\int_e\alpha_1(\norml{\jump{\uD}}^2+\norml{\jump{\vD}}^2)\dd s$.
In particular $R_1=0$ when $\alpha_1\equiv0$.

\smallskip\noindent\textbf{Gravitational spatial terms.}
The terms of $R_2$ are the pressure and source terms and the gravity face
corrections:
\begin{equation}
\begin{aligned}
  R_2
  &=\sum_{K\in\Th}\sum_{j=1}^{\kpce}\int_K\widehat q^{x}_{\Delta,j}\bigl(g(\widehat h_{\Delta,j}+\widehat B_{\Delta,j})\bigr)_x\dd\bm x
  +\sum_{K\in\Th}\sum_{j=1}^{\kpce}\int_K\widehat q^{y}_{\Delta,j}\bigl(g(\widehat h_{\Delta,j}+\widehat B_{\Delta,j})\bigr)_y\dd\bm x\\
  &\quad-
  \sum_{K\in\Th}\sum_{j=1}^{\kpce}\int_K
  \bigl(g\Pmat((\hD+\BD)_x)\widehat u_{\Delta,j}\hD\bigr)_j\dd\bm x\\
  &\quad-\sum_{K\in\Th}\sum_{j=1}^{\kpce}\int_K
  \bigl(g\Pmat((\hD+\BD)_y)\widehat v_{\Delta,j}\hD\bigr)_j\dd\bm x\\
  &\quad+
  \sum_{e\in\Ehi}\sum_{j=1}^{\kpce}\int_e
  g\bigl(\Pmat(\jump{\hD+\BD})\avg{\widehat u_{\Delta,j}\hD}\bigr)_j n_x\dd s
  \\
  &\quad+\sum_{e\in\Ehi}\sum_{j=1}^{\kpce}\int_e
  g\bigl(\Pmat(\jump{\hD+\BD})\avg{\widehat v_{\Delta,j}\hD}\bigr)_j n_y\dd s -
  \sum_{e\in\Ehi}\int_eg\fq\cdot\jump{\hD+\BD}\dd s.
\end{aligned}
\end{equation}
The reduction of $R_2$ uses the velocity projection and the commutativity
of $\Pmat$. In the volume terms,
\begin{align*}
  \Pmat((\hD + \BD)_x)\hD\cdot\uD
  &=\Pmat(\hD)(\hD+\BD)_x\cdot\uD
    &&\text{by } \Pmat(\widehat a)\widehat b=\Pmat(\widehat b)\widehat a,\\
  &=\Pmat(\hD)\uD\cdot(\hD + \BD)_x
    &&\text{by Lemma~\ref{lem:comm-dot}},\\
  &=\qxD\cdot(\hD+\BD)_x
    &&\text{by } \eqref{eq:velocity-projection},
\end{align*}
so that the two $x$-volume gravitational terms cancel; the $y$-volume terms
cancel in the same way. Note that the velocity projection
\eqref{eq:velocity-projection} may be applied in the last step because
$(\hD+\BD)_x\in[\mathbb V_\Delta^{k-1}]^{\kpce}\subset[\Vk]^{\kpce}$ is an
admissible test function on every element. For the interface terms, the
same commutativity property gives
\begin{equation}
\begin{aligned}
  \sum_{j=1}^{\kpce}\bigl(\Pmat(\jump{\hD+\BD})\avg{\widehat u_{\Delta,j}\hD}\bigr)_j
  &=\frac12\bigl(\Pmat(\jump{\hD + \BD})\hD^+\cdot\uD^+
  +\Pmat(\jump{\hD+\BD})\hD^-\cdot\uD^-\bigr)\\
  &=\frac12\bigl(\Pmat(\hD^+)\uD^+
  +\Pmat(\hD^-)\uD^-\bigr)\cdot\jump{\hD+\BD}\\
  &=\avg{\Pmat(\hD)\uD}\cdot\jump{\hD+\BD},
\end{aligned}
\end{equation}
and analogously
\[
  \sum_{j=1}^{\kpce}\bigl(\Pmat(\jump{\hD+\BD})\avg{\widehat v_{\Delta,j}\hD}\bigr)_j
  =\avg{\Pmat(\hD)\vD}\cdot\jump{\hD+\BD}.
\]
Therefore $R_2$ reduces to
\begin{align}\label{eq:R2-before-flux}
  R_2
  &=\sum_{e\in\Ehi}\int_e
  \left(
    \avg{g\Pmat(\hD)\uD}n_x
    +\avg{g\Pmat(\hD)\vD}n_y
    -g\fq
  \right)\cdot\jump{\hD+\BD}\dd s .
\end{align}
Substituting the height flux \eqref{eq:flux-h} into \eqref{eq:R2-before-flux} gives
\begin{equation}\label{eq:R2-diss}
  R_2=-\frac{g}{2}\sum_{e\in\Ehi}\int_e\alpha_1
  \norml{\jump{\hD+\BD}}^2\dd s\le 0.
\end{equation}
In particular $R_2=0$ when $\alpha_1\equiv0$.

\smallskip\noindent\textbf{Conclusion of the entropy estimate.}
Combining \eqref{eq:R1-diss} and \eqref{eq:R2-diss},
\[
  \frac{\dd}{\dd t}E_\Delta=R_1+R_2\le 0,
  \qquad \alpha_1\ge 0,
\]
so the scheme is entropy stable for the Lax--Friedrichs-type fluxes. In
the central-flux case $\alpha_1\equiv0$ both $R_1$ and $R_2$ vanish,
\[
  \frac{\dd}{\dd t}E_\Delta=0,
\]
and the scheme is entropy conservative.
\end{proof}

\bibliographystyle{plain}
\bibliography{ref}

\end{document}